\documentclass{article}

\usepackage{amsmath}
\usepackage{amssymb}
\usepackage{amsbsy}
\usepackage{bbm}
\usepackage{hyperref}
\usepackage{subcaption}
\usepackage{algorithm}
\usepackage{algpseudocode}
\usepackage{lineno}
\usepackage{xcolor}
\usepackage{PRIMEarxiv}
\usepackage[numbers]{natbib}
\usepackage[utf8]{inputenc} 
\usepackage[T1]{fontenc}    
\usepackage{hyperref}       
\usepackage{url}            
\usepackage{booktabs}       
\usepackage{amsfonts}       
\usepackage{nicefrac}       
\usepackage{microtype}      
\usepackage{lipsum}
\usepackage{fancyhdr}       
\usepackage{graphicx}       
\graphicspath{{media/}}     

\newtheorem{assumption}{Assumption}
\newtheorem{theorem}{Theorem}
\newtheorem{lemma}{Lemma}
\newtheorem{corollary}{Corollary}
\title{ Optimal Computing Budget Allocation  for Data-driven Ranking and Selection
}

\author{
  Yuhao Wang \\
  School of Industrial and Systems Engineering \\
  Georgia Institute of Technology\\
  Atlanta\\
  \texttt{yuhaowang@gatech.edu} \\
   \And
  Enlu Zhou \\
  School of Industrial and Systems Engineering \\
  Georgia Institute of Technology\\
  Atlanta\\
  \texttt{enlu.zhou@isye.gatech.edu} \\
}

\begin{document}
\maketitle

\begin{abstract}
In a fixed budget ranking and Selection (R\&S) problem, one aims to identify the best design among a finite number of candidates by efficiently allocating the given computing budget to evaluate design performance. Classical methods for R\&S usually assume the distribution of the randomness in the system is exactly known.  
In this paper, we consider the practical scenario where the true distribution is unknown but can be estimated from streaming input data that arrive in batches over time.  We formulate the R\&S problem in this dynamic setting as a multi-stage problem, where we adopt the Bayesian approach to estimate the distribution and formulate a stage-wise optimization problem to allocate the computing budget.
We characterize the optimality conditions for the stage-wise problem by applying the large deviations theory to maximize the decay rate of probability of false selection. Based on the optimality conditions and combined with the updating of distribution estimates, we design two sequential budget allocation procedures for R\&S under streaming input data. We theoretically guarantee the consistency and asymptotic optimality of the proposed procedures. We demonstrate the practical efficiency through numerical experiments in comparison with the equal allocation policy and an extension of the  Optimal Computing Budget Allocation algorithm.
\end{abstract}


\section{Introduction}
\label{sec:intro}

In many applications, performance of a complex stochastic system is often evaluated through time-consuming and expensive experimentation or simulation. Comparison of multiple system designs and selection of the best one is referred to as the ranking and selection (R\&S) problem. R\&S dates back to the 1950s in agricultural and clinical applications (\cite{bechhofer1954single,gupta1956decision}) and has since been used in many application problems arising from a variety of areas such as healthcare and manufacturing. 

While R\&S has traditionally assumed access to a fixed simulator, the recent application of \textbf{ Digital Twin} has been calling for data-driven techniques for R\&S. As indicated by the US National Academies report \cite{national2023foundational}, a digital twin is “a set of virtual information constructs that mimics the structure, context, and behavior of a natural, engineered, or social system (or system of systems), is dynamically updated with data from its physical twin, has a predictive capability, and informs decisions that realize value.” The  physical system continually provides the digital twin with real-time data for better representation of the real system, whereas the digital twin informs predictions and decisions (through simulation optimization such as R\&S) for the real system. This repeated process poses a significant challenge to conventional R\&S with a fixed simulation model due to its fast recurrence, requiring efficient data-driven R\&S that accommodates dynamic updates of the simulation model.
 
To define data-driven R\&S more specifically, consider that the random factors in the system are modeled by some probability distributions, which are usually called ``input distributions'' or ``input models''. 
The input distributions are often estimated with input data observed from the real system. While in many existing works of R\&S, the input distributions are only estimated once with a fixed set of input data prior to the performance evaluation, in many application problems (such as digital twin applications), input data are often collected in batches over time. It is then beneficial to incorporate the sequentially arrived input data to get better input distributions that are more accurate estimates of the underlying distributions of the real system, and subsequently refine the performance evaluations that are used to compare potential system designs. We illustrate the problem of data-driven R\&S with the following  concrete examples.
\begin{enumerate} 
     \item {\bf Selecting the best budget allocation scheme for an activity network:} Activity networks are widely used for project management, computer system security, and digital circuit design (\cite{yang2022optimal,wan2023crashing}). An activity network is represented by a directed acyclic graph with one starting node and one
ending node. The directions of arcs in an activity network represent the precedence relationships between different activities, which all have certain duration to complete. However, over
the course of the project (i.e., to complete all activities), a random disruption (such as natural disasters, electrical outages, workers strikes) may occur, which makes the project completion time (PCT) stochastic. The goal is to select the best budget allocation scheme, which allocates budget to activities to reduce their completion times, in order to minimize the expected PCT. Because of the complexity of large-scale activity networks, evaluating the PCT of a given scenario is usually computationally expensive. Multiple scenarios are generated from the distribution of random disruption, and PCTs are evaluated in these scenarios to form a estimate for the expected PCT. The distribution of the random disruption can be updated when more data become available, and therefore, more evaluations of PCTs can be done under the updated distributions to refine the selection of the best budget allocation scheme. \label{example: activity network}

\item {\bf Selecting the best COVID policy:} The goal is to compare several COVID policies to select the best one in order to minimize the expected number of deaths over a certain time period. Simulation of the spread of an infectious disease, such as COVID, is often based on the so-called compartment models that are widely used in epidemiology \citep{brauer2008compartmental}. It is essentially a system of ordinary differential equations (ODEs) or partial differential equations (PDEs) that might have random initial conditions (e.g., initial infected population size) or unknown parameters (e.g., transmission rate, recovery rate). These unknown factors can be estimated with data that arrive frequently (e.g., daily test data) through methods such as Bayesian inference (\cite{chen2021bayesian}.  Evaluation of each COVID policy requires simulating the compartment model on multiple scenarios of the unknown factors (i.e., samples from the distribution on the initial infected population size, or from the estimated distribution of the transmission rate). For a given scenario, simulation of the compartment model is equivalent to numerically solving a system of ODEs or PDEs till the end of the time period, and hence is computationally expensive. Moreover, with new data collected over time, evaluations of the policies should be updated by running more simulations under the updated distributions.
  \label{example: covid}
\item \textbf{Selecting the best investment strategy:} The goal is to select the best asset portfolio to optimize the trade-off between expected return and risk. Evaluating a portfolio often requires simulating a stochastic model of the asset prices. However, some parameters of the asset price model, such drift and volatility, are unknown and can be estimated from the market data such as the past asset price or interest rates. Given periodically available data, the stochastic model is updated and new simulations are run to update the evaluations of different investment strategies.
\end{enumerate}

These examples above motivate us to consider a fixed budget ranking and selection (R\&S) problem with streaming data where new data arrive over time in batches of possibly varying sizes. The computing budget, which is the amount of total computation that can be used for system performance evaluations, between arrivals of two successive batches of input data is usually limited and determined by external factors such as the inter-arrival time of the data and the computational cost of each performance evaluation. {\it The goal of this paper is to identify the best  allocation scheme of the limited given computing budget to evaluate candidate designs in the setting of streaming input data. }

A procedure for fixed budget R\&S aims to achieve a probability of correct selection (PCS) as high as possible with a given computing budget. Our proposed procedure for the streaming data set is built on the Optimal Computing Budget Allocation (OCBA) algorithm, which is one of the most widely applied and studied algorithms for fixed budget R\&S with a fixed input distribution. It computes the budget allocation rule by repeatedly maximizing an approximate PCS objective in each iteration with plug-in estimators of design performances. OCBA was first proposed in \cite{chen2000simulation} and was shown to converge asymptotically to the optimal allocation rule in \cite{Li2021on}. The statistical validity of OCBA crucially relies on the stationarity of the underlying input distribution, which implies that the performance evaluations are independent and identically distributed (i.i.d.) samples even though they are generated in different iterations of the algorithm. Hence, the performance estimation error  diminishes as more samples are generated over iterations, leading to convergence of the allocation policy to the optimal policy. OCBA have been extended in the past years to various problems such as subset selection (\cite{chen2008efficient} and \cite{gao2015efficient}),  contextual R\&S (\cite{gao2019selecting} and \cite{jin2019optimal}), multi-objective (\cite{lee2004optimal}), finding simplest good designs (\cite{jia2013efficient} and \cite{yan2012efficient}), maximizing opportunity cost (\cite{gao2017new}), robust R\&S under input uncertainty (\cite{gao2017robust}), and many others. We refer the reader to \cite{chen2011stochastic} for a comprehensive tutorial on OCBA.  All of these works either assume known input distributions or consider an empirical input distribution estimated with a fixed set of input data, where the performance evaluations for the same design follow the same input distribution over iterations. 

The setting of streaming data considered in this paper is more challenging than the setting of fixed input distribution in the aforementioned OCBA works. Most notably, the input distribution is no longer fixed but is updated with new input data at each time stage, and hence it breaks the i.i.d. condition of the performance evaluations for each fixed design. 
To address this challenge,  we adopt a Bayesian approach to estimate the unknown parametrized input distribution. We begin with a finite parameter space, 
where each input parameter represents a simulation scenario, so that we can evaluate the system performance under a fixed input parameter (which will be referred to as ``design-input pair''), and aggregate these evaluations according to the current posterior probability for each input parameter. Since the simulation is conducted under fixed input parameter,  we are then able to generate i.i.d. performance evaluations (samples) under the same design-input pair. As a result, we can estimate the expected performance for a fixed design through a Bayesian average estimator, which is the weighted sum of the estimated design-input performance multiplied by its corresponding posterior probability. We derive the optimal budget allocation policy among different design-input pairs and propose two fully sequential procedures with provable performance guarantees.
Then, we extend the proposed procedures to the more general setting where the input parameter space can be continuous, by partitioning the entire parameter space into a finite set of subspace.  Despite the challenge of non-stationary simulation outputs posed by continuous input parameter subspace as opposed to fixed input parameter, we show the statistical validity of the proposed algorithms holds with the same budget allocation rules given by the proposed algorithms. 

We summarize the contributions of this paper as follows.

1. This paper, along with our earlier conference version \cite{wang2022fixed}, is the first to consider streaming input data in fixed budget R\&S problems and design a data-driven approach. This paper differs from the conference version in three main ways. First, we relax the assumption of a discrete input distribution with finite support, and instead take a Bayesian approach to estimate the parametrized input distribution, which can have continuous support. Second, we develop a new procedure 
that is designed to solve the original reformulated problem rather than compute an approximate solution. Third, we obtain a stronger convergence result that characterizes the speed at which the allocation policies given by the procedures converge to the optimal solution.

2. We propose a new framework for conducting simulation sequentially with the presence of streaming input data, where we run the simulation for a fixed design-input pair each time and adopt a  Bayesian average estimator, defined as the weighted sum of design-input sample means weighted by the corresponding posterior probability, to estimate the expected performance of each design. 
To calculate the rate function of this performance estimator, we apply the Gartner-Ellis theorem (see \cite{dembo1994large}). Unlike many other works derived from \cite{glynn2004large}, where they directly use the same rate function, we need to recalculate it due to samples from different distributions in the performance estimator. We formulate a stage-wise rate maximization problem using the recalculated rate function and derive the corresponding optimality conditions. Compared with the optimality conditions in \cite{glynn2004large}, we obtain an additional ``Input Balance" 
condition that characterizes the allocation rule among different input parameters according to the current posterior distribution.

3.  We develop two fully sequential procedures of data-driven OCBA (DD-OCBA), namely DD-OCBA-approx(-C) and DD-OCBA-balance(-C), based on different approaches to solve the optimality conditions with either finite or continuous input parameter space. We prove the statistical consistency and asymptotic optimality of both procedures. Specifically, for the setting of finite parameter space, we provide a stronger convergence result by characterizing the asymptotic convergence rate of the allocation policies. 
 

Next, we briefly review the relevant literature with an emphasis on the relation to our work.

\subsection*{Literature Review}
The research on R\&S largely falls into two related yet different categories. The fixed confidence R\&S procedures aim to achieve a pre-specified probability of correct selection (PCS) using the least possible amount of simulation effort, whereas the fixed budget R\&S procedures typically tend to attain a PCS as high as possible with a given simulation budget. For fixed confidence, a large body of literature goes to the indifference zone (IZ) formulation. An IZ procedure guarantees selecting the best design with at least a pre-specified confidence level, given that the difference between the top-two designs is sufficiently large. Existing IZ procedures in the R\&S literature include but are not limited to the KN procedure in \cite{kim2001fully}, the  KVP and UVP procedures in \cite{jeff2006fully}, and the BIZ procedure in \cite{frazier2014fully}. We refer the reader to \cite{kim2007recent} for a comprehensive review of IZ formulations. In addition, the Bayesian approaches in \cite{chick2001input,chick2012sequential} and the probably approximately correct (PAC) selection in \cite{ma2017efficient} has also been studied in this stream of works. 

In this paper, we focused on the fixed budget R\&S. As discussed in Section \ref{sec:intro}, OCBA was originally derived under a normality assumption and an approximate PCS objective.  The allocation rule of OCBA can be justified from a rigorous perspective of the large deviations theory in \cite{glynn2004large}, after which lots of works followed this large deviations formulation. For instance, \cite{chen2019complete} designed a fully sequential budget allocation algorithm for general distributions using the optimality conditions in \cite{glynn2004large};  \cite{hunter2013optimal} and \cite{pasupathy2014stochastically} applied the large deviations theory to constrained R\&S; \cite{gao2019selecting} and \cite{cakmak2022contextual} extended the large deviations approach to contextual R\&S;  \cite{he2007opportunity} applied OCBA procedure to optimize the opportunity cost as opposed to PCS; \cite{gao2017robust} computed the large deviations rate (LDR) function with respect to a worst-case performance estimator. In this paper, we compute the LDR function of a performance estimator aggregating samples across different input distributions. Other well-known fixed budget R\&S procedures include the expected value of information (EVI) approach proposed by \cite{chick2010sequential} and the knowledge-gradient (KG) approach proposed by \cite{frazier2009knowledge}, where EVI is derived by asymptotically minimizing a bound of the expected loss and KG determines the optimal sampling allocation policy by maximizing the so-called acquisition function. We refer the reader to  \cite{hong2021review} for a recent overview of the R\&S literature.


All the aforementioned works assume the underlying distribution is known. While there are extensive studies on the impact of estimated input distributions on the simulation outputs (e.g., \citep{chick2001input,zouaoui2004accounting,ng2006reducing,xie2014bayesian,lam2017empirical}), which we refer the reader to \cite{corlu2020stochastic} for a recent review, R\&S with estimated input distributions have only been studied in recent years.  \cite{corlu2013subset,corlu2015subset} aimed to eliminate as many inferior designs as possible and return a subset of superior designs with a fixed amount of input data. \cite{wu2017ranking} formulated a fixed budget problem under OCBA framework to simultaneously allocate the effort to carry out stimulation and the effort to obtain input data; this work is followed by \cite{xu2020joint}, which proposed a general framework that integrates input data collection and simulation in which the data collection and simulation costs themselves can be random.   \cite{gao2017robust}, \cite{xiao2018simulation}, \cite{xiao2020optimal} took a fixed budget formulation with a robust approach, aiming to select a design with the best worst-case performance over an uncertainty set of finite distributions that contains the true input distribution; \cite{fan2020distributionally} also used this worst-case criterion but took an IZ formulation.   \cite{song2019input} derived confidence bands to account for both estimation error in input distribution estimation and stochastic error in simulation output in R\&S. 
Despite the assumption of estimated input distributions, these works focused on a fixed set of input data. 

More recently, \cite{Wu2022data} considered R\&S with streaming input data, similar to the setting in this paper, but used a fixed confidence formulation. They proposed a moving average performance estimator to aggregate simulation outputs from different input distributions over time stages and designed sequential elimination procedures to screen out the inferior designs until one is left to be the optimal one with at least a pre-specified confidence level. The methodology in \cite{Wu2022data} is fundamentally different from our approach, as their focus is on deriving the valid concentration bound to reach the given confidence level while we focus on optimizing budget allocation. In \cite{Wu2022data}, each remaining design is simulated once at every iteration, which allows the usage of common random numbers to reduce the variance of simulation output and narrow the concentration bound. However, from the perspective of budget allocation, \cite{Wu2022data} trivially allocates the simulation budget equally to all designs, which can be inefficient if the simulation budget itself is limited.
 The fixed budget R\&S with streaming data is first considered in \cite{wang2022fixed} and \cite{Kim2022optimizing}, where \cite{wang2022fixed} is the early conference version of this paper as discussed in the previous section. \cite{Kim2022optimizing} considers the setting where input data can be actively collected, which is similar to the setting of \cite{wu2017ranking} but with periodic data collection and a selecting criterion called ``most probable best". 

The rest of the paper is organized as follows. We describe the problem setting and present the overall framework of the proposed data-driven procedures in Section \ref{sec: problem formulation}. In Section \ref{sec:LDR} we explicitly solve the stage-wise budget allocation problem by applying the large deviations theory to calculate the rate function for the performance estimator and characterize the stage-wise optimal allocation policy. In Section \ref{sec:procedure} we propose sequential procedures for the R\&S problem under streaming data. We show the statistical consistency and asymptotic optimality of the two procedures in Section \ref{sec:consistency}. In Section \ref{sec: continuous space}, we extend the framework to continuous input parameter space, while preserving the provable statistical validity. We present numerical results in Section \ref{sec:numerical} and conclude in Section \ref{sec:conclusion}. 

\section{Problem Statement} \label{sec: problem formulation}
We first give some basic notations. Suppose we have a set of finite number of designs $\mathcal{I} = \{1,2,\cdots,K\}$, and the goal is to find the design with the highest expected performance. The performance of each design $i \in \mathcal{I}$ is evaluated through repeated simulations.
The computing budget, which is the total number of replications we can run on all designs, is often limited by computational time or expense. 
The core of the fixed budget R\&S problem is to devise procedures that maximize the probability of correct selection (PCS) of the optimal design  when exhausting the computing budget.

In classical R\&S, the input distributions, $\{F_i\}_{i\in\mathcal{I}}$, that capture various sources of system randomness are assumed to be known. However, in practice the true underlying distributions are seldom known and need to be estimated from input data, which are a finite amount of real-world observations.  Sources of system randomness, such as disruption delay in the activity network, are often shared among all designs. Therefore, throughout the paper we assume that all the designs share the same input distribution $F^c$ (thus, dropping the subscript $i$) and consequently common input data from these distributions. Here the superscript $c$ stands for ``correct", meaning that the input distribution is exactly the same as the true distribution. Note that here the assumption of common input distribution does not rule out the existence of possibly design-specific distributions, since we can simply incorporate them into the common input distribution even though they may not affect the simulation of all designs.  


We assume the underlying input distributions belong to some parametric family.
\begin{assumption} \label{assump: parametric input}
\begin{enumerate}
    \item  The cumulative density function (cdf) of the true input distribution, $F^c = F_{\theta^c}$, belongs to a known parametric family $\{F_\theta: \theta \in \Theta \}$. 
    \item  Furthermore, the parameter space $\Theta = \{\theta_1,\ldots,\theta_D\}$ is finite. 
\end{enumerate}
  
\end{assumption}
The parametric family in Assumption \ref{assump: parametric input} enables us to estimate the input distribution $F^c$ through estimating the true parameter $\theta^c$. The finite assumption on parameter space can be regarded as approximation of the original space, which can be obtained by either discretization or sampling from the prior distribution. While for now we focus on the finite parameter space, we will show in Section \ref{sec: continuous space} our method can be extended to a general continuous parameter space.

We take a Bayesian approach to estimate the unknown parameter $\theta^c$. Specifically, let $\pi_0$ denote the prior distribution and $f_\theta$ denote the density (likelihood) function of $F_{\theta}$. We make the following assumption on the prior distribution and the likelihood function, which are standard assumptions to guarantee the strong consistency of the Bayesian posterior distribution.  

\begin{assumption} \label{assump: prior likelihood}
    \begin{enumerate}
        \item The prior distribution $\pi^0$ satisfies $\pi^0(\theta) > 0, \theta\in\Theta$. 
        \item The likelihood function $f_\theta$ satisfies
        $f_\theta(\xi) > 0$ almost surely for all $\theta \in \Theta$, where $\xi \sim F_{\theta^c}$.
    \end{enumerate}
\end{assumption}
Given $m$ independent and identically distributed (i.i.d.) data $\{\xi_1,\ldots,\xi_m \}$ with $\xi_i \sim F_{\theta^c}$, the posterior distribution $\pi = [\pi_1,\ldots,\pi_D]$ 
is computed as  
$$ \pi_j := \mathbb{P}\left(\theta = \theta_j | \xi_1,\ldots,\xi_m\right) = \frac{\pi^0(\theta_j) \prod_{\ell=1}^m f_{\theta_j}(\xi_\ell)}{\sum_{j'=1}^D \pi^0(\theta_{j'}) \prod_{\ell=1}^m f_{\theta_{j'}}(\xi_\ell) }, j=1,\ldots,D.$$
By Assumption~\ref{assump: prior likelihood}, $\pi_j > 0$ almost surely.

The posterior distribution provides a density estimate of the unknown parameter $\theta$, and naturally leads to an estimate for the expected performance. Specifically, let $X_i(\theta)$ denote the random performance of design $i$ under input parameter $\theta$ and $\mu_i(\theta) = \mathbb{E}[X_i(\theta)|\theta]$ the expected performance of design $i$ under $\theta$. Moreover, with finite parameter space $\Theta$, we can write $X_{ij} := X_i(\theta_j)$ and $\mu_{ij} := \mathbb{E}[X_{ij}]$. The expected performance under the unknown true input distribution, $\mu_i(\theta^c)$, can be estimated by the Bayesian average performance, which is expected performance under the posterior:
\begin{equation}\label{BayesianAvg}
    \Bar{\mu}_{i}  = \mathbb{E}_\pi \left[ \mu_i(\theta)\right] = \sum_{j=1}^D \pi_j \mu_{ij}.
\end{equation}


We make the following assumption on the Bayesian average performance as well as the true expected performance, to ensure the uniqueness of the (estimated and true) best design.

\begin{assumption} \label{assump: unique b}
    Let $b := \arg\max_i \Bar{\mu}_i$ denote the design that maximizes the Bayesian average performance,  which is unique for almost every $\pi$.  Moreover, the true best design $b^c := \arg\max_i \mu_i(\theta^c)$ is also unique.
\end{assumption}

By \eqref{BayesianAvg}, Bayesian average performance $\Bar{\mu}_{i}$ allows aggregation of $\mu_{ij}$'s under different input distributions $F_{\theta_j}$'s by simply taking the weighted average with weights $\pi_j$'s.
We will refer to simulating design $i$ under input distribution $F_{\theta_j}$ as simulating the $(i,j)$ ``design-input'' pair 
throughout the paper. A great advantage of simulating under fixed design-input pair is that the simulation outputs for the fixed design-input pair are i.i.d. across different time stages, despite the varying posterior distribution.

For each design $i$, let $X_{i,j}^{\ell}$ denote the $\ell^{th}$ sample of $X_{i,j}$. 
We make the following assumption of Gaussian simulation output.

 \begin{assumption}\label{assump:output}\  
 \begin{enumerate}
     \item The simulation output $X_{i,j}$ follows a normal distribution with unknown mean $\mu_{i,j}$ and unknown variance $\sigma_{i,j}^2$.  
     \item The simulation output  $\{X_{i,j}^{(\ell)}\}$ are independent for all $i,j$ and $\ell$.
 \end{enumerate}
 \end{assumption}

Assumption \ref{assump:output}.1 models the simulation error as Gaussian noise, which is common in the R\&S literature as one can use batched simulation output. 
Assumption \ref{assump:output}.2 can be guaranteed since now we simulate on fixed distribution. 


\subsection{Data-driven Framework}
\begin{figure}
     \centering
     \includegraphics[width=0.8\textwidth]{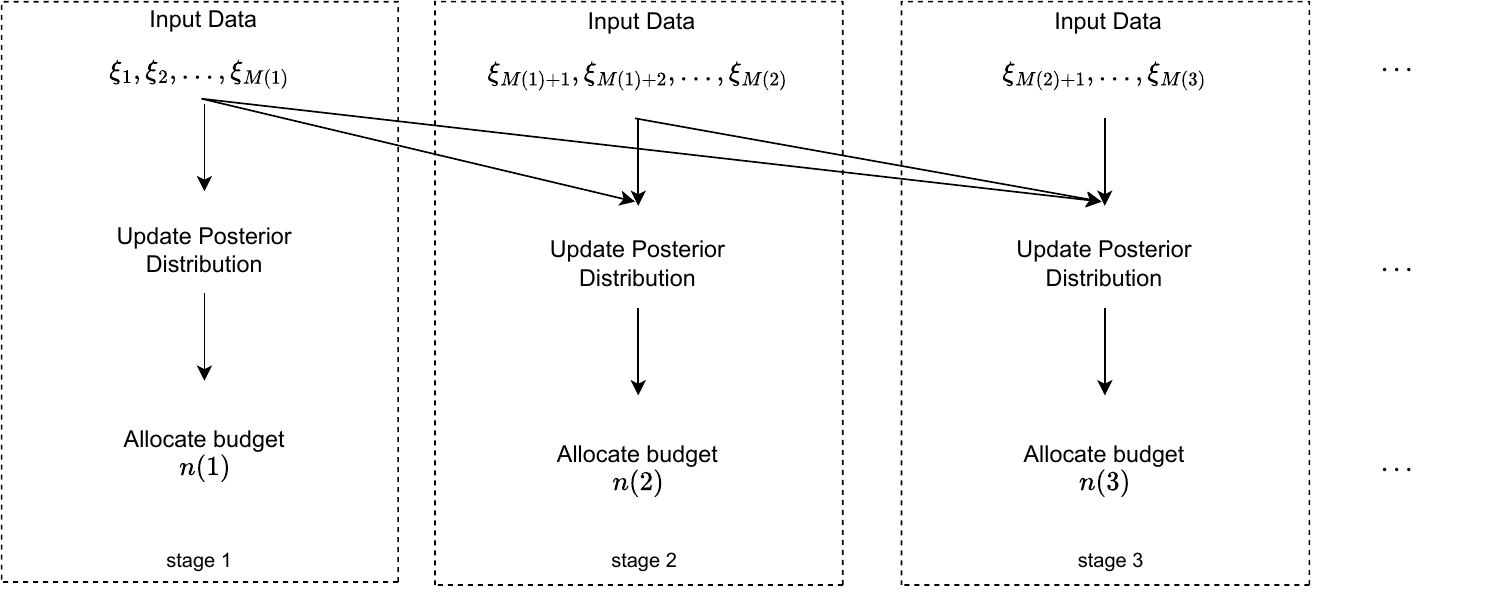}
     \caption{Illustration of budget allocation with streaming input data}
     \label{fig:procedure}
 \end{figure}

Next, we describe the overall framework of our proposed data-driven budget allocation procedures, where the posterior distribution $\pi$ is updated periodically with streaming data.
Specifically, at time stage $t$, new input data of batch size $m(t)$ can be obtained and used to update the estimate of the input distribution, and then we allocate computing budget $n(t)$ to design-input pairs according to the current estimated input distribution.  We assume both $n(t)$ and $m(t)$ are given. This process is illustrated in Figure \ref{fig:procedure}, where $M(t) = \sum_{\tau=1}^t m(\tau)$ is the total amount of input data collected up to stage $t$. 
For input data, we make the following assumption on the input data to help guarantee the Bayesian consistency of the posterior distribution $\pi$.
\begin{assumption} \label{assump:input data}
 The input data, $\{\xi_s\}_{s=1}^\infty$, are identically and independently distributed.
\end{assumption}
 To find the budget allocation rule for each stage, we apply the large deviations theory to formulate an optimization problem under the current estimated input distribution and characterize its optimality condition to derive the stage-wise optimal budget allocation rule in Section \ref{sec:LDR}. Then combining with the updating of posterior distribution, we develop two data-driven budget allocation procedures for the multi-stage setting in Section \ref{sec:procedure}.

   \section{Rate-Optimal Budget Allocation}\label{sec:LDR}
In this section, we formulate and solve a static optimal budget allocation problem 
under the current posterior distribution $\pi$. Let $n$ denote the total simulation budget and $n_{i,j}$ denote the simulation budget allocated to design $i$ under input parameter $j$ (i.e., $\theta_j$).  Let $\alpha_i = (\alpha_{i,1},\ldots,\alpha_{i,D})^\intercal$ be the ratio of budget allocated to design $i$ and across all input parameters, $\alpha_{i,j} = \frac{n_{i,j}}{n}$. Let $$ \hat{\mu}_{i,j}(\alpha_{i,j}n) = \frac{1}{\alpha_{i,j}n}\sum_{s=1}^{\alpha_{i,j}n} X_{i,j}^{(s)} \qquad \text{and} \qquad \hat{{\mu}}_{i}(\alpha_i,n) = \sum_{j=1}^D \pi_j \hat{\mu}_{i,j}(\alpha_{i,j}n),$$ denote the estimated performance for the $(i,j)$  (design-input) pair and the estimated performance for design $i$, respectively. Ignoring the minor issue of $\alpha_{i,j}n$ not being an integer, we formulate an optimization problem from the  large deviations perspective as in \cite{glynn2004large} and define the rate function 
\[ \mathbf{G}_i(\alpha_b,\alpha_i):= -\lim_{n\rightarrow \infty} \frac{1}{n} \log \mathbf{P}\left(\hat{{\mu}}_b(\alpha_b,n) \le  \hat{{\mu}}_i(\alpha_i,n)\right).\]

That is, $n\mathbf{G}_i(\alpha_b,\alpha_i) $ is  the exponential rate  of the probability of the event $\{ \hat{{\mu}}_b(\alpha_b,n) \le  \hat{{\mu}}_i(\alpha_i,n)\}$ going to zero as $n$ goes to infinity.  Since the probability of false selection (PFS), which is defined as 
$$ PFS=\mathbf{P} \left(\hat{{\mu}}_b(\alpha_b,n) \le  \max_{ i \neq b} \hat{{\mu}}_i(\alpha_i,n)\right),$$ can be bounded by
$$\max_{ i \neq b} \mathbf{P}\left(\hat{{\mu}}_b(\alpha_b,n) \le  \hat{{\mu}}_i(\alpha_i,n)\right) 
\le PFS \le
(K-1) \max_{ i \neq b} \mathbf{P}\left(\hat{{\mu}}_b(\alpha_b,n) \le  \hat{{\mu}}_i(\alpha_i,n)\right),
 $$
 we have 
 \begin{equation}\label{eq:LDR}
    \lim_{n\rightarrow \infty} \frac{1}{n} \log PFS = -\min_{ i \neq b}\mathbf{G}_i(\alpha_b,\alpha_i) .
\end{equation}
That is, $\min_{ i \neq b}\mathbf{G}_i(\alpha_b,\alpha_i)$ is the asymptotically exponential decay rate of PFS. To maximize this decay rate of PFS, we consider the following optimization problem:
\begin{equation}\label{eq:opt LDR}
    \begin{aligned}
    \max_{\alpha_{i,j}, 1\le i\le K, 1\le j \le D} &\quad z &\\
   s.t. \qquad &\mathbf{G}_i(\alpha_b,\alpha_i)  - z \ge 0& \quad & i \neq b\\
      &\sum_{i=1}^K \sum_{j=1}^D \alpha_{i,j} = 1 & \\
      &\alpha_{i,j} \ge 0 & \quad &1 \le i \le K, 1\le j\le D.
    \end{aligned}
\end{equation}

Optimizing the large deviations rate of PFS is first studied in \cite{glynn2004large}, where input distribution is assumed to be known. Their formulation can be seen as a special case of $D=1$ in our setting. Our approach is an extension of \cite{glynn2004large}, which incorporates simulation samples under different input parameter in calculating the rate function. Due to this difference, we cannot directly apply their result. Instead, we take a similar approach using the Gartner-Ellis Theorem (see \cite{dembo1994large}), and the detailed calculation is shown in the next section.



\subsection{Calculation of the Rate Function}
In this section we give the explicit form of the rate function $\mathbf{G}_i(\alpha_b,\alpha_i)$. 

For a fixed $i$, let $\Lambda_{i,j}(\cdot)$ denote the log-moment generating function of $X_{i,j}$ and $\Lambda_n(\cdot,\cdot)$ denote the log-moment generating function of $Z_n = (\hat{{\mu}}_b(\alpha_b,n),\hat{{\mu}}_i(\alpha_i,n))$. That is,
\begin{align*}
    \Lambda_{n} (\lambda_b,\lambda_i) =& \log\mathbf{E}\left[\mathrm{e}^{\lambda_b\hat{{\mu}}_b(\alpha_b,n)+\lambda_i\hat{{\mu}}_i(\alpha_i,n)}\right] \\
    =& \log \mathbf{E}\left[ \exp \left(\lambda_b\sum_{j=1}^D {\pi_j}\sum_{s=1}^{n\alpha_{b,j}} \frac{X_{b,j}^{(s)}}{n\alpha_{b,j}} + \lambda_i\sum_{j=1}^D {\pi_j}\sum_{s=1}^{n\alpha_{i,j}} \frac{X_{i,j}^{(s)}}{n\alpha_{i,j}}\right)  \right]\\
    =&\sum_{j=1}^D  n\alpha_{b,j} \log \mathbf{E}\left[\exp\left(\frac{\lambda_b {\pi_j}}{n\alpha_{b,j} }X_{b,j} \right) \right] +\sum_{j=1}^D  n\alpha_{i,j} \log \mathbf{E}\left[\exp\left(\frac{\lambda_i {\pi_j}}{n\alpha_{i,j} }X_{i,j} \right) \right]\\
    =& \sum_{j=1}^D n\alpha_{b,j}\Lambda_{b,j}\left(\frac{\lambda_b{\pi_j}}{n\alpha_{b,j}}\right) +\sum_{j=1}^D n\alpha_{i,j}\Lambda_{i,j}\left(\frac{\lambda_i{\pi_j}}{n\alpha_{i,j}}\right).
\end{align*}
The third equality follows from Assumption \ref{assump:output}.2 that all simulation outputs are independent across designs and simulation outputs are identically distributed for the same design-input pair.
Then, substituting $\lambda_b$ and $\lambda_i$ with $n\lambda_b$ and $n\lambda_i$, respectively, we obtain
$$ \frac{1}{n} \Lambda_n(n\lambda_b,n\lambda_i) = \sum_{j=1}^D \alpha_{b,j}\Lambda_{b,j}\left(\frac{\lambda_b{\pi_j}}{\alpha_{b,j}}\right) +\sum_{j=1}^D \alpha_{i,j}\Lambda_{i,j}\left(\frac{\lambda_i{\pi_j}}{\alpha_{i,j}}\right).$$
Since $X_{i,j}$ follows a normal distribution with mean $\mu_{i,j}$ and variance $\sigma_{i,j}^2$, $\Lambda_{i,j}(\lambda) = \lambda \mu_{i,j} + \frac{1}{2}\lambda^2\sigma_{i,j}^2$.
Let $I(x_b,x_i)$ be the Fenchel-Legendre transform of $\Lambda_n$. Then,
$$
\begin{aligned}
I(x_b,x_i) &= \sup_{\lambda_b,\lambda_i} \left\{ \lambda_bx_b+\lambda_ix_i - \sum_{j=1}^D \alpha_{b,j}\Lambda_{b,j}\left(\frac{\lambda_b{\pi_j}}{\alpha_{b,j}}\right) -\sum_{j=1}^D \alpha_{i,j}\Lambda_{i,j}\left(\frac{\lambda_i{\pi_j}}{\alpha_{i,j}}\right)\right\}\\
&= \sup_{\lambda_b} \left\{ \lambda_bx_b - \sum_{j=1}^D \alpha_{b,j}\Lambda_{b,j}\left(\frac{\lambda_b{\pi_j}}{\alpha_{b,j}}\right) \right\}+\sup_{\lambda_i} \left\{\lambda_ix_i-\sum_{j=1}^D \alpha_{i,j}\Lambda_{i,j}\left(\frac{\lambda_i{\pi_j}}{\alpha_{i,j}}\right)\right\}\\
& = \sup_{\lambda_b} \left\{ \lambda_bx_b - \sum_{j=1}^D \left(\lambda_b{\pi_j}\mu_{1,j} + \frac{1}{2}\frac{\sigma_{b,j}^2\lambda_b^2{\pi_j}^2}{\alpha_{b,j}} \right) \right\}+\sup_{\lambda_i} \left\{ \lambda_ix_i - \sum_{j=1}^D \left(\lambda_i{\pi_j}\mu_{i,j} + \frac{1}{2}\frac{\sigma_{i,j}^2\lambda_i^2{\pi_j}^2}{\alpha_{i,j}} \right) \right\}\\
&=\underbrace{\frac{1}{2} \frac{(x_b - \Bar{\mu}_{b})^2}{\sum_{j=1}^D\frac{\sigma_{b,j}^2{\pi_j}^2}{\alpha_{b,j}}}}_{\textstyle \mathstrut =I_b}+
\underbrace{\frac{1}{2} \frac{(x_i - \Bar{\mu}_{i})^2}{\sum_{j=1}^D\frac{\sigma_{i,j}^2{\pi_j}^2}{\alpha_{i,j}}}}_{\textstyle \mathstrut =I_i}.
\end{aligned}
$$
By the Gartner-Ellis Theorem,
$\mathbf{G}_i(\alpha_b,\alpha_i) = \inf_{x_b \le x_i} I(x_b,x_i). $
It is easy to see that  $I_b$ is decreasing for $x_b \le \Bar{\mu}_{b}$ and increasing for $x_b \ge \Bar{\mu}_{b}$, and $I_i$ is decreasing for $x_i \le \Bar{\mu}_{i}$ and increasing for $x_i \ge \Bar{\mu}_{i}$. Since $\Bar{\mu}_{b} > \Bar{\mu}_{i}$, we must have
\begin{equation} \label{eq:rate formula}
    \mathbf{G}_i(\alpha_b,\alpha_i) = \inf_{\Bar{\mu}_{i} \le x \le \Bar{\mu}_{b}} I(x,x) = \frac{(\Bar{\mu}_{b} - \Bar{\mu}_{i})^2}{2\left( \sum_{j=1}^D \frac{\sigma_{b,j}^2{\pi_j}^2}{\alpha_{b,j}} + \sum_{j=1}^D\frac{\sigma_{i,j}^2{\pi_j}^2}{\alpha_{i,j}}\right)}.
\end{equation}  
When $D=1$, we recover exactly the same rate function as in \cite{glynn2004large}. The following lemma summarizes some important properties of $\mathbf{G}_i(\alpha_b,\alpha_i)$.
\begin{lemma} \label{lem:rate function}\ Suppose Assumption \ref{assump:output} holds. Then,
\begin{enumerate}
    \item $\mathbf{G}_i(\alpha_b,\alpha_i)$ is strictly increasing in $\alpha_{b,j}$ and $\alpha_{i,j}$ for $\alpha_{b,j},\alpha_{i,j}>0,     j=1,2,\ldots,D $. Moreover, $\mathbf{G}_i(\alpha_b,\alpha_i) = 0$ if there exists $j_0$ such that $\min(\alpha_{b,j_0},\alpha_{i,j_0})=0$.
    \item $\mathbf{G}_i(\alpha_b,\alpha_i)$ is concave in $(\alpha_b,\alpha_i)$ for $(\alpha_b,\alpha_i) > 0$.
\end{enumerate}
\end{lemma}

Lemma \ref{lem:rate function}.1  implies that any design-input pair  must be allocated with a positive ratio of the simulation budget; otherwise, the rate will be zero. Lemma \ref{lem:rate function}.2 claims the concavity of $\mathbf{G}_i$, which guarantees the optimality with the Karush–Kuhn–Tucker (KKT) condition for the optimization problem \eqref{eq:opt LDR} in the following section.
\subsection{Optimal Allocation Policy} \label{sec:optimality}
In this section we derive the optimality conditions for problem \eqref{eq:opt LDR}, shown in the following theorem.

\begin{theorem} \label{thm:optimal allocation}
Suppose Assumption \ref{assump:output} holds. Let $\alpha \ge 0$ be a feasible allocation policy. Then $\alpha$ is the optimal solution to (\ref{eq:opt LDR}) if and only the following three conditions hold:
\begin{flalign}
    &1.\text{(Input Balance) }\ \frac{\partial\mathbf{G}_i(\alpha_b,\alpha_i)}{\partial\alpha_{i,j}} = \frac{\partial\mathbf{G}_i(\alpha_b,\alpha_i)}{\partial\alpha_{i,j'}} \quad  i \neq b \text{ and }  1 \le j < j' \le D;  && \label{thm:input balance} \\
    &2.\text{(Total Balance) } \ \sum_{i \neq b} \frac{\partial \mathbf{G}_i(\alpha_b,\alpha_i)/\partial\alpha_{b,j}}{\partial\mathbf{G}_i(\alpha_b,\alpha_i)/\partial\alpha_{i,j}} = 1 \quad  1\le j\le D;&& \label{thm:total balance}\\
    &3.\text{(Local Balance) } \ \mathbf{G}_i(\alpha_b,\alpha_i) = \mathbf{G}_{i'}(\alpha_b,\alpha_{i'}) \quad   i \neq i'\neq b.&& \label{thm:local balance}
\end{flalign}
Or equivalently in the explicit form:
\begin{flalign}
    &1.\ \frac{\alpha_{i,j}}{\sigma_{i,j}{\pi_j}} = \frac{\alpha_{i,j'}}{\sigma_{i,j'}\pi_{j'}} \quad i \neq b, \  1 \le j \le D, && \label{eq:input balance} \\
    &2. \ \left(\frac{\alpha_{b,j}}{\sigma_{b,j}}\right)^2 = \sum\limits_{i\neq b} \left(\frac{\alpha_{i,j}}{\sigma_{i,j}}\right)^2 \quad 1 \le j \le D,&& \label{eq:total balance}\\
    &3. \ \frac{(\Bar{\mu}_{b} - \Bar{\mu}_{i})^2}{\sum_{j=1}^D \frac{\sigma_{b,j}^2{{\pi_j}^2}}{\alpha_{b,j}}+\sum_{j=1}^D\frac{\sigma_{i,j}^2{{\pi_j}^2}}{\alpha_{i,j}}} = \frac{(\Bar{\mu}_{b} - \Bar{\mu}_{i'})^2}{\sum_{j=1}^D \frac{\sigma_{b,j}^2{{\pi_j}^2}}{\alpha_{b,j}}+\sum_{j=1}^D\frac{\sigma_{i',j}^2{{\pi_j}^2}}{\alpha_{i',j}}} \quad i \neq i'\neq b.&& \label{eq:local balance}
\end{flalign}
Furthermore, the optimal solution $\alpha^*$ to (\ref{eq:opt LDR}) is unique.
\end{theorem}

{\it {\bf Remark:} Compared with the optimality condition in \cite{glynn2004large}, in addition to the ``total balance"  condition that characterizes the relation between the optimal design and the non-optimal designs and the ``local balance" conditions that characterize the relation between two non-optimal designs.
Here we have the additional optimality condition (\ref{thm:input balance}), the ``input balance" condition. It states that within the allocation for a certain design $i$, the partial derivative of the rate function $\mathbf{G}_i$ with respect to $\alpha_{i,j}$ is the same for all $j$'s. That is, simulation for each fixed input parameter should provide the same improvement to identify that design $b$ is better than $i$. Furthermore, with normally distributed simulation errors, equation \eqref{eq:input balance} indicates that for a fixed design $i$ the optimal allocation ratio $\alpha_{i,j}$  should be proportional to the posterior probability mass ${\pi_j}$ and the standard deviation $\sigma_{i,j}$, which quantitatively characterizes how input uncertainty affects the optimal allocation policy.}

Also notice for fixed $i$, \eqref{eq:input balance} only depends on $i$, which means the relative allocation ratios among different input parameters for a fixed design do not depend on other designs.  On the other hand, \eqref{eq:local balance} indicates that the relative allocation ratios among designs under the same input parameter $j$ are affected by all ${\pi_j}$'s, which implies directly applying OCBA to designs under a fixed input parameter $j$ may perform poorly since it does not take information from other design-input pairs into consideration.  
Moreover, notice that the three optimality conditions \eqref{thm:input balance}-\eqref{thm:local balance} not only hold for Gaussian simulation noise but also  hold as long as the rate function $G_i$ has the properties shown in Lemma \ref{lem:rate function}.

\section{Sequential Procedure with Streaming Input Data} \label{sec:procedure}
In deriving Theorem \ref{thm:optimal allocation} above, we assume a fixed posterior distribution and full knowledge of simulation output distribution. In this section, by trying to satisfy the optimality conditions in Theorem \ref{thm:optimal allocation} with the current posterior distribution, we develop two data-driven optimal budget computing budget allocation (DD-OCBA) procedures, namely DD-OCBA-approx and DD-OCBA-balance, for simulation budget allocation in the multi-stage setting with streaming input data. The two procedures mainly differ in how to satisfy the optimality conditions: DD-OCBA-approx solves the optimality conditions approximately, while DD-OCBA-balance tries to balance the two sides of the optimality equations. 
In the implementation phase, all unknown parameters, such as the posterior probability \(\pi_j\), the Bayesian average performance \(\bar{\mu}_i\), and the design-input simulation variance \(\sigma_{i,j}\), will be replaced by their respective estimators. 

A major difficulty of solving the optimality conditions is that the optimality equations \eqref{eq:input balance}-\eqref{eq:local balance} do not have closed-form solutions, and it is usually computationally expensive to solve them using numerical methods such as gradient descent. To improve computational efficiency, we design the two procedures tackling the optimality conditions in different ways. The DD-OCBA-approx procedure tries to directly solve the optimization problem \eqref{eq:opt LDR} at each iteration but approximating \eqref{eq:local balance} by assuming that a weighted ratio of allocation budget assigned to the optimal design is much larger than that assigned to other designs, which is a similar assumption taken by \cite{chen2000simulation}. This approximation enables us to compute the solution in a much simpler way. Alternatively, by taking a similar approach in \cite{chen2019complete}, the DD-OCBA-balance procedure avoids directly solving the optimality equations and instead balances the two sides of the equations, i.e., reduces the difference between two sides of the equations in each iteration when allocating the budget. Plausibly, DD-OCBA-approx is expected to converge faster since we solve the equations every time, while DD-OCBA-balance may converge slower since we only balance instead of solving the equations. However, DD-OCBA-approx approximates the optimality conditions, meaning that the ``optimal solution" we get may not be really optimal in the original problem. DD-OCBA-balance, instead, targets at the original problem and will eventually converge to the true optimal solution as more and more data are collected. The empirical comparison of these two methods will be carried out numerically in Section \ref{sec:numerical}.  

\subsection{DD-OCBA-approx}
In this section we derive the DD-OCBA-approx procedure. Let $\beta_i = \frac{\alpha_{i,j}}{ {\pi_j}\sigma_{i,j}} $,  $\ 1\le i\le K$, which is independent of $j$ by \eqref{eq:input balance} for optimal $\alpha$. Plugging $\alpha_{i,j}=\beta_i \sigma_{i,j} {\pi_j}$ into (\ref{eq:local balance}), we have
\begin{equation*}
    \frac{(\Bar{\mu}_b - \Bar{\mu}_i)^2}{ \frac{\sum_{j=1}^D\sigma_{b,j}{{\pi_j}}}{\beta_{b}}+\frac{\sum_{j=1}^D\sigma_{i,j}{{\pi_j}}}{\beta_{i}}} = \frac{(\Bar{\mu}_b - \Bar{\mu}_{i'})^2}{ \frac{\sum_{j=1}^D\sigma_{b,j}{{\pi_j}}}{\beta_{b}}+\frac{\sum_{j=1}^D\sigma_{i',j}{{\pi_j}}}{\beta_{i'}}}, \quad  i \neq i' \neq b.
\end{equation*}
Assume $\beta_b \gg \beta_i, \ \forall i\neq b$,  i.e., $ \frac{\alpha_{b,j}}{\sigma_{b,j}} \gg \frac{\alpha_{i,j}}{\sigma_{i,j}}\ \forall i\neq b, \forall j$,  the simulation budget assigned to the optimal design-input pair divided by its standard deviation is much larger than that of other designs.  Then we have $\frac{\beta_i}{\beta_{i'}} \approx \frac{\sum_{j=1}^D\sigma_{i,j}{{\pi_j}}/(\Bar{\mu}_b - \Bar{\mu}_i)^2}{\sum_{j=1}^D\sigma_{i',j}{{\pi_j}} /(\Bar{\mu}_b - \Bar{\mu}_{i'})^2} \ \forall i \neq i' \neq b$.
 Plugging $\beta_i=\frac{\alpha_{i,j}}{{\pi_j} {\sigma}_{i,j}}$ back with this approximation, we have
\begin{equation} \label{eq:approx opt cond}
    \frac{\alpha_{i,j}}{\alpha_{i^\prime, j^\prime}} = \frac{{\pi_j}\sigma_{i,j}\sum_{k=1}^D\sigma_{i,k}{\pi_{k}}/(\Bar{\mu}_b - \Bar{\mu}_i)^2}{{\pi_{j'}}{\sigma}_{i^\prime, j^\prime} \sum_{k=1}^D\sigma_{i',k}{{\pi_{k}}} /(\Bar{\mu}_b - \Bar{\mu}_{i'})^2}, \quad i,i^\prime \neq b.
\end{equation}
Furthermore, with (\ref{eq:total balance}) and $\sum_{i=1}^K\sum_{j=1}^K \alpha_{i,j} =1$, we can calculate $\alpha_{i,j}$ explicitly. 
 On a related note,
when input distributions are assumed to be known, i.e., $D=1$ and $\theta_1=\theta^c$,  \eqref{eq:approx opt cond} simplifies to:
$$ \frac{\alpha_{i,1}}{\alpha_{i',1}} = \frac{\sigma_{i,1}^2/(\mu_{b,1} -\mu_{i,1})^2}{\sigma_{i',1}^2/(\mu_{b,1} -\mu_{i',1})^2},$$
indicating that the allocation rule for a sub-optimal design is proportional to the ratio of its simulation variance to the square of its performance difference from the optimal design. This allocation principle aligns with the OCBA procedure described by \cite{chen2000simulation}. In contrast, when input distributions are unknown, the allocation rule for a sub-optimal design-input pair depends on not only the simulation variance under the same input parameter but also the variances under different input parameters. This is because the simulation outcomes under the same design but varying input parameters jointly influence the Bayesian average estimator.

 \textbf{DD-OCBA-approx }

\begin{enumerate}
    \item {\bf Input.} Number of designs $K$, input parameter space  $\Theta = \{\theta_1,\ldots,\theta_D \}$, initial sample size $n_0$, total simulation budget $n$, input data batch size $\{m(t)\}_{t=1}^\infty$ and stage-wise simulation budget $\{n(t)\}_{t=1}^\infty$, prior distribution $\pi^0$. 
    \item  {\bf Initialization.} Time stage counter $t \leftarrow 0 $, iteration counter $\ell\leftarrow 0$, total input data $M(t) \leftarrow 0$, posterior distribution $\pi\leftarrow \pi^0$. Collect $n_0$ initial samples for each design-input  pair $(i,j)$. Set $N_{i,j}^{(\ell)} = n_0 $. Compute the initial sample mean $\hat{\mu}_{i,j}^{(\ell)} = \frac{1}{N_{i,j}^{(\ell)}}\sum_{s=1}^{N_{i,j}^{(\ell)}} X_{i,j}^{(s)}$, and sample standard deviation $ \hat{\sigma}_{i,j}^{(\ell)} = \sqrt{\frac{1}{N_{i,j}^{(\ell)}-1}\sum_{s=1}^{N_{i,j}^{(\ell)}} (X_{i,j}^{(s)}-\hat{\mu}^{(\ell)}_{i,j})^2 } $.  $t\leftarrow t+1$.
    \item  {\bf WHILE} $\sum_{i=1}^K\sum_{j=1}^D N_{i,j}^{(\ell)} < n $ {\bf DO}
    \item Given input data of batch size $m(t)$, let $M(t) = \sum_{\tau=1}^t m(\tau)$ and update posterior distribution $\pi $. \label{line: input update} 
    \item Compute $\hat{\mu}_{i}^{(\ell)} = \sum_{j=1}^D \pi_j\hat{\mu}_{i,j}^{(\ell)}$. 
    \item {\bf REPEAT} n(t) {\bf TIMES} 
    \item   $\hat{b}^{(\ell)} \leftarrow \arg\max_i \hat{\mu}_i^{(\ell)}$. 
         \item Update $\hat{\alpha}^{(\ell)}_{i,j}$ using \eqref{eq:approx opt cond}, \eqref{eq:total balance} and $\sum_{i,j} \hat{\alpha}_{i,j} = 1$, with $ \Bar{\mu}_i, \sigma_{i,j}$ replaced by $\hat{\mu}_i^{(\ell)}$ and $\hat{\sigma}_{i,j}^{(\ell)}$, respectively. Calculate $\hat{N}_{i,j}^{(\ell)} = \hat{\alpha}^{(\ell)}_{i,j} \left( 1 + \sum\limits_{i=1}^K\sum\limits_{j=1}^D N_{i,j}^{(\ell)}\right), \quad \forall 1\le i\le K,\ 1\le j\le D.$
    \item Find the design-input pair $({I},J) = \arg\max_{i,j} \left(\hat{N}_{i,j}^{(\ell)} - N_{i,j}^{(\ell)}\right)$.  Simulate the pair $(I,J)$ once. Update $\hat{\mu}^{(l+1)}_{I, J}$, $\hat{\sigma}_{I,J}^{(l+1)}$ and $ \hat{\mu}^{(l+1)}_{I}$ using the new simulation output,  and set $\hat{\mu}^{(l+1)}_{i,j} = \hat{\mu}^{(\ell)}_{i,j}$, $\hat{\sigma}_{i,j}^{(l+1)} =\hat{\sigma}_{i,j}^{(\ell)}  $ and $ \hat{\mu}^{(l+1)}_{i}=\hat{\mu}^{(\ell)}_{i} $ for $i\neq I, j\neq J$. Let $N^{(l+1)}_{I, J} = N^{(\ell)}_{I, J}+1$ and $N^{(l+1)}_{i,j} = N^{(\ell)}_{i,j}$ for all $i\neq I, j\neq J$. 
    \item $l \leftarrow l+1$. $t\leftarrow t+1$.
    \item {\bf END REPEAT}
    \item {\bf END WHILE}
    \item {\bf Output:} Output $i_b = \arg\max_i \hat{\mu}_i^{(\ell)}$ as the best design.
\end{enumerate}            
 \subsection{DD-OCBA-balance}
Unlike DD-OCBA-approx where we try to directly solve for the optimal solutions, DD-OCBA-balance only requires to evaluate both sides of the three optimality equations given the current number of replications for each design-input pair. The procedure selects a design-input pair each time to reduce the difference (balance) of at least one of the optimality equations. In particular, at each iteration, the procedure will first decide whether to simulate the estimated best design or one of the non-optimal designs to balance the ``total balance" conditions. If the estimated best design is not selected, then the procedure selects a non-optimal design to balance the ``total balance" conditions.
After selecting the design, an input realization is chosen by balancing the ``input balance" conditions. Notice that although in \eqref{eq:input balance} the ``input balance" conditions are only for non-optimal designs, \eqref{eq:input balance} also holds for $i = b$ by \eqref{eq:total balance} with $\frac{\alpha_{i,j}}{\sigma_{i,j}} $ replaced by $ \frac{\alpha_{i,j'}}{\sigma_{i,j'}} \frac{{\pi_j}}{{\pi_{j'}}}, \ \forall j\neq j'$. The balancing approach utilizes the monotonicity of both sides of all optimality equations in terms of the allocation policy $\alpha_{i,j}$. For example, if we have in one of the equations in (\ref{eq:input balance}) violated by $\frac{\alpha^{(\ell)}_{i,j}}{\sigma_{i,j}{\pi_j}} < \frac{\alpha_{i,j'}^{(\ell)}}{\sigma_{i,j'}{\pi_{j'}}} $, then we may want to simulate the design-input pair $(i,j)$ to make the left hand side larger. The DD-OCBA-balance procedure is presented as follows:

\textbf{DD-OCBA-balance}
\begin{enumerate}
    \item {\bf Input.} Number of designs $K$, input parameter space  $\Theta = \{\theta_1,\ldots,\theta_D \}$, initial sample size $n_0$, total simulation budget $n$, input data batch size $\{m(t)\}_{t=1}^\infty$, and stage-wise simulation budget $\{n(t)\}_{t=1}^\infty$, prior distribution $\pi^0$. 
   \item  {\bf Initialization.} Time stage counter $t \leftarrow 0 $, iteration counter $\ell\leftarrow 0$, total input data $M(t) \leftarrow 0$, posterior distribution $\pi\leftarrow\pi^0$. Collect $n_0$ initial samples for each design-input  pair $(i,j)$. Set $N_{i,j}^{(\ell)} = n_0 $. Compute the initial sample mean $\hat{\mu}_{i,j}^{(\ell)} = \frac{1}{N_{i,j}^{(\ell)}}\sum_{s=1}^{N_{i,j}^{(\ell)}} X_{i,j}^{(s)}$, and sample standard deviation $ \hat{\sigma}_{i,j}^{(\ell)} = \sqrt{\frac{1}{N_{i,j}^{(\ell)}-1}\sum_{s=1}^{N_{i,j}^{(\ell)}} (X_{i,j}^{(s)}-\hat{\mu}^{(\ell)}_{i,j})^2 } $.  $t\leftarrow t+1$.
    \item  {\bf WHILE} $\sum_{i=1}^K\sum_{j=1}^D N_{i,j}^{(\ell)} < n $ {\bf DO} 
    \item Given input data of batch size $m(t)$, let $M(t) = \sum_{\tau=1}^t m(\tau)$ and update posterior distribution $\pi$.
    \item Compute $\hat{\mu}_{i}^{(\ell)} = \sum_{j=1}^D \pi_j \hat{\mu}_{i,j}^{(\ell)}$.
    \item {\bf REPEAT} n(t) {\bf TIMES} 
    \item   $\hat{b}^{(\ell)} \leftarrow \arg\max_i \hat{\mu}_i^{(\ell)}$. 
        \item Let ${j^*} = \arg\max_j\left| \left(\frac{N^{(\ell)}_{\hat{b}^{(\ell)},j}}{\hat{\sigma}^{(\ell)}_{\hat{b}^{(\ell)},j}}\right)^2 -\sum\limits_{i\neq \hat{b}^{(\ell)}} \left(\frac{N^{(\ell)}_{i,j}}{\hat{\sigma}^{(\ell)}_{i,j}}\right)^2   \right| $.
    \item \textbf{IF} $\left(\frac{N^{(\ell)}_{\hat{b}^{(\ell)},{j^*}}}{\hat{\sigma}^{(\ell)}_{\hat{b}^{(\ell)},{j^*}}}\right)^2 -\sum\limits_{i\neq \hat{b}^{(\ell)}} \left(\frac{N^{(\ell)}_{i,{j^*}}}{\hat{\sigma}^{(\ell)}_{i,{j^*}}}\right)^2 < 0$,  set $I = \hat{b}^{(\ell)}$, $J = \arg\min_j \frac{N^{(\ell)}_{\hat{b}^{(\ell)},j}}{\hat{\sigma}^{(\ell)}_{\hat{b}^{(\ell)},j} \pi_j}$. \label{item:hatb,J}
    \item \label{item:suboptimal,J}\textbf{ELSE} set $I = \arg\min\limits_{i \neq \hat{b}^{(\ell)}} \frac{(\hat{\mu}^{(\ell)}_{\hat{b}^{(\ell)}} - \hat{\mu}_i^{(\ell)})^2}{\sum\limits_{j=1}^D \frac{(\hat{\sigma}^{(\ell)}_{\hat{b}^{(\ell)},j})^2\pi_j^2}{N_{\hat{b}^{(\ell)},j}^{(\ell)}} + \sum\limits_{j=1}^D\frac{(\hat{\sigma}_{i,j}^{(\ell)})^2 \pi_j^2 }{N_{i,j}^{(\ell)}}} $, $J = \arg\min_j \frac{N_{I,j}^{(\ell)}}{\hat{\sigma}^{(\ell)}_{I,j}\pi_j}. $
    \item \textbf{END IF}
    \item Simulate the pair $(I,J)$ once. Update $\hat{\mu}^{(l+1)}_{I, J}$, $\hat{\sigma}_{I,J}^{(\ell)}$ and $ \hat{\mu}^{(l+1)}_{I}$ using the new simulation output,  and set $\hat{\mu}^{(l+1)}_{i,j} = \hat{\mu}^{(\ell)}_{i,j}$, $\hat{\sigma}_{i,j}^{(l+1)} =\hat{\sigma}_{i,j}^{(\ell)}  $ and $ \hat{\mu}^{(l+1)}_{i}=\hat{\mu}^{(\ell)}_{i} $ for $i\neq I, j\neq J$. Let $N^{(l+1)}_{I, J} = N^{(\ell)}_{I, J}+1$ and $N^{(l+1)}_{i,j} = N^{(\ell)}_{i,j}$ for all $i\neq I, j\neq J$. 
    \item $l \leftarrow l+1$.  $t\leftarrow t+1$.
    \item {\bf END REPEAT}
    \item {\bf END WHILE}
\end{enumerate}  


\section{CONSISTENCY AND ASYMPTOTIC OPTIMALITY} \label{sec:consistency}


\subsection{Consistency}
A R\&S algorithm is consistent if it selects the true optimal design $b^c$ (as defined in Assumption \ref{assump: unique b}) as time stage $t$ goes to infinity almost surely. This happens if the estimated Bayesian average performance $\hat{\mu}_i^{(\ell)}$ converges to its true value $\mu_i(\theta^c)$ for each $i$ almost surely. Let $\pi^t = (\pi^t_j)_{j=1}^D$ denote the posterior distribution at stage $t$. In order to show the consistency, we need to (i) show  the convergence of posterior distribution $\pi^t$ that is estimated from input observations and (ii) the convergence of $\hat{\mu}_{i,j^c}^{(\ell)}$, where $j^c$ is the index such that $\theta_{j^c} = \theta^c$. 

 We first make the following assumption on the identifiability on the input parameters.
\begin{assumption} \label{assump: identifiable finite}
\textbf{(Identifiability)} 
For $\theta \neq \theta' \in \Theta$, $F_\theta \neq F_{\theta'}$.
\end{assumption}

Together with \ref{assump: prior likelihood}.1, Assumption \ref{assump: identifiable finite} guarantees the Bayesian consistency of the posterior distribution $\pi$ by Doob's consistency theorem (see, e.g., Theorem 10.10 in \cite{van2000asymptotic}). The Bayesian consistency states if the number of i.i.d. input data goes to infinity, then the posterior distribution $\pi$ converges to a Dirac measure centered at $\theta^c$. 


Moreover, we make the following assumption about the input data batch size and simulation budget in each stage for the DD-OCBA-approx and DD-OCBA-balance procedures.
\begin{assumption} \label{assump:batch limit}
The stage-wise input data batch size $m(t)$ and simulation budget $n(t)$ satisfy 
$$ \lim_{T\rightarrow \infty} \sum_{t=1}^T n(t) = \infty, \quad  \lim_{T\rightarrow \infty} \sum_{t=1}^T m(t) = \infty.$$ 
\end{assumption}
Assumption \ref{assump:batch limit} ensures that both the total amount of input data and  total simulation replications go to infinity as time stage $T$ goes to infinity, which helps guarantee the consistency of the posterior distribution $\pi^t$ as well as the design-input performance estimate $\hat{\mu}_{i,j^c}^{(\ell)}$ under the true parameter $\theta^c$. Recall that $\ell$ is the iteration counter (i.e. the amount of simulation budget that has been assigned from stage $1$ by the algorithm). 
The following theorem shows the consistency of DD-OCBA-approx and DD-OCBA-balance.

\begin{theorem} \textbf{(Consistency)} \label{thm:consistency}
     Suppose Assumptions \ref{assump: parametric input}-\ref{assump: unique b}, Assumption \ref{assump:output}.2 and  
 Assumption \ref{assump:input data}-\ref{assump:batch limit} hold. Then,
     \begin{enumerate}
         \item DD-OCBA-approx selects the optimal design $b^c$ almost surely as $\ell \rightarrow\infty$.
         \item DD-OCBA-balance selects the optimal design $b^c$ almost surely as $\ell \rightarrow \infty$.
     \end{enumerate}
\end{theorem}

\subsection{Asymptotic Optimality}
In the last section, we establish the consistency of the two proposed algorithms, ensuring they select the optimal design as the simulation budget approaches infinity. However, it remains unclear if the algorithm's allocation policy will converge to the optimal one. Thus, we examine the algorithms' asymptotic optimality in this section. An algorithm is said to be asymptotically optimal if the allocation policy given by the algorithm converges to the \textbf{limiting} optimal policy as the total simulation budget goes to infinity. Here the allocation policy is defined as $\{\alpha_{i,j}^{(\ell)} = \frac{N_{i,j}^\ell}{N^{(\ell)}}\}$, where $N_{i,j}^{(\ell)}$ is the simulation budget assigned to $(i,j)$ pair up to iteration $\ell$ and $N^{(\ell)}$ is the total simulation budget used up to iteration $\ell$.  Notably, the optimal policy is computed by maximizing the decay rate in \eqref{eq:LDR} given a current posterior distribution $\pi$, which is updated periodically and converges to a Dirac measure, i.e., $\pi_{j^c} \rightarrow 1$. 
Hence, the \textbf{limiting} optimal policy $\alpha^*$ is defined to be the limit of that optimal policy  that maximizes \eqref{eq:LDR} as $\pi_{j^c} \rightarrow 1$. 

\subsubsection*{Limiting optimal policy.} From \eqref{eq:rate formula}, we have 
\begin{equation} \label{eq: limit rate}
    \mathbf{G}_i(\alpha_b,\alpha_i) =  \frac{(\Bar{\mu}_b - \Bar{\mu}_i)^2}{2\left( \sum_{j=1}^D \frac{\sigma_{b,j}^2{\pi_j}^2}{\alpha_{b,j}} + \sum_{j=1}^D\frac{\sigma_{i,j}^2{\pi_j}^2}{\alpha_{i,j}}\right)} \rightarrow  \frac{({\mu}^c_{b^c} - {\mu}^c_i)^2}{2\left(  \frac{\sigma_{b^c,j^c}^2}{\alpha_{b^c,j^c}} +\frac{\sigma_{i,j^c}^2}{\alpha_{i,j^c}}\right)} \text{ as } \pi_{j^c} \rightarrow 1,
\end{equation}
where $\mu_i^c:=\mu_i(\theta^c), i\in\mathcal{I}$ denote the true expected performance under $\theta^c$, $b$ is the design that maximizes the Bayesian average performance and $b^c$ is the design that maximizes the true expected performance.
\eqref{eq: limit rate} indicates the limiting rate function only depends on allocation ratio $\{\alpha_{i,j^c}\}_{i\in\mathcal{I}}$ as the Bayesian average performance converges to the true expected performance, which only depends on the design-input pair under the true input parameter $\theta^c$. Moreover,
the limiting rate function \eqref{eq: limit rate} coincides with the classic R\&S with known input distribution (e.g., in \cite{chen2000simulation,glynn2004large}), since the uncertainty of the input distribution decreases to $0$ in the limit. Then, the limiting optimal policy can be easily computed as a special case ($D=1$) of Theorem \ref{thm:optimal allocation}, which we formally states in the following theorem that guarantees the asymptotic optimality of DD-OCBA-approx.

\begin{theorem} \textbf{(Asymptotic Optimality for DD--OCBA-approx)} \label{thm: optimality approx}
Suppose Assumptions \ref{assump: parametric input}- \ref{assump: unique b}, \ref{assump:output}.2, \ref{assump:input data}-\ref{assump:batch limit} hold. Then, almost surely, for DD-OCBA-approx, 
\begin{enumerate}
\item $\alpha^{(\ell)}_{i,j} \rightarrow 0, \forall i \in \mathcal{I}, j \neq j^c$ as $\ell \rightarrow \infty$,
    \item  $\lim_{\ell\rightarrow \infty} \alpha_{i,j^c}^{(\ell)} = \alpha^*_{i,j^c}, \forall i \in \mathcal{I}$ where $\alpha^*_{i,j^c}, i\in \mathcal{I}$ satisfy
    \begin{flalign}
    &\romannumeral 1. \ \left(\frac{\alpha^*_{b^c,j^c}}{\sigma_{b^c,j^c}}\right)^2 = \sum_{i\neq b^c} \left(\frac{\alpha_{i,j^c}^*}{\sigma_{i,j^c}^*} \right)^2&& \label{thmeq:total balance approx}\\
     &\romannumeral 2.\  \frac{\alpha_{i,j^c}^*}{\alpha_{i',j^c}^*} = \frac{\sigma_{i,j^c}^2/(\mu^c_{b^c}-\mu^c_i)^2}{\sigma_{i',j^c}^2/(\mu^c_{b^c} - \mu^c_{i'})^2} \quad \forall i\neq i' \neq b^c &&  \label{thmeq: local balance approx}
\end{flalign}
\end{enumerate}
\end{theorem}

Theorem \ref{thm: optimality approx} guarantees the asymptotic optimality for DD-OCBA-approx. For DD-OCBA-balance, we introduce the following additional assumption regarding the parametric log likelihood ratio and the average batch size of input data and average stage-wise simulation budget.
\begin{assumption} \label{assump: optimality balance} \ 
\begin{enumerate}
    \item The strong law of large number (SLLN) holds for the log likelihood ratio sequence $\left\{\log \frac{f_\theta(\xi_k)}{f_{\theta^c}(\xi_k)}\right\}_{k=1}^\infty, \forall \theta\neq \theta^c \in \Theta$, where $\xi_k \sim F_{\theta^c}, k=1,2,\ldots$ are i.i.d. sequence. That is, almost surely for every sequence $\{\xi_k\}_{k=1}^\infty$,
    $$ \lim_{m\rightarrow\infty} \frac{1}{m} \sum_{k=1}^m \log\frac{f_\theta(\xi_k)}{f_{\theta^c}(\xi_k)} = \mathbb{E}_{\theta^c}\left[ \log \frac{f_\theta(\xi)}{f_{\theta^c} (\xi)}\right] = -\operatorname{KL}(\theta^c\|\theta) <0,$$
    where $\operatorname{KL}(\cdot\|\cdot)$ is the Kullback–Leibler divergence (KL) divergence.
    \item There exist $\bar{n}$, $\bar{m} >0$, such that $\lim_{T\rightarrow \infty} \frac{1}{T}\sum_{t=1}^T n(t) =\bar{n}$, $\lim_{T\rightarrow \infty} \frac{1}{T}\sum_{t=1}^T m(t) =\bar{m}$ almost surely.
\end{enumerate}
    
\end{assumption}
Assumption \ref{assump: optimality balance}.1 imposes a stronger assumption on the likelihood function, whereas Assumption \ref{assump: optimality balance}.2 imposes an additional constraint on the data batch size, which ensures that the total input data and simulation replications increase at the same rate $O(T)$ as time stage $T$ increases to infinity. The purpose of introducing Assumption \ref{assump: optimality balance} is to guarantee a certain (almost sure) convergence rate for both the posterior distribution and the estimators of design-input performance and variance, in addition to  consistency of the posterior distribution and the estimators. 

There are several technical reasons why we need the extra Assumption \ref{assump: optimality balance} to prove the asymptotic optimality for DD-OCBA-balance. We provide the intuition for the most important reason in the following. In DD-OCBA-approx, the stage-wise optimal allocation policy can be computed explicitly as a smooth and bounded function of parameters including the current posterior distribution, design-input performance estimator and design-input variance estimator. When all of these parameters varies only in a small neighborhood of their limiting value (consistency), the computed stage-wise optimal policy also only varies little around the limiting optimal policy. This is a key property in the proof of Theorem \ref{thm: optimality approx}. Nonetheless, in DD-OCBA-balance, we do not have an explicit-form solution to guide the allocation policy, but instead only evaluate the two sides of the optimality equations in Theorem \ref{thm:optimal allocation}. Unlike DD-OCBA-approx where the stage-wise optimal policy is a smooth function of the aforementioned parameters (posterior, mean and variance), values of two sides of Input Balance condition is not smooth or bounded in terms of the posterior probability $\pi_j$ for $j\neq j^c$. To be specific, as $\pi^{t_\ell}_j \rightarrow 0$, the limit of 
$\frac{\alpha_{i,j}^{(\ell)}}{\sigma^{(\ell)}_{i,j} \pi^{t_\ell}_j}$  can be either $0,\infty$ or a finite value,
depending on how fast of both $\alpha_{i,j}^{(\ell)}$ and $\pi^{t_\ell}_j$ converge to zero. Moreover, the dependence on the past allocation policy $\alpha_{i,j}^{(\ell)}$ and the complex form of the Local Balance condition in \eqref{eq:local balance} also extremely complicates the proof. Hence, we require some stronger assumptions to guarantee the same convergence rate of the different estimators.

The following lemma states the converge rate of the posterior probability with Assumption \ref{assump: optimality balance}. 
\begin{lemma} \label{lem: posterior convergence}
Suppose Assumptions \ref{assump: parametric input}, \ref{assump:input data}, \ref{assump: optimality balance} hold. Then, there exists $\kappa > 0$, for almost every sequence of input data $\{\xi_s\}_{s=1}^\infty$,
\begin{equation*}
    \pi_{j^c}^t = 1 - O\left( \exp\left( -\kappa t\right)\right).
\end{equation*}
As a result, for $j\neq j^c$, $0<\pi_j^t \le 1 - \pi_{j^c}^t = O\left( \exp\left( -\kappa t\right)\right)$.
\end{lemma}
Lemma \ref{lem: posterior convergence} implies that the posterior probability converges exponentially fast almost surely, which can help guarantee that any design-input pair under input parameter $\theta_j \neq \theta^c$ will only be simulated at most finitely many times.

We are now ready to provide the asymptotic optimality for DD-OCBA-balance. In addition, we can further characterize the convergence rate of the allocation policy to the limiting optimal allocation policy, for both DD-OCBA-approx and DD-OCBA-balance..
\begin{theorem} \label{thm: convergence speed}
    Suppose Assumptions  \ref{assump: parametric input}-\ref{assump: unique b}, \ref{assump:output}.2,
 \ref{assump:input data} and \ref{assump: optimality balance} hold. Then, 
    \begin{enumerate}
        \item For DD-OCBA-approx, $|\alpha^{(\ell)}_{i,j^c} - \alpha^*_{i,j^c}| = O\left(\sqrt{\frac{\log \log \ell}{\ell}}\right), \ \forall 1\le i\le K$ almost surely, where  $\alpha_{i,j^c}^*$ satisfies the optimality conditions  \eqref{thmeq:total balance approx} and \eqref{thmeq: local balance approx}.
        \item For DD-OCBA-balance, 
         \begin{flalign}
    &\romannumeral 1. \ \left|\left(\frac{\alpha_{b^c,j^c}^{(\ell)}}{\sigma_{b^c,j^c}}\right)^2 - \sum\limits_{i\neq b^c}^D \left(\frac{\alpha_{i,j^c}^{(\ell)}}{\sigma_{i,j^c}}\right)^2\right| = O\left(\sqrt{\frac{\log \log \ell}{\ell}}\right) \quad \text{ almost surely } ,&& \label{thmeq:error total balance}\\
    &\romannumeral 2. \ \left|\frac{(\mu^c_b - \mu^c_i)^2}{ \frac{\sigma_{ b^c,j^c}^2 }{\alpha_{ b^c,j^c}^{(\ell)}}+\frac{\sigma_{i,j^c}^2 }{\alpha_{i,j^c}^{(\ell)}}} - \frac{(\mu^c_b - \mu^c_{i'})^2}{ \frac{\sigma_{ b^c,j^c}^2 }{\alpha_{ b^c,j^c}^{(\ell)}}+\frac{\sigma_{i',j^c}^2 }{\alpha_{i',j^c}^{(\ell)}}} \right| =  O\left(\sqrt{\frac{\log \log \ell}{\ell}}\right)  \quad \forall i \neq i' \neq  b^c \quad \text{ almost surely. } && \label{thmeq:error local balance}
\end{flalign}
    \end{enumerate}
\end{theorem}

\section{Extension to Continuous Parameter Space} \label{sec: continuous space}
So far we have assumed that the input parameter space $\Theta$ is a finite set, which can be regarded as a discretization of the original parameter space. Nonetheless, this may lead to the issue of model mis-specification, which means the true input parameter $\theta^c$
may not belong to the descretized finite set. Although the posterior distribution $\pi$ still converges to a Dirac delta measure concentrated on $\theta_{j^c}$, which minimizes the Kullback-Leibler divergence from the true distribution $F_{\theta^c}$ to the set of distributions $\{F_{\theta_j}, j=1,\ldots,D\}$ (e.g., see \cite{lian2009rates}), the best design under $F_{\theta_{j^c}}$ may not be the true optimal design under $F_{\theta^c}$ due to the discrepancy between $\theta^c$ and $\theta_{j^c}$. To address the issue of model mis-specification, in this section we extend the proposed methods to the general continuous input parameter space by generalizing the discretization approach.

To be specific, let $\Theta_1,\Theta_2,\ldots,\Theta_D$ to be $D$ sub-space of $\Theta$ (entire parameter space) such that $\Theta_j\cap \Theta_{j'} = \O$ and $\cup_{j=1}^D \Theta_j = \Theta$, i.e., $\Theta_1,\Theta_2,\ldots,\Theta_D$ is a partition of $\Theta$. Then, the true parameter $\theta^c \in \Theta_{j^c}$ for some $1\le j^c\le D$. In addition, when $\Theta_j =\{\theta_j\}$ is a singleton set, then this reduces to the setting in Section \ref{sec: problem formulation}.

Since now the parameter space $\Theta$ can be continuous, the posterior (density) distribution is computed as
$$ \pi(\theta|\xi_1,\ldots,\xi_m) = \frac{\pi_0(\theta) \prod_{\ell=1}^m f_{\theta}(\xi_\ell)}{\int_{\Theta} \pi_0(\theta') \prod_{\ell=1}^m f_{\theta'}(\xi_\ell) \mathrm{d}\theta'}.$$

We further define
$$ \pi_j := \pi(\Theta_j|\xi_1,\ldots,\xi_m) = \frac{\int_{\Theta_j} \pi_0(\theta) \prod_{\ell=1}^m f_{\theta}(\xi_\ell)\mathrm{d}\theta}{\int_{\Theta} \pi_0(\theta) \prod_{\ell=1}^m f_{\theta}(\xi_\ell) \mathrm{d}\theta}$$
and 
$$ \Tilde{\pi}_j(\theta) = \frac{\pi_0(\theta)}{\pi_j} \prod_{\ell=1}^m f_{\theta}(\xi_\ell).$$
That is, $\pi_j$ is the posterior probability of $\{\theta^c \in \Theta_j\}$ and $\Tilde{\pi}_j$ is the conditional posterior probability of $\theta^c$ conditioned on the event $\theta^c \in \Theta_j$.
Then the Bayesian average performance can be expressed as 
$$ \Bar{\mu}_{i} = \mathbb{E}_\pi[\mu_i(\theta)] = \sum_{j=1}^D \pi_j \mathbb{E}_\pi[\mu_i(\theta)|\theta\in\Theta_j] = \sum_{j=1}^D \pi_j \mathbb{E}_{\Tilde{\pi}_j}[\mu_i(\theta)].$$
Instead of simulating a design under a fixed input parameter $\theta_j$, we will then simulate a design under a fixed subspace $\Theta_j$. To generate one simulation output under the design-input pair $(i,j)$  (here input $j$ refers to subspace $\Theta_j$, as opposed to a single input parameter, $\theta_j$, previously), we first generate $\theta \sim \Tilde{\pi}_j$ and then under the input distribution $F_\theta$  we generate one simulation output. Notably, both of the proposed algorithms, DD-OCBA-approx and DD-OCBA-balance, can then be directly applied to allocate simulation budget among the design-input pairs, with the only differences on the posterior updating and generation of simulation outputs. To differentiate with the setting of finite parameter space, we name the two procedures as DD-OCBA-approx-C and DD-OCBA-balance-C, respectively, where ``C" stands for continuous.

\textbf{Remark:} In practice, computing the exact posterior probability $\pi_j$ and sampling from the conditional posterior distribution $\Tilde{\pi}_j$ can be computationally expensive. To improve the computation efficiency, given the current posterior distribution $\pi$, one can first compute the Bayesian mean estimator $\Bar{\theta} = \mathbb{E}_\pi[\theta]$
and then let $\theta_j = \operatorname{Proj}_{\Theta_j} (\Bar{\theta})$, i.e., $\theta_j \in \Theta_j$ is the closest point in $\Theta_j$ to $\Bar{\theta}$ (with respect to, e.g., $l_2$ distance). Then, one can run the simulation for $(i,j)$ pair under input parameter $\theta_j$. Furthermore, we approximate the posterior probability $\pi_j$ by the the posterior density at $\theta_{j'},j'=1,\ldots,D$ as $\pi_j \approx \frac{\pi(\theta_j)}{\sum_{j'=1}^D \pi(\theta_{j'})} $.

\subsection*{Convergence Analysis}
In this section, we show that with continuous parameter space, Theorem \ref{thm:consistency} and Theorem \ref{thm: optimality approx} still hold. That is, we can still prove the consistency for DD-OCBA-approx-C and DD-OCBA-balance-C and asymptotic optimality for DD-OCBA-approx-C. 

Before jumping into the analysis, we first give some intuitive explanation on the difference of proving the consistency/asymptotic optimality with continuous parameter space. Recall the only differences between DD-OCBA-approx(balance) and DD-OCBA-approx(balance)-C are the way of generating simulation output and updating the posterior distribution, which affects the convergence of posterior probability $\pi_j$, performance estimator $\hat{\mu}_{i,j}^{(\ell)}$ and variance estimator $\hat{\sigma}_{i,j}^{(\ell)}$. As a result, the same analysis for consistency and asymptotic optimality can be applied if these aforementioned estimators or posterior probabilities satisfy certain properties as required in the proof of Theorem \ref{thm:consistency}-\ref{thm: optimality approx}, which include the strong consistency of the posterior distribution as well as the performance/variance estimators.
With the finite parameter space, for a fixed design-input pair the simulation outputs are i.i.d.. However, with continuous parameter space, as the posterior distribution is being updated, the conditional posterior distribution $\Tilde{\pi}_j$ also varies. Consequently, the simulation outputs under a fixed design-input (subspace) pair is not i.i.d.. across stages, requiring us to strengthen the current analysis.

 Let $\pi^t$, $\pi_j^t$ and $\Tilde{\pi}^t_j$ be the posterior distribution, posterior probability for $\{\theta^c \in \Theta_j\}$, and conditional posterior distribution for  $\{\theta^c: \theta^c \in \Theta_j\}$ at stage $t$, respectively. We make the following assumptions for proving consistency of the algorithms.
\begin{assumption} \label{assump: consistency continuous space} \ 
\begin{enumerate}

    \item $\Theta \in \mathbb{R}^d, d<\infty$ is convex and compact. Furthermore, there exists $1\le j^c\le D$, such that $\theta^c \in \operatorname{int}(\Theta_{j^c})$.
    \item $y_i^{(k)}(\theta) := \mathbb{E}[X^k_i(\theta)|\theta]$ is continuous in $\theta$ for $i\in\mathcal{I}$ and $k=1,2,4$.
    \item At iteration $\ell$, the simulation output for some design-input pair  $(i,j)$ is generated by first sampling $\theta^{(\ell)} \sim \Tilde{\pi}_j^{(t_\ell)}$ and then simulating design $i$ under $\theta^{(\ell)}$. Moreover, conditioned on $\Tilde{\pi}_j^{(t_\ell)}$, the simulation output is independent of the past input data and past simulation outputs.
\end{enumerate}
\end{assumption}
Assumption \ref{assump: consistency continuous space}, together with Assumption \ref{assump: prior likelihood} and Assumption \ref{assump: identifiable finite}, guarantees the consistency of  posterior distribution as well as the consistency of performance and variance estimator, as formally stated in the following Lemma \ref{lem: consistency continuous}. 

\begin{lemma} \label{lem: consistency continuous}
Suppose Assumption \ref{assump: parametric input}.1, \ref{assump: prior likelihood}, \ref{assump: identifiable finite} and \ref{assump: consistency continuous space} hold.
\begin{enumerate}
    \item If $\lim_{t\rightarrow\infty} M(t) = \infty$, then
$\pi_{j^c}^t \rightarrow 1$ as $t\rightarrow\infty $ for almost every sequence of input data.
\item If  also $\lim_{\ell\rightarrow\infty} N_{i,j^c}^{(\ell)}=\infty$ for some design $i\in\mathcal{I}$, then $\lim_{\ell \rightarrow \infty} \hat{\mu}_{i,j^c}^{(\ell)} =\mu_i(\theta^c)$ and $\lim_{\ell \rightarrow \infty} \hat{\sigma}_{i,j^c}^{(\ell)} =  \sigma_i(\theta^c)$ almost surely.
\end{enumerate}
 \end{lemma}
 With Lemma \ref{lem: consistency continuous}, we can then follow the same proofs of Theorem \ref{thm:consistency} and Theorem \ref{thm: optimality approx} to get the following consistency result for DD-OCBA-approx-C and DD-OCBA-balance-C.
 
\begin{corollary}\label{cor: consistency continuous}
    Suppose Assumption \ref{assump: parametric input}.1,  \ref{assump: prior likelihood}, \ref{assump: unique b}, \ref{assump:output}.2, \ref{assump:input data}, \ref{assump: identifiable finite}, \ref{assump:batch limit}, \ref{assump: consistency continuous space} hold. Then, both DD-OCBA-approx-C and DD-OCBA-balance-C selects the true optimal design $b^c := \arg\max_{i}\mu_i(\theta^c)$ almost surely. Furthermore, the allocation rule $\alpha_{i,j}^{(\ell)}, \forall $ given by
    DD-OCBA-approx-C converges to the limiting optimal allocation policy defined in Theorem \ref{thm: optimality approx}.
\end{corollary}

\section{NUMERICAL EXPERIMENT}\label{sec:numerical}
\subsection{Comparison Baselines}

We test the performance of DD-OCBA-approx(-C) and DD-OCBA-balance(-C) by comparing with (i) Equal Allocation, which allocates an equal simulation budget to all design-input pairs, and (ii) Adaptive OCBA, which is an extension of the OCBA (see \cite{chen2000simulation}) algorithm to our data-driven setting. OCBA sequentially decides the next system to simulate based on the past simulation samples, whereas Adaptive OCBA uses the same budget allocation rule as OCBA, but with simulation samples generated under the current estimated input distribution each time as opposed to a fixed input distribution in OCBA. In our context, with the current posterior distribution $\pi$, a simulation output for design $i$ is generated by first sampling $\theta \sim \pi$ and then running the simulation procedure under this parameter $\theta$ to get a sample of $X_i(\theta)$.  In fact, Adaptive OCBA can be viewed as a special case of DD-OCBA-approx-C with number of sub-space $D=1$ and the entire space $\Theta$ to be finite.
On a related note, we do not explicitly compare with the SEIU algorithm by \cite{Wu2022data}, which also considers the R\&S with streaming input data but took a fixed confidence formulation. From the perspective of budget allocation, they simply equally allocate the budget to all designs with elimination until only one is left with a given confidence level. Hence, applying their method to the setting of fixed budget R\&S is equivalent to using Equal Allocation.  

We test the procedures on two different problems, beginning with a quadratic problem with finite input parameter space and then on a portfolio problem with continuous parameter space.
\subsection{Quadratic Problem}
Consider the following optimization problem.
$$ \min_{i\in \mathcal{I}}~ \mathbb{E}_{\theta^c} \left[(x_i-\zeta)^2\right],$$ where $x_i = \theta^c + i-1, i \in \mathcal{I} = \{1,2,\ldots,20\}$ and $\zeta$ follows an exponential distribution with  unknown mean $\theta^c = 3$. Hence, the true best design $b^c = 0$. The input parameter space $\Theta = \{\theta_j: \theta_j = 0.5\cdot j, j=1,2,\ldots,20\}$. In the following we set the initial number of simulation for each design-input pair $n(0) = 10$. The initial prior is uninformative, that is, $\pi^0_j = \frac{1}{20}, j=1,\ldots,20.$ We first test with constant stage-wise input data batch size and simulation budget. We set the stage-wise simulation budget $n(t) = 100$ and vary the input data batch size $m(t) \in \{5,10,15\}$. 
\subsubsection*{Experiment Results}

\begin{figure}[h]
        \begin{subfigure}[t]{0.3\textwidth}
        \centering
        \includegraphics[height=1.5in]{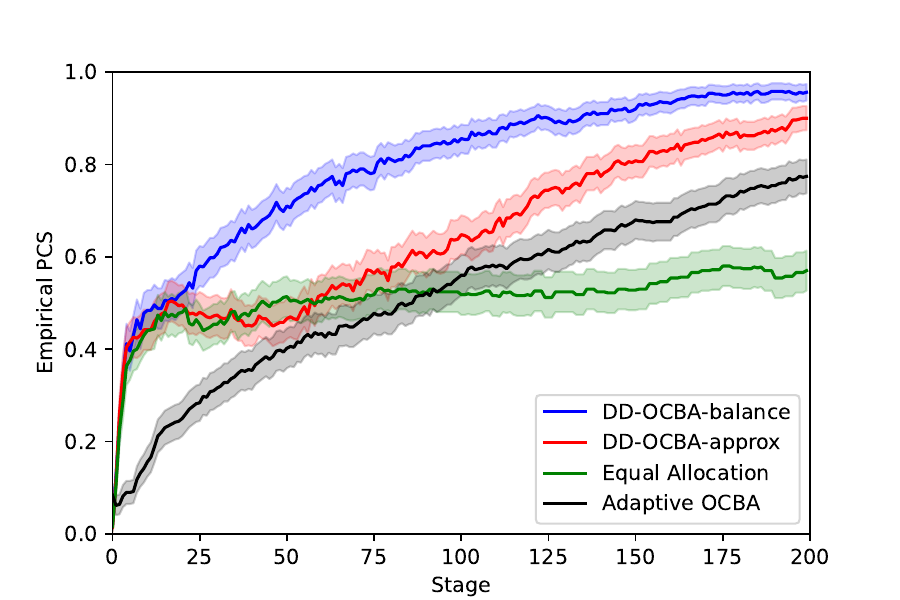}
        \caption{$m(t)=5$}
        \label{fig:quad_const_mt_5}
    \end{subfigure}
    ~~
        \begin{subfigure}[t]{0.3\textwidth}
        \centering
        \includegraphics[height=1.5in]{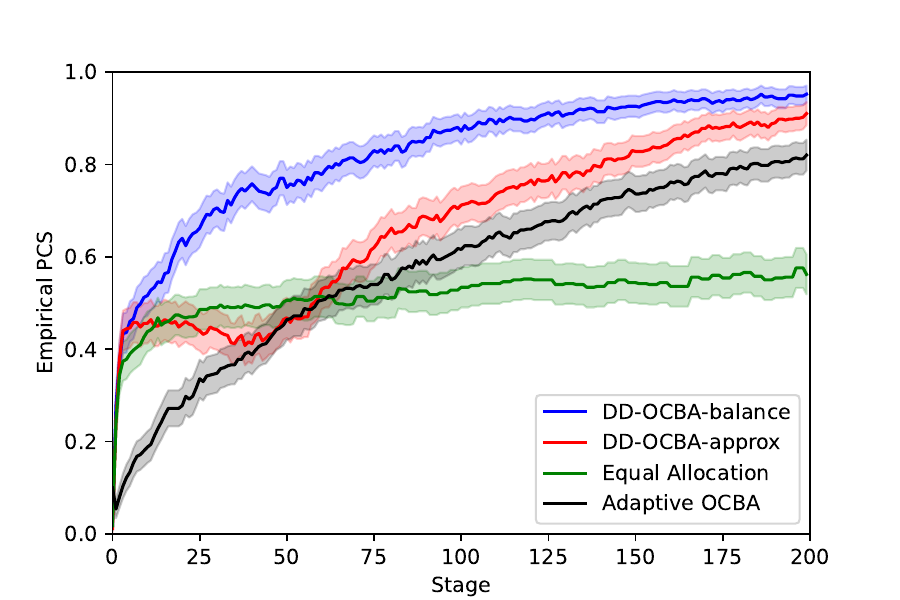}
        \caption{$m(t)=10$}
        \label{fig:quad_const_mt_10}
    \end{subfigure}
    ~~
        \begin{subfigure}[t]{0.3\textwidth}
        \centering
        \includegraphics[height=1.5in]{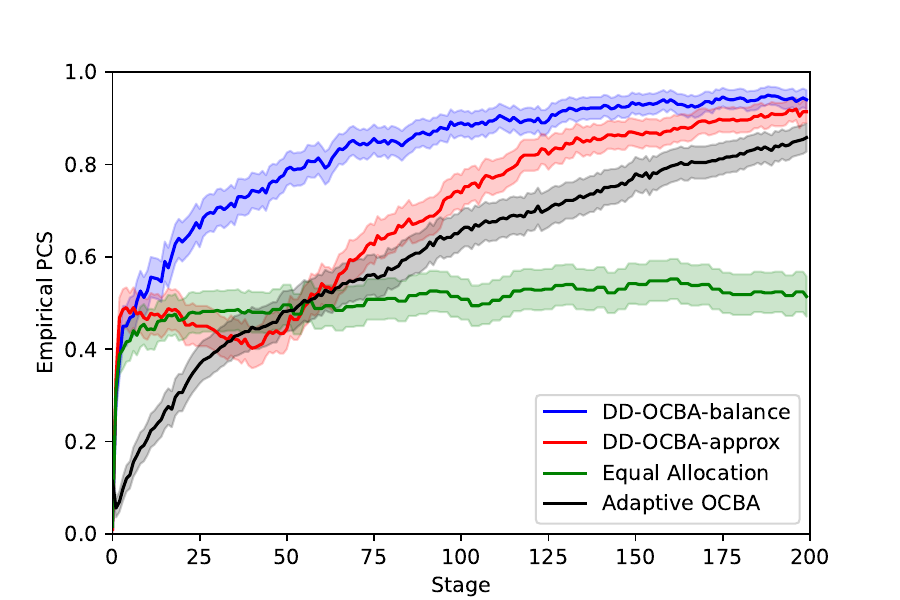}
        \caption{$m(t)=15$}
        \label{fig:quad_const_mt_15}
    \end{subfigure}
    \caption{Quadratic example with constant batch size and different choices of $m(t)$}
    \label{fig:quad_const_diff_mt}
\end{figure}

In Figure \ref{fig:quad_const_diff_mt}, We plot the empirical PCS of each algorithm at the end of each stage (with a total of 200) along with its $95\%$ confidence interval. The empirical PCS and the confidence interval are calculated by running $500$ macro-replications.  

The observations from Figure \ref{fig:quad_const_diff_mt} are summarized as follows:
\begin{enumerate}
    \item Across all configurations of $m(t)$, DD-OCBA-approx achieves a final empirical PCS around $90\%$
    and DD-OCBA-balance achieves a final empirical PCS around $95\%$, both of which surpass the performance of the other algorithms. The Equal Allocation method, in comparison, is markedly less efficient (achieving a final empirical PCS lower than $60\%$ in all 3 scenarios), primarily due to its overallocation of resources to design-input pairs characterized by incorrect input parameters. Furthermore, among the two highlighted methods, DD-OCBA-balance demonstrates superior performance over DD-OCBA-approx. This advantage stems from DD-OCBA-balance's direct approach to satisfy the original optimality conditions, in contrast to DD-OCBA-approx, which relies on an approximation of the "local balance" condition.
    \item In the comparison between DD-OCBA-approx and Adaptive OCBA, Adaptive OCBA always achieves a lower final empirical PCS aorund $80\%$ in all $3$ scenarios. Despite the different empirical performance by the two procedures, note that both allocation policies given by the two procedures will  converge to the limiting optimal policy defined in Theorem \ref{thm: optimality approx}. The convergence of Adaptive OCBA is guaranteed by Corollary \ref{cor: consistency continuous}, where Adpative OCBA can be regarded as a special case of DD-OCBA-approx-C with the number of subspace $D=1$. Nonetheless, Adaptive OCBA is outperformed by DD-OCBA-approx-C. This is attributed to how Adaptive OCBA aggregates simulation outputs, which vary across distributions, leading to an initial performance estimator with significant bias. As a result, at early stages, when estimation is inaccurate, the simulation output has a large bias.  This early bias diminishes over time with additional input data and simulations (Lemma \ref{lem: consistency continuous}). 
    On the contrary, DD-OCBA-approx maintains an unbiased estimator for each design-input pair, as simulations for these pairs are consistently conducted under the same input parameter, which implies more simulation outputs are always beneficial as it reduces the simulation error and does not introduce larger bias to the design-performance estimator. The only source of bias in DD-OCBA-approx comes from the posterior probability $\pi_j, j=1,\ldots,D$, which is independent of the simulation outputs.   This fundamental distinction renders DD-OCBA-approx more effective, particularly in scenarios with streaming input data and smaller batch sizes, where Adaptive OCBA's estimators are prone to larger biases.

\end{enumerate}

\subsection{Portfolio Optimization with Continuous Parameter Space and Random Batch Size}
In the previous example, we test the performance of DD-OCBA-approx and DD-OCBA-balance with finite parameter space. In this example, we test DD-OCBA-approx-C and DD-OCBA-balance-C, the generalized versions of the proposed procedures, 
on a more general problem of portfolio optimization where the input parameter space is continuous.
An investor invests a certain amount of capital in a riskless asset with interest rate $r$ and a risky asset, whose price per share at time $t$ is denoted as $S_t$.  $\{S_t\}$ is often assumed to follow a Geometric Brownian motion with initial price $S_0$, which admits the following expression for any fixed $t$:
\begin{equation} \label{eq: GBM}
    S_t = S_0 \exp\left[{\left(\theta^c-\frac{\sigma^2}{2}\right) t + \sigma B(t)}\right],
\end{equation}
where $\sigma$ is the volatility parameter, $\theta^c$ is the drift, and $\{B(t):t\ge 0\}$ is a standard Brownian motion. At time 0, the investor makes a one-time decision $x \in [0,1]$, which is the proportion of investment in the risky asset. Then, the total wealth at time t, denoted by $W_t$, is
\begin{equation*}
    W_t = x W_0 \exp\left[{\left(\theta^c-\frac{\sigma^2}{2}\right) t + \sigma B(t)}\right]+ (1-x)W_0\mathrm{e}^{r t}.
\end{equation*}
A risk-averse investor aims to maximize the mean-variance of the total asset after $T$ length of time  with a risk-averse parameter $\rho$, as follows:
\begin{equation*}
        F(x) = \mathbf{E} \left\{ x W_0 \exp\left[{\left(\theta^c-\frac{\sigma^2}{2}\right) T + \sigma B(T)}\right]+ (1-x)W_0\mathrm{e}^{r T}\right\} -\rho x^2W_0^2\mathrm{e}^{ 2\theta^c T  }(\mathrm{e}^{\sigma^2T}-1).
\end{equation*}
Here the variance term is calculated explicitly using the distribution of log-normal random variable, and the expectation term needs to be estimated. Furthermore, we assume the drift $\theta^c$ (also known as the risky return rate) is unknown but can be estimated with streaming data. Specifically, suppose we have a sequence of observations $\{S_{k\tau}\}_{k\in\mathbb{N}}$, where $S_{k\tau}$ is the price of the risky asset at time $k\tau$. From \eqref{eq: GBM}, we know 
$Z_k := \log{\left(\frac{S_{(k+1)\tau}}{S_{k\tau}}\right)} \sim \mathcal{N}\left((\theta^c - \frac{\sigma^2}{2})\tau, \tau\sigma^2\right)$ and $\{Z_k\}_{k\in \mathbb{N}}$ is serially independent. Hence, we can employ a Normal-Normal conjugate prior to estimate the unknown drift $\theta^c$.

As for implementation details, we set the initial wealth $W_0 = 1$, the interest rate of riskless asset $r=0.5$, volatility $\sigma = 1$, drift $\theta^c = 2.5$, risk-averse parameter $\rho = 0.1$. The candidate set $\mathcal{I}= \{0,1,\ldots,10\}$ with $i$th candidate being the solution $x_i = 0.1\cdot i$. We partition the entire parameter space (of $\theta^c$) $\Theta = \cup_{j=-10,10} \Theta_j$ where $\Theta_j = [\theta^c + 0.2 \cdot j - 0.1,\theta^c + 0.2 \cdot j + 0.1], -9\le j\le9$, $ \Theta_{-10} = (-\infty, 0.6]$
and $\Theta_{10} = [4.4,\infty)$. The prior distribution is set to be $\mathcal{N}(0,\tau \sigma^2)$.
 Furthermore, the stage-wise input data batch size and simulation budget is set to be random. Specifically, the input data batch size $m(t) = \tilde{m} * Z$ and the simulation budget $n(t) = \tilde{n} *Z$ with $Z$ being a random variable equally distributed among $\{1,2,3,4,5\}$ and $\Tilde{m}, \Tilde{n} \in \mathbb{N}_+$.
 We set the average stage-wise simulation budget $\Bar{n} = 2.5 \Tilde{n} = 50$ and varies the stage-wise input data batch size $\Bar{m} = 2.5 \Tilde{n}$ from $\{5,10,20\}$. In Figure \ref{fig:asset_random_diff_mt}, We plot the empirical PCS of each algorithm at the end of each stage (with a total of 200 stages) along with its $95\%$ confidence interval. The empirical PCS and the confidence interval are calculated by running $500$ macro-replications.

\subsubsection*{Experiment Result}
\begin{figure}[h]
\begin{minipage}[t]{0.3\textwidth}
        \subcaptionbox{$\Bar{m}=5$}{\includegraphics[height=1.5in]{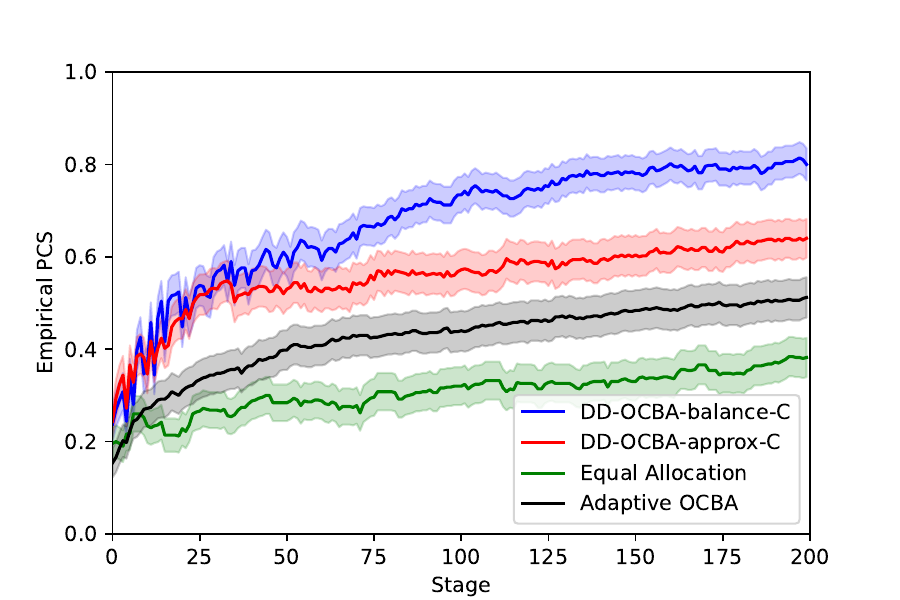}
        \label{fig:asset_random_mt_5}}
\end{minipage}
\hspace{0.1in}
\begin{minipage}[t]{0.3\textwidth}
        \subcaptionbox{$\Bar{m}=10$}{\includegraphics[height=1.5in]{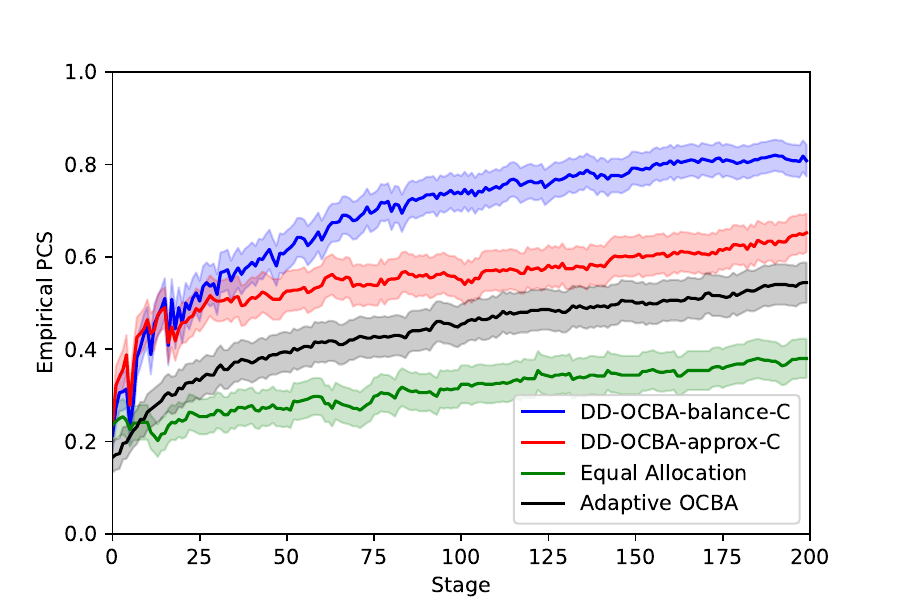}
        \label{fig:asset_random_mt_10}}
\end{minipage}
\hspace{0.1in}
\begin{minipage}[t]{0.3\textwidth}
        \subcaptionbox{$\Bar{m}=20$}{\includegraphics[height=1.5in]{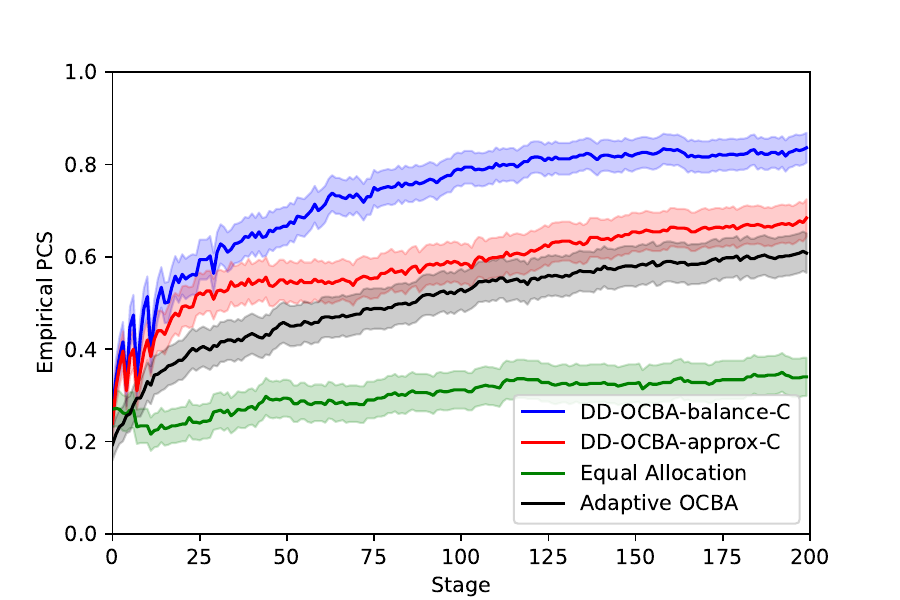}
        \label{fig:asset_random_mt_20}}
\end{minipage}
    \caption{Portfolio example with random batch size}
    \label{fig:asset_random_diff_mt}
\end{figure}
Similar conclusions can be drawn when comparing Figure \ref{fig:asset_random_diff_mt} with Figure \ref{fig:quad_const_diff_mt}. The observations are as follows:

\begin{enumerate}
\item In this portfolio example, which features a continuous parameter space and random batch sizes, the generalized extension procedures DD-OCBA-approx-C and DD-OCBA-balance-C continue to outperform other methods, achieving empirical PCS of approximately $60\%$ and $80\%$ within $200$ stages, respectively. In contrast, the Adaptive OCBA algorithm achieves a final empirical PCS of less than $50\%$ for $\Bar{m} = 5, 10$ and less than $60\%$ for $\Bar{m} = 20$. The equal allocation procedure performs the worst, with a final empirical PCS below $40\%$ in all three scenarios.

\item Unlike in the quadratic example where DD-OCBA-approx generated i.i.d. unbiased simulation outputs of the design-input performance, in this example of continuous parameter space, both DD-OCBA-approx-C and Adaptive OCBA suffer from bias due to variations in the posterior distribution. However, DD-OCBA-approx-C still outperforms Adaptive OCBA. This superiority stems from its partitioning of the entire parameter space, which helps reduce bias. Specifically, in Adaptive OCBA, when the posterior distribution is roughly estimated, the simulation outputs for a fixed design are generated under input parameters that may deviate significantly from the true value $\theta^c$, introducing substantial bias. Conversely, in DD-OCBA-approx-C, the simulation outputs for a fixed design $i$ under $\Theta_0$ (which includes the true parameter $\theta^c$) are generated under input parameters $\theta \in \Theta_0$. Since $\Theta_0$ represents a ``small" neighborhood around $\theta^c$ with a radius of $0.1$, as opposed to the entire (unbounded) parameter space, the bias in the simulation outputs under $(i,\Theta_0)$ generated in DD-OCBA-approx-C remains relatively small, even with a roughly estimated posterior distribution.
\end{enumerate} 

\section{Conclusion and Future Work} \label{sec:conclusion}
In this paper we consider a fixed budget ranking and selection (R\&S) problem, where the common input distribution across designs is unknown but can be estimated with streaming input data that come sequentially in time.
We initially assume a finite parametrization of the input distribution and utilize a Bayesian approach for estimation, which is updated at the beginning of each stage as batched data arrives.  Subsequently, the stage-wise computing budget is allocated to facilitate new simulations for assessing design performances. We apply the large deviations theory to obtain the optimal stage-wise budget allocation policy for design-input pairs. Based on the optimality equations, we design two fully sequential algorithms that achieve consistency (i.e., select the best design with probability 1 as times go to infinity) and asymptotic optimality (i.e., converge to the optimal budget allocation policy under the true input distribution). 
 We also extend our approach to accommodate a continuous input parameter space, while maintaining provable statistical validity.  Our numerical experiments demonstrate  superior performance of the proposed procedures over the equal allocation rule and an extension of OCBA when dealing with unknown input distributions with streaming input data.

For future research, it would be intriguing to explore methods for adaptively adjusting the finite input parameter space as additional input data becomes available. This adjustment could help correct any potential discretization errors in the original parameter space. One possible approach could be based on criteria or indices like the Bayesian Information Criterion (BIC) or the Davies-Bouldin Index, which would guide the refinement process.
Another promising research direction is to incorporate a non-parametric Bayesian framework. This approach could provide a flexible way to address and correct model mismatches arising from the parametric assumptions initially made about the input distribution. By using non-parametric methods, researchers could potentially capture a wider variety of distribution shapes and complexities, thereby enhancing the robustness and accuracy of the modeling process.  

\section*{ACKNOWLEDGMENT}
The authors gratefully acknowledge the support by the Air Force Office of Scientific Research under Grant FA9550-22-1-0244, the National Science Foundation under Grant NSF-DMS2053489, and AI Institute for Advances in Optimization (AI4OPT).

\section*{BIOGRAPHY}

\noindent {\bf YUHAO WANG} is a Ph.D. candidate at the H. Milton Stewart School of Industrial and Systems Engineering at Georgia Institute of Technology. He received his B.S. degree from the Department of Mathematics at Nanjing University, China, in 2021. His research interests include simulation, stochastic optimization, and reinforcement learning. \\

\noindent {\bf ENLU ZHOU} is a Professor in the H. Milton Stewart School of Industrial and Systems Engineering at Georgia Institute of Technology. She received the B.S. degree with highest honors in electrical engineering from Chu Kochen Honors College, Zhejiang University, China, in 2004, and the Ph.D. degree in electrical engineering from the University of Maryland, College Park, in 2009. Prior to joining Georgia Tech in 2013, she was an assistant professor in the Department of Industrial and Enterprise Systems Engineering at the University of Illinois Urbana-Champaign from 2009 to 2013. She is a recipient of the AFOSR Young Investigator award in 2012, NSF CAREER award in 2014, the INFORMS Outstanding Simulation Publication award in 2020, and the Best Theoretical Paper award at the Winter Simulation Conference twice in 2009 and 2022. Her research interests lie in theory, methods, and applications of simulation optimization, stochastic optimization, and stochastic control.

\bibliographystyle{plainnat}  
\bibliography{arxiv_main}  

\newcommand{\tl}{{t_\ell}}
\newcommand{\ptl}{{(t_\ell)}}
\newcommand{\pl}{{(\ell)}}
\newcommand{\jc}{{j^c}}
\newcommand{\Opar}[1]{O\left({#1}\right)}
\newcommand{\pars}[1]{\left({#1}\right)}
\newcommand{\logll}{\sqrt{\frac{\log\log\ell}{\ell}}}
\newcommand{\llogl}{\sqrt{\ell \log\log\ell}}
\newcommand{\logrr}{\sqrt{\frac{\log\log r}{r}}}
\newcommand{\rlogr}{\sqrt{r \log\log r}}
\appendix
\section*{Technical Proof}
Recall $\ell$ is the iteration counter. Denote by $t_\ell$ the stage where $\ell$th simulation is run. That is 
\begin{equation}\label{eq: def t_l}
    t_\ell : = \arg \max_t \{ t:  N(t-1) < \ell\},
\end{equation}
where $N(t) = \sum_{\tau=1}^t n(t)$ is total simulation budget up to stage $t$. 
\subsection{ Proof of Lemma \ref{lem:rate function}}
\textbf {Proof.}
 Lemma \ref{lem:rate function}.1 is easily seen from \eqref{eq:rate formula}. To prove  Lemma \ref{lem:rate function}.2, it suffices to show the concavity of the function for $x>0$ with form $ f(x) =1/( \sum_{i=1}^n \frac{a_i}{x_i}),$  where $a_i > 0$ for $i=1,2,\dots,n$. We prove the concavity of the multivariate function by proving the concavity along all lines. For any $y \in \mathbb{R}^n$, let $g(t) = f(x+ty)$ where $t \in \mathbb{R}$ such that $x+ty > 0$. We have 
 \begin{equation*}
     g^{\prime\prime}(t) = \frac{2}{(\sum_{i=1}^n\frac{a_i}{x_i+ty_i})^3  } \left\{ \left[ \sum_{i=1}^n \frac{a_iy_i}{(x_i+ty_i)^2}  \right]^2 - \sum_{i=1}^n \frac{a_i y_i^2}{(x_i+ty_i)^3} \sum_{i=1}^n \frac{a_i}{x_i+ty_i} \right\}  \le 0,
 \end{equation*}
 where the inequality uses the Cauchy inequality. Hence, $f$ is concave in $x > 0$. \hfill $\blacksquare$

\subsection{Proof of Theorem \ref{thm:optimal allocation}}
\textbf{Proof.}
 We first show the existence of $\alpha$. The existence follows from the continuity of $G_i$ with respect to $\alpha \in \Delta^{KB-1}$, where $\Delta^n$ denotes the $n$-dimensional simplex. Furthermore, by Lemma \ref{lem:rate function}.1, $G_i$ is strictly increasing in $\alpha_{1,j'}$ and $\alpha_{i,j'}$. Since $\alpha_{i,j} = \frac{1}{KB} \ \forall i,j$ is a feasible solution and the corresponding objective value is strictly positive, the optimal solution $\alpha$ must satisfy $\alpha_{i,j} > 0$ for all $i,j$.
 
 Now we show the necessity of the three optimality conditions. By Lemma \ref{lem:rate function}.2, the optimization problem (\ref{eq:opt LDR}) is a concave maximization problem, and therefore the KKT conditions are both sufficient and necessary for the optimality. With $\alpha$ strictly positive, the KKT conditions can be written as
 \begin{flalign}
 1 - \sum_{i \neq b} \lambda_i &= 0, \label{eq:KKT z}\\
\lambda_i \frac{\partial G_i}{\partial \alpha_{i,j}}(\alpha_b,\alpha_i) &= \gamma \quad i \neq b, \  1\le j \le D,\label{eq:KKT derivative}\\
\sum_{i \neq b}  \lambda_i \frac{\partial G_i}{\partial \alpha_{b,j}}(\alpha_b,\alpha_i) &= \gamma \quad 1\le j\le D, \label{eq:KKT total} \\
\lambda_i(G_i(\alpha_b,\alpha_i)-z)&=0 \quad i \neq b, \label{eq:KKT local}
\end{flalign}
for some $\gamma$ and $\lambda_i \ge 0, \ i \neq b$. By (\ref{eq:KKT z}) there exists at least one $i_0$ such that $\lambda_{i_0} > 0$. Then since $G_i$ is increasing in $\alpha_{i,j}$, we have $\frac{\partial G_{i_0}}{\partial \alpha_{i_0,j}}(\alpha_b,\alpha_{i_0}) > 0 $. This implies $\gamma > 0$ by (\ref{eq:KKT derivative}). Hence, we must have $\lambda_i > 0$ for all $i \neq b$. Then we have $\frac{\partial\mathbf{G}_i(\alpha_b,\alpha_i)}{\partial\alpha_{i,j}} = \frac{\gamma}{\lambda_i}, \  1\le j\le D,\ i \neq b$, which proves (\ref{thm:input balance}). Since $\lambda_i > 0$,
$G_i(\alpha_b,\alpha_i)=z,  \ i \neq b$ by (\ref{eq:KKT local}). Hence, (\ref{thm:local balance}) holds. To see why (\ref{thm:total balance}) holds, solving for $\lambda_i =\frac{\gamma}{\frac{\partial G_i}{\partial \alpha_{i,j}}(\alpha_b,\alpha_i)}$ in (\ref{eq:KKT derivative}) and substituting $\lambda_i$ in (\ref{eq:KKT total}), we get the desired result.

For sufficiency, first let $\lambda_i = \frac{1}{\partial G_i(\alpha_b,\alpha_i) / \partial \alpha_{i,j}}/(\sum\limits_{i \neq b} \frac{1}{\partial G_k(\alpha_b,\alpha_k) / \partial \alpha_{k,j}})$ for $i \neq b$. Notice that $\lambda_i > 0$ and does not depend on the choice of $j$ by (\ref{thm:input balance}). Moreover, $\{\lambda_i\}_{i \neq b}$ satisfy condition (\ref{eq:KKT z}). Further let $\gamma = (\sum\limits_{i \neq b} \frac{1}{\partial G_k(\alpha_b,\alpha_k) / \partial \alpha_{kj}})^{-1}$, which is also independent of $j$. We can easily verify that both (\ref{eq:KKT derivative}) and (\ref{eq:KKT total}) hold. (\ref{eq:KKT z}) also holds by setting $z = G_i(\alpha_b,\alpha_i)$, which is independent of $i$ by (\ref{thm:local balance}).

Now we are only left to show the uniqueness of $\alpha$. First notice that from (\ref{eq:input balance})  and (\ref{eq:total balance}), we have
$$\frac{\alpha_{b,j}}{\sigma_{b,j}{\pi_j}} = \sqrt{\sum\limits_{i \neq b} (\frac{\alpha_{i,j}}{\sigma_{i,j}{\pi_j}})^2}= \sqrt{\sum\limits_{i \neq b} (\frac{\alpha_{i,j'}}{\sigma_{i,j'}\pi_{j'}})^2} =\frac{\alpha_{1,j'}}{\sigma_{1,j'}\pi_{j'}} \quad    1 \le j<j' \le D. $$ 
Letting $\beta_i = \frac{\alpha_{i,j}}{\sigma_{i,j} {\pi_j}}$ which is independent of $j$, we can write $\alpha_{i,j} = {\pi_j}\sigma_{i,j} \beta_i$ for all $i=1,2,\cdots,K$ and $j=1,2,\cdots,D$. Since $\alpha$ and $\beta = (\beta_b,\cdots,\beta_K)$ are bijective, it is sufficient to show the uniqueness of $\beta$. Plugging $\alpha_{i,j}$ into (\ref{eq:local balance}) and (\ref{eq:total balance}), we have
\begin{equation*}
    \frac{({\mu}_b - {\mu}_i)^2}{ \frac{\sum_{j=1}^D\sigma_{b,j}{{\pi_j}}}{\beta_{b}}+\frac{\sum_{j=1}^D\sigma_{i,j}{{\pi_j}}}{\beta_{i}}} = \frac{({\mu}_b - {\mu}_{i'})^2}{ \frac{\sum_{j=1}^D\sigma_{b,j}{{\pi_j}}}{\beta_{b}}+\frac{\sum_{j=1}^D\sigma_{i',j}{{\pi_j}}}{\beta_{i'}}} \quad i \neq i' \neq b
\end{equation*} 
with
$
\beta_b^2 = \sum_{i \neq b} \beta_i^2.
$
 Let $\eta = \frac{\beta}{\beta_b}$. Then $\eta$ satisfies
 \begin{equation}\label{eq:eta balance}
\frac{({\mu}_b - {\mu}_i)^2}{ \sum_{j=1}^D\sigma_{b,j}{{\pi_j}}+\frac{\sum_{j=1}^D\sigma_{i,j}{{\pi_j}}}{\eta_{i}}} = \frac{({\mu}_b - {\mu}_{i'})^2}{ \sum_{j=1}^D\sigma_{b,j}{{\pi_j}}+\frac{\sum_{j=1}^D\sigma_{i',j}{{\pi_j}}}{\eta_{i'}}} \quad i \neq i' \neq b      
\end{equation}
with
$
    1 = \sum_{i \neq b} \eta_i^2.
$
If there exists $\eta' \neq \eta$ satisfying these two conditions, then there must be $i \neq k \neq 1$ such that $\eta_i < \eta'_i$ and $\eta_k > \eta'_k$. Then, we have 
\\$\frac{({\mu}_b - {\mu}_i)^2}{ \sum_{j=1}^D\sigma_{b,j}{{\pi_j}}+\frac{\sum_{j=1}^D\sigma_{i,j}{{\pi_j}}}{\eta_{i}}} <
\frac{({\mu}_b - {\mu}_i)^2}{ \sum_{j=1}^D\sigma_{b,j}{{\pi_j}}+\frac{\sum_{j=1}^D\sigma_{i,j}{{\pi_j}}}{\eta'_{i}}}=
\frac{({\mu}_b - {\mu}_{k})^2}{ \sum_{j=1}^D\sigma_{b,j}{{\pi_j}}+\frac{\sum_{j=1}^D\sigma_{k,j}{{\pi_j}}}{\eta'_{k}}}<
\frac{({\mu}_b - {\mu}_{k})^2}{ \sum_{j=1}^D\sigma_{b,j}{{\pi_j}}+\frac{\sum_{j=1}^D\sigma_{k,j}{{\pi_j}}}{\eta_{k}}}$,  
which contradicts \eqref{eq:eta balance}. Hence, $\eta$ is unique, which implies $\beta = C*\eta$ for some constant $C$. Then if there exists $\beta' \neq \beta$ and both are optimal, we have $\beta >(<) \beta'$. This implies the corresponding $\alpha >(<) \alpha'$, which contradicts $\sum_{i,j}\alpha_{i,j} = \sum_{i,j}\alpha'_{i,j} =1$.    \hfill $\blacksquare$

\subsection{Proof of Theorem \ref{thm:consistency}}
\textbf{Proof}
 Denote by $\omega$ a sample path of one simulation process which contains all the simulation outputs and input data observation. 
 To prove for either DD-OCBA-approx or DD-OCBA-balance, it suffices to show $N_{i,j^c}^{(\ell)} \rightarrow \infty, $ as $ \ell \rightarrow \infty$ almost surely for all $i$. To see this, recall $X_{i,j^c}^{(r)}$ is the $r$th simulation output for $(i,j^c)$ pair (the iteration at which $X_{i,j^c}^{(r)}$ is simulated is still random). By Assumption \ref{assump:output}.3, $\{X_{i,j^c}^{(r)}\}$ are i.i.d. for $r=1,2, \ldots $. Hence $\frac{1}{n}\sum_{r=1}^n X_{i,j^c}^{(r)}$ converges to $\mu_{i,j^c}$ almost surely by LLN. Take any sample path $\omega$, then almost surely $N_{i,j^c}^{(\ell)} |\omega \rightarrow \infty$ as $\ell \rightarrow \infty$. Under the same sample path $\frac{1}{n}\sum_{r=1}^n X_{i,j^c}^{(r)} | \omega \rightarrow \mu_{i,j^c}$ as $n \rightarrow \infty$. Hence,  $\frac{1}{N_{i,j^c}^{(\ell)}|\omega}\sum_{r=1}^{N_{i,j}^{(\ell)}|\omega} X_{i,j^c}^{(r)} | \omega \rightarrow \mu_{i,j^c}$ as $\ell \rightarrow \infty$.   Furthermore,  $\hat{\sigma}_{i,j^c}^{(\ell)} \rightarrow \sigma_{i,j^c} \ \text{ almost surely }$ by the same argument and $ \pi_{j^c}^{(t_\ell)} \rightarrow 1 \ \text{ almost surely }$ by Doob's consistency theorem. Therefore, we have $\hat{\mu}_{i}^{(\ell)} \rightarrow  \mu^ci$ as $\ell \rightarrow \infty$ almost surely. We fix a sample path $\omega$ in the following proof of showing $N_{i,j^c}^{(\ell)} \rightarrow \infty, $ as $ \ell \rightarrow \infty$ almost surely. Denote by $A = \{ i | N_{i,j^c}^{(\ell)} \rightarrow \infty\}$.  
 Notice
$\hat{\mu}_{i,j}^\pl, \hat{\sigma}^{(\ell)}_{i,j}$ will converge (to some random variable) almost surely, no matter whether $N_{i,j}^{(\ell)}$ tend to infinity. This is because if $N_{i,j}^{(\ell)}$ is at most finite, then $\hat{\mu}_{i,j}^{(\ell)}$ and $\hat{\sigma}^{(\ell)}_{i,j}$ will remain unchanged after finite iterations. Since $\pi_j^\tl$ converges almost surely, $\hat{\mu}_{i}^{(\ell)}$ will also  converge almost surely. This implies $\hat{b}^\pl$ also converge (to some random design index) almost surely. We use $\hat{\mu}_{i},\hat{\sigma}_{i,j},\hat{b}$ to denote the limit of $\hat{\mu}_{i,j}^\pl, \hat{\sigma}^{(\ell)}_{i,j}$ and $\hat{b}^\pl$, respectively.
\\


\textbf{Proof of Theorem \ref{thm:consistency}.1.} 
 Denote by $N^{(\ell)} = \sum_{i,j} N_{i,j}^{(\ell)}$. Then there exists an allocation policy $\{ \Tilde{\alpha}_{i,j} \}$  satisfying $\lim_{\ell\rightarrow \infty} \frac{\hat{N}_{i,j}^{(\ell)}}{N^{(\ell)}}  = \Tilde{\alpha}_{i,j}$. Furthermore, $\Tilde{\alpha}_{i,j} = 0$ if $j\neq j^c$ and $\Tilde{\alpha}_{i,j^c} > 0$ if $j=j^c$. This is because
for $i,i'\neq \hat{b}$, and $j\neq j^c$,
$$  \frac{\hat{N}_{i,j}^\pl}{\hat{N}_{i^\prime, j^c}} = \frac{{\pi_j^\tl}\hat{\sigma}^\pl_{i,j}\sum_{k=1}^D\hat{\sigma}^\pl_{i,k}{\pi_{k}^\tl}/(\hat{\mu}^\pl_{\hat{b}} - \hat{\mu}^\pl_i)^2}{{\pi^\tl_{j^c}}{\hat{\sigma}}^\pl_{i^\prime, j^c} \sum_{k=1}^D\hat{\sigma}^\pl_{i',k}{{\pi^\tl_{k}}} /(\hat{\mu}^\pl_{\hat{b}} - \hat{\mu}^\pl_{i'})^2} \rightarrow 0,$$
by the convergence of $\hat{\mu}_{i,j}^\pl,\hat{\sigma}^\pl_{i,j}$ and the fact that $\pi_j^\tl \rightarrow 0, \pi_{j^c}^\tl \rightarrow 1$ by Doob's consistency theorem. This implies ${\Tilde{\alpha}_{i,j}} = 0$. Furthermore, since $\hat{N}_{\hat{b},j}^\pl = \hat{\sigma}_{\hat{b},j}^\pl \sqrt{\sum_{i\neq \hat{b}} \left(\frac{\hat{N}_{i,j}^\pl}{\hat{\sigma}_{i,j}^\pl}\right)^2}$, we know $\Tilde{\alpha}_{\hat{b},j} = 0$. Moreover, we have
$$ \frac{\Tilde{\alpha}_{i,j^c}}{\Tilde{\alpha}_{i',j^c}}= \lim_{\ell\rightarrow\infty} \frac{\hat{N}_{i,j^c}^\pl}{\hat{N}_{i^\prime, j^c}} \rightarrow \frac{\hat{\sigma}^2_{i,j^c} \left(\hat{\mu}_{i'} -\hat{\mu}_{\hat{b}}\right)^2}{\hat{\sigma}^2_{i',j^c} \left(\hat{\mu}_{i} -\hat{\mu}_{\hat{b}}\right)^2} , \Tilde{\alpha}_{\hat{b},j^c} = \hat{\sigma}_{b,j^c} \sqrt{\sum_{i\neq \hat{b}} \left(\frac{\Tilde{\alpha}_{i,j^c}}{\hat{\sigma}_{i,j^c}}\right)^2} .$$
This implies $\Tilde{\alpha}_{i,j^c} >0, \forall i \in \mathcal{I}$.
Suppose there exists $i_0 \not \in A$. Let $\varepsilon = \frac{\tilde{\alpha}_{i_0,j^c}}{KD}$, there exists $L$ large enough, such that $\forall \ell \ge L, \forall i,j$, $|\frac{\hat{N}_{i,j}^{(\ell)}}{N^{(\ell)}} - \tilde{\alpha}_{i,j}| <\frac{\varepsilon}{2}$. Since $N_{i_0,j^c}^{(\ell)} $ remains unchanged after some iteration $\ell_0$, 
we have $\frac{\hat{N}_{i_0,j^c}^{(\ell)} - N_{i_0,j^c}^{(\ell)}}{N^{(\ell)}} > 0$ for $\ell>L$ and $\ell$ large enough. 
Notice for any $i,j$, if $N_{i,j}^{(\ell)} > (\Tilde{\alpha}_{i,j}+\frac{\varepsilon}{2}) N^{(\ell)} > \hat{N}_{i,j}^{(\ell)} $, then $(i,j)$ will not be simulated since $ \hat{N}_{i,j}^{(\ell)} - N_{i,j}^{(\ell)} < 0 < \hat{N}_{i_0,j^c}^{(\ell)} - N_{i_0,j^c}^{(\ell)} $. Hence, we must have $N_{i,j}^{(\ell)} \le (\tilde{\alpha}_{i,j}+\frac{\varepsilon}{2})N^{(\ell)} +1 \le (\tilde{\alpha}_{i,j}+\varepsilon)N^{(\ell)}$ for $\ell>L$ and $\ell$ large enough. Then, we have
\begin{equation} \label{eq:approx consistency proof eq1}
    N^{(\ell)} = \sum_{i,j} N^{(\ell)}_{i,j} = N^{(\ell)}_{i_0,j^c} + \sum_{(i,j)\neq (i_0,j_0)} N_{i,j}^{(\ell)} \le N^{(\ell)}_{i_0,j^c} + N^{(\ell)}\sum_{(i,j)\neq (i_0,j_0)} (\tilde{\alpha}_{i,j}+\varepsilon) 
\end{equation}
Divided by $N^{(\ell)}$ on both sides and let $\ell\rightarrow\infty$, \eqref{eq:approx consistency proof eq1} implies
\begin{equation*}
    1 \le \sum_{(i,j)\neq (i_0,j^c)} \tilde{\alpha}_{i,j}+(KD-1)\varepsilon = \sum_{(i,j)\neq (i_0,j_0)} \tilde{\alpha}_{i,j}+ \frac{KD-1}{KD}\tilde{\alpha}_{i_0,j^c} <1,
\end{equation*}
a contradiction to $i_0\not\in A$. The proof is complete.
\\ 
\textbf{ Proof of Theorem \ref{thm:consistency}.2.}  It suffices to prove $\lim_{\ell\rightarrow\infty} N_{i,j^c}^{(\ell)} \rightarrow \infty$ almost surely. We first prove $A\neq \O$ by contradiction. Suppose $i\not\in A, \forall i \in \mathcal{I}$. We know there exists $i_0,j\neq j^c$, such that $N_{i_0,j}^\pl \rightarrow \infty$. This implies $\frac{N^\pl_{i_0,j}}{\hat{\sigma}_{i_0,j} \pi_j^\tl} \rightarrow \infty$, which can happens only if $\frac{N^\pl_{i_0,j'}}{\hat{\sigma}_{i_0,j'} \pi_{j'}^\tl} \rightarrow \infty, \forall j'$ otherwise $(i,j)$ cannot be sampled after some iteration. Hence, we have $\frac{N^\pl_{i_0,j^c}}{\hat{\sigma}_{i_0,j^c}\pi_{j^c}} \rightarrow \infty$, which further implies $N^\pl_{i_0,j^c} \rightarrow \infty$ since $\pi_{j^c} \rightarrow 1$. Hence we prove $A\neq \O$.

Next, we prove there exists $j$, $N_{\hat{b},j}^\pl \rightarrow \infty$. Suppose not, then for all $j$, 
$$ \left(\dfrac{N_{\hat{b},j}^{(\ell)}}{\hat{\sigma}^{(\ell)}_{\hat{b},j}} \right)^2 -\sum_{i\neq \hat{b}} \left(\dfrac{N_{i,j}^{(\ell)}}{\hat{\sigma}^{(\ell)}_{i,j}} \right)^2 \le  \left(\dfrac{N_{\hat{b},j}^{(\ell)}}{\hat{\sigma}^{(\ell)}_{\hat{b},j}} \right)^2$$
is upper bounded by a constant. Furthermore we know there exists $i' \in A$, hence 
$ \left(\dfrac{N_{\hat{b},j^c}^{(\ell)}}{\hat{\sigma}^{(\ell)}_{\hat{b},j^c}} \right)^2-\sum_{i\neq \hat{b}} \left(\dfrac{N_{i,j^c}^{(\ell)}}{\hat{\sigma}^{(\ell)}_{i,j^c}} \right)^2 \rightarrow -\infty $. This implies $i'$ cannot be simulated after some iteration, a contradiction to $N_{i',j^c}^\pl \rightarrow \infty$. Hence, we know there exists $j$, $N_{\hat{b},j}^\pl \rightarrow \infty $. This implies $\frac{N_{\hat{b},j}^\pl}{\pi_j^\tl, \hat{\sigma}_{\hat{b},j}^\pl} \rightarrow \infty$, which further implies $ N_{\hat{b},j^c}^\pl \rightarrow \infty$ wit the similar argument for proving $A \neq \O$.  Hence, we prove $\hat{b} \in A$.  

Next, we prove there exists $i'\neq \hat{b}, j'$, such that $N_{i',j'}^\pl \rightarrow \infty$. Suppose not, then 
$$ \left(\dfrac{N_{\hat{b},j}^{(\ell)}}{\hat{\sigma}^{(\ell)}_{\hat{b},j}} \right)^2 -\sum_{i\neq \hat{b}} \left(\dfrac{N_{i,j}^{(\ell)}}{\hat{\sigma}^{(\ell)}_{i,j}} \right)^2 \ge  -\sum_{i\neq \hat{b}} \left(\dfrac{N_{i,j}^{(\ell)}}{\hat{\sigma}^{(\ell)}_{i,j}} \right)^2$$
is lower bounded by some constant and 
$ \left(\dfrac{N_{\hat{b},j^c}^{(\ell)}}{\hat{\sigma}^{(\ell)}_{\hat{b},j^c}} \right)^2-\sum_{i\neq \hat{b}} \left(\dfrac{N_{i,j^c}^{(\ell)}}{\hat{\sigma}^{(\ell)}_{i,j^c}} \right)^2 \rightarrow \infty$. This implies $(\hat{b},j), \forall j$ cannot be simulated after some iteration, which contradicts to $b\in A$. Hence, there exists $i',j'$, such that $N_{i',j'} ^\pl \rightarrow \infty$. Again, this implies $N_{i',j^c}^\pl \rightarrow \infty$.

Finally we prove for all $i\neq \hat{b}, i\in A$. Since $N_{\hat{b},j^c}^\pl \rightarrow \infty, N_{i',j^c}^\pl \rightarrow \infty$, $\pi_j^\tl \rightarrow0, j\neq j^c$, we know 
$$
    \dfrac{(\hat{\mu}^{(\ell)}_{\hat{b}^{(\ell)}} - \hat{\mu}_{i'}^{(\ell)})^2}{\sum\limits_{j=1}^D \frac{(\hat{\sigma}^{(\ell)}_{\hat{b}^{(\ell)},j})^2(\pi_j^\tl)^2}{N_{\hat{b}^{(\ell)},j}} + \sum\limits_{j=1}^D\frac{(\hat{\sigma}^{(\ell)}_{i',j})^2(\pi_j^\tl)^2}{N_{i_0,j}^{(\ell)}}} \rightarrow \infty
$$, we must have $\forall i\neq \hat{b}$, $
    \dfrac{(\hat{\mu}^{(\ell)}_{\hat{b}^{(\ell)}} - \hat{\mu}_{i}^{(\ell)})^2}{\sum\limits_{j=1}^D \frac{(\hat{\sigma}^{(\ell)}_{\hat{b}^{(\ell)},j})^2(\pi_j^\tl)^2}{N_{\hat{b}^{(\ell)},j}} + \sum\limits_{j=1}^D\frac{(\hat{\sigma}^{(\ell)}_{i,j})^2(\pi_j^\tl)^2}{N_{i,j}^{(\ell)}}} \rightarrow \infty
$. This implies that $N_{i,j^c}^{(\ell)} \rightarrow \infty,\ \forall i \neq \hat{b}$. So far the proof is complete. \hfill $\blacksquare$

\subsection{Proof of Theorem \ref{thm: optimality approx}}
\textbf{Proof.}
Again denote by $\omega$ any sample path of one simulation process and   we fix a sample path $\omega$. 
 Since $\pi^\tl_j \rightarrow 0, j\neq j^c$, $\pi^\tl_{j^c} \rightarrow 0$ and $N_{i,j^c}^\pl \rightarrow \infty$ as $\ell \rightarrow \infty$, we have $\hat{\alpha}_{i,j}^{(\ell)} \rightarrow 0$ for $i\in\mathcal{I}, j\neq j^c$ and $\hat{\alpha}_{i,j^c}^\pl \rightarrow \alpha_{i,j^c}^*, i \in \mathcal{I}$ as $\ell \rightarrow \infty$. For simplicity, let $\alpha_{i,j}^* = 0$ for $j\neq j^c$. Let $A= \{(i,j): N_{i,j}^\pl \rightarrow \infty\}$. Clearly $(i,j^c) \in A, \forall i\in\mathcal{I}$ by the consistency result. If $(i,j) \not\in A,$ then we know $\alpha_{i,j}^\pl \rightarrow \alpha_{i,j}^* = 0$. Hence, $\forall \varepsilon > 0,\  \exists \Tilde{L}$ such that for all $\ell\ge\Tilde{L}$, $|\hat{\alpha}_{i,j}^{(\ell)} -\alpha_{i,j}^*| < \varepsilon \  \forall i,j$, $\alpha_{i,j}^\pl \le \varepsilon \ \forall (i,j) \not \in A$ and $\frac{1}{N^{(\ell)}} < \varepsilon $.
 Let $\Tilde{L}_{i,j} = \min\{\ell>\Tilde{L}: N_{i,j}^{(\ell)} = N_{i,j}^{(\Tilde{L})} + 1\} \ \forall (i,j) \in A$, the first time $(i,j)$ is sampled after $\Tilde{L}$. Let $L = \max_{(i,j) \in A} \Tilde{L}_{i,j} < \infty$ by the definition of $A$. Then for any $\ell>L$, let $D_l = \{(i,j): \hat{\alpha}_{i,j}^{(\ell)} - \frac{N_{i,j}^{(\ell)}}{N^{(\ell)}} < 0\}$. Then if $(i,j) \in D_l \cap A$, let  $L_{i,j} = \max\{s<\ell:N_{i,j}^{(s)} = N_{i,j}^{(\ell)}-1\}$. Then we have $L_{i,j} \ge \Tilde{L}$ by the definition of $L$. Furthermore,
\begin{equation} \label{eq:L_{i,j}}
    \hat{\alpha}_{i,j}^{(\ell)} - \frac{N_{i,j}^{(\ell)}}{N^{(\ell)}} \ge \hat{\alpha}_{i,j}^{(\ell)} - \frac{N_{i,j}^{(L_{i,j} + 1)}}{N^{(L_{i,j} + 1)}} = \underbrace{\left[\hat{\alpha}_{i,j}^{(L_{i,j})} - \frac{N_{i,j}^{(L_{i,j})}}{N^{(L_{i,j})}}\right]}_{\textstyle \mathstrut =E_1} + \underbrace{\left[\hat{\alpha}_{i,j}^{(\ell)} - \hat{\alpha}^{(L_{i,j})}_{i,j} \right]}_{\textstyle \mathstrut =E_2} + \underbrace{\left[ \frac{N_{i,j}^{(L_{i,j})}}{N^{(L_{i,j})}} -   \frac{N_{i,j}^{(L_{i,j} + 1)}}{N^{(L_{i,j} + 1)}}   \right]}_{\textstyle \mathstrut =E_3},  
\end{equation}
where the first inequality follows from $N_{i,j}^{(\ell)} = N_{i,j}^{(L_{i,j}+1)}$ and $N^{(L_{i,j}+1)} < N^{(\ell)}$ by the definition of $L_{i,j}$. Since $(i,j)$ is sampled at $L_{i,j} $, we must have $E_1\ge 0$. Further since $\ell,L_{i,j} \ge \Tilde{L}$, we have $E_2 \ge -|\alpha_{i,j}^* - \hat{\alpha}_{i,j}^{(\ell)}| - |\alpha_{i,j}^*-\hat{\alpha}^{(L_{i,j})}_{i,j}| > -2\varepsilon $ and $E_3 = \frac{N_{i,j}^{(L_{i,j})}}{N^{(L_{i,j})}} -   \frac{N_{i,j}^{(L_{i,j})+1}}{N^{(L_{i,j}+1)}} \ge  \frac{N_{i,j}^{(L_{i,j})}}{N^{(L_{i,j})}} -   \frac{N_{i,j}^{(L_{i,j})}+1}{N^{(L_{i,j})}+1} = \frac{N_{i,j}^{(L_{i,j})}-N^{(L_{i,j})}}{N^{(L_{i,j})}(N^{(L_{i,j})}+1)} > \frac{-N^{(L_{i,j})}}{N^{(L_{i,j})}(N^{(L_{i,j})}+1)} = \frac{-1}{N^{(L_{i,j})}+1} > -\varepsilon  $. Hence, we have $ \hat{\alpha}_{i,j}^{(\ell)} - \frac{N_{i,j}^{(\ell)}}{N^{(\ell)}} \ge -3\varepsilon$ for $(i,j) \in D_l\cap A$. For $(i,j)\in D_l \backslash A$, we also have $\hat{\alpha}^\pl_{i,j} - \frac{N_{i,j}^\pl}{N^\pl} \ge - \frac{N_{i,j}^\pl}{N^\pl} \ge -\varepsilon > -3\varepsilon$. This implies for all $(i,j) \in D_l$,   $ \hat{\alpha}_{i,j}^{(\ell)} - \frac{N_{i,j}^{(\ell)}}{N^{(\ell)}} \ge -3\varepsilon$.

As a result,
$
    0 = \sum_{i,j}\hat{\alpha}_{i,j}^{(\ell)} - \sum_{i,j} \frac{N_{i,j}^{(\ell)}}{N^{(\ell)}} = \sum_{(i,j)\in D_l}\left(\hat{\alpha}_{i,j}^{(\ell)} - \frac{N_{i,j}^{(\ell)}}{N^{(\ell)}}\right) + \sum_{(i,j) \in D_l^c}\left(\hat{\alpha}_{i,j}^{(\ell)} - \frac{N_{i,j}^{(\ell)}}{N^{(\ell)}}\right) \ge -3|D_l|\varepsilon + \sum_{(i,j) \in D_l^c}\left(\hat{\alpha}_{i,j}^{(\ell)} - \frac{N_{i,j}^{(\ell)}}{N^{(\ell)}}\right).
$
Hence, $ 0 \le \max_{i,j} \left \{\hat{\alpha}_{i,j}^{(\ell)} - \frac{N_{i,j}^{(\ell)}}{N^{(\ell)}} \right\}\le \sum_{(i,j) \in D_l^c}\left(\hat{\alpha}_{i,j}^{(\ell)} - \frac{N_{i,j}^{(\ell)}}{N^{(\ell)}}\right) \le3|D_l|\varepsilon \le 3KD\varepsilon $. By arbitrary $\varepsilon > 0$, we get $\lim_{\ell\rightarrow \infty} \max_{i,j} \left(\hat{\alpha}_{i,j}^{(\ell)} - \frac{N_{i,j}^{(\ell)}}{N^{(\ell)}} \right) = 0 $. Since for any $(i_0,j_0) $ we have 
$$\max_{i,j} \left\{\hat{\alpha}_{i,j}^{(\ell)} - \frac{N_{i,j}^{(\ell)}}{N^{(\ell)}} \right\} \ge \hat{\alpha}_{i_0,j_0}^{(\ell)} - \frac{N_{i_0,j_0}^{(\ell)}}{N^{(\ell)}} = -\sum_{i\neq i_0,j\neq j_0}\left(\hat{\alpha}_{i,j}^{(\ell)} - \frac{N_{i,j}^{(\ell)}}{N^{(\ell)}} \right) \ge - (KB-1)\max_{i,j} \left\{\hat{\alpha}_{i,j}^{(\ell)} - \frac{N_{i,j}^{(\ell)}}{N^{(\ell)}} \right \},  $$
we obtain $\lim_{\ell \rightarrow \infty} \hat{\alpha}_{i_0,j_0}^{(\ell)} - \frac{N_{i_0,j_0}^{(\ell)}}{N^{(\ell)}} = 0 = \alpha_{i_0,j_0}^* - \lim_{\ell \rightarrow \infty}\frac{N_{i_0,j_0}^{(\ell)}}{N^{(\ell)}} $ as desired. 
The proof is complete. \hfill $\blacksquare$

\subsection{Proof of Lemma \ref{lem: posterior convergence}}
\textbf{Proof.} Recall $M(t)$ is the total number of input data up to stage $t$. We have 
$$\pi_{j^c}^t = \frac{\pi_0(\theta_{j^c}) \prod_{\ell=1}^{M(t)} f_{\theta_{j^c}}(\xi_\ell)}{\sum_{j=1}^D \pi_0(\theta_{j}) \prod_{\ell=1}^{M(t)} f_{\theta_{j}}(\xi_\ell) }.$$
Take the inverse of both sides we obtain
$$ \frac{1}{\pi_{j^c}^t} = \frac{\sum_{j=1}^D \pi_0(\theta_{j}) \prod_{\ell=1}^{M(t)} f_{\theta_{j}}(\xi_\ell) }{\pi_0(\theta_{j^c}) \prod_{\ell=1}^{M(t)} f_{\theta_{j^c}}(\xi_\ell)}=1+ \sum_{j\neq \jc} \frac{\pi_0(\theta_j)}{\pi_0(\theta_\jc)} \prod_{\ell=1}^{M(t)} \frac{f_{\theta_j}(\xi_\ell)}{f_{\theta_\jc}(\xi_\ell)}.$$
We first prove $  \prod_{\ell=1}^{M(t)} \frac{f_{\theta_j}(\xi_\ell)}{f_{\theta_\jc}(\xi_\ell)} = \Opar{e^{-ct}}$ for some $c>0$. To see this, take log on both sides.
$$  \log \left(\prod_{\ell=1}^{M(t)} \frac{f_{\theta_j}(\xi_\ell)}{f_{\theta_\jc}(\xi_\ell)}\right) = \sum_{\ell=1}^{M(t)} \log \left(\frac{f_{\theta_j}(\xi_\ell)}{f_{\theta_\jc}(\xi_\ell)} \right)= M(t) \frac{1}{M(t)} \sum_{\ell=1}^{M(t)} \log \left(\frac{f_{\theta_j}(\xi_\ell)}{f_{\theta_\jc}(\xi_\ell)} \right).$$
By Assumption \ref{assump: optimality balance}.1, we obtain almost surely, $ \frac{1}{M(t)} \sum_{\ell=1}^{M(t)} \log \left(\frac{f_{\theta_j}(\xi_\ell)}{f_{\theta_\jc}(\xi_\ell)} \right) \rightarrow -\operatorname{KL}(\theta^c\|\theta_j) < 0$. Hence, there exists $0<c_1 <\operatorname{KL}(\theta^c\|\theta_j)$, such that for all sufficiently large $t$, 
$$  \log \left(\prod_{\ell=1}^{M(t)} \frac{f_{\theta_j}(\xi_\ell)}{f_{\theta_\jc}(\xi_\ell)}\right) \le -c_1 M(t).$$
Hence we obtain for all sufficiently large $t$
$$  \prod_{\ell=1}^{M(t)} \frac{f_{\theta_j}(\xi_\ell)}{f_{\theta_\jc}(\xi_\ell)} \le e^{-c_1 M(t)}.$$
Since $\frac{M(t)}{t} \rightarrow \bar{m}$ by Assumption \ref{assump: optimality balance}.2, $\frac{\pi_0(\theta_j)}{\pi_0(\theta_\jc)} >0$ is a constant and the fact that there are only $D_1$ different $j$s such that $j\neq j^c$. We can find $\kappa >0$, such that 
$$\sum_{j\neq \jc} \frac{\pi_0(\theta_j)}{\pi_0(\theta_\jc)} \prod_{\ell=1}^{M(t)} \frac{f_{\theta_j}(\xi_\ell)}{f_{\theta_\jc}(\xi_\ell)} = \Opar{e^{-\kappa t}}.$$ Hence
$$\pi_\jc^t = \frac{1}{1+ \Opar{e^{-\kappa t}}} = 1 -  \Opar{e^{-\kappa t}}.$$ The proof is complete.

\subsection{Proof of Theorem \ref{thm: convergence speed}.}
 Several lemmas are needed to complete the proof. First, 
Lemma \ref{lem: finite j} ensures that the design-input pairs under an incorrect input parameter will only be sampled finitely many times. \ref{lem:ratio2}-\ref{lem:ratio4} ensures that a positive ratio of the budget will be assigned to all design-input pairs under correct input parameter. In the proof of the following, a sample path of the simulation process is fixed. Furthermore, since the consistency hold, we use $b \rightarrow b^c$ almost surely hence, we drop the superscript $c$ in the following.

\begin{lemma} \label{lem: finite j}
    \item design-input pair $(i,j), \forall i\in\mathcal{I}, j \neq j^c$ will only be simulated finitely many times.
\end{lemma}
\textbf{Proof}
By input balance conditions, we have for any $j\neq j^c,i\in\mathcal{I}$ and $\ell$ sufficiently large, almost surely
$$\frac{\alpha_{i,j}^{(\ell)}}{\hat{\sigma}_{i,j}^{(\ell)}\pi_j^{t_\ell}} \ge \frac{n_0}{\hat{\sigma}_{i,j}^{(\ell)}} \frac{N^{(\ell)}}{\pi_j^{t_\ell}}\ge \frac{n_0 \Bar{n}}{2\hat{\sigma}^{(\ell)}_{i,j}} \times  O\left( \frac{1}{t_\ell} \exp\left( \kappa t_\ell \right)\right) \rightarrow \infty,$$ where the second inequality is because  $\frac{N^{(\ell)}}{t_\ell} \rightarrow \bar{n}$ by Assumption \ref{assump: optimality balance}.2, which implies for large $\ell$ (or $t_\ell$), $N^{(\ell)} \ge \frac{\Bar{n}}{2} t_\ell$. 
 At the same time, for the same design $i$ and sufficiently large $\ell$,
 $$ \frac{\alpha_{i,j^c}^{(\ell)}}{\hat{\sigma}^\pl_{i,j^c}\pi_{j^c}^{t_\ell}} \le \frac{4}{\sigma_{i,j^c}} <\infty,$$ 
 where the inequality holds since $\alpha^{(\ell)} \le 1, \pi^{t_\ell}_{j^c} = 1 -  O\left({e^{-\kappa t}}\right)  \ge \frac{1}{2}$ and $\hat{\sigma}_{i,j^c} \le 2 \sigma_{i,j^c}$ as $(i,j^c)$ will be simulated infinitely many times by Theorem \ref{thm:consistency}.2. This implies that the design-input pair $(i,j)$ with $j\neq j^c$ will only be simulated finitely many times.
\hfill $\blacksquare$
\begin{lemma} \label{lem: global O(1/l)}
(i) $(i,j^c), i\neq b^c$   is simulated at iteration $\ell$ implies 
$\left(\frac{\alpha_{b,j^c}^\pl}{\hat{\sigma}^\pl_{b,j^c}}\right)^2 - \sum_{i\neq b} \left(\frac{\alpha_{i,j^c}^\pl}{\hat{\sigma}_{i,j^c}^\pl} \right)^2 \ge -\Opar{\frac{1}{\ell^2}}$. (ii) Conversely, $(b,j^c)$ is simulated at iteration $\ell$ implies 
$\left(\frac{\alpha_{b,j^c}^\pl}{\hat{\sigma}^\pl_{b,j^c}}\right)^2 - \sum_{i\neq b} \left(\frac{\alpha_{i,j^c}^\pl}{\hat{\sigma}_{i,j^c}^\pl} \right)^2 \le \Opar{\frac{1}{\ell^2}}$.
\end{lemma}
\textbf{Proof} By Lemma \ref{lem: finite j}, we know $(i,j) i\in \mathcal{I},j\neq j^c$ will only be simulated finitely many times. Suppose after $\ell_0$, no $(i,j),j\neq j^c$ will be simulated. Then we have  
$$\left|\left(\frac{\alpha_{b,j}^\pl}{\hat{\sigma}^\pl_{b,j}}\right)^2 - \sum_{i\neq b} \left(\frac{\alpha_{i,j}^\pl}{\hat{\sigma}_{i,j}^\pl} \right)^2 \right|=\frac{1}{(N^\pl)^2} \left| \left(\frac{N_{b,j}^{(\ell_0)}}{\hat{\sigma}^{(\ell_0)}_{b,j}}\right)^2 - \sum_{i\neq b} \left(\frac{N_{i,j}^{(\ell_0)}}{\hat{\sigma}_{i,j}^{(\ell_0)}} \right)^2\right| = O(\frac{1}{\ell^2}).$$
Hence, if $(i,j^c)$ is simulated at $\ell$ and ${j}^* = \arg\max_{j}\left| \left(\frac{\alpha_{b,j}^\pl}{\hat{\sigma}^\pl_{b,j}}\right)^2 - \sum_{i\neq b} \left(\frac{\alpha_{i,j}^\pl}{\hat{\sigma}_{i,j}^\pl} \right)^2 \right| \neq j^c$, we know 
$$\left(\frac{\alpha_{b,j^c}^\pl}{\hat{\sigma}^\pl_{b,j^c}}\right)^2 - \sum_{i\neq b} \left(\frac{\alpha_{i,j^c}^\pl}{\hat{\sigma}_{i,j^c}^\pl} \right)^2 \ge \left(\frac{\alpha_{b,j^*}^\pl}{\hat{\sigma}^\pl_{b,j^*}}\right)^2 - \sum_{i\neq b} \left(\frac{\alpha_{i,j^*}^\pl}{\hat{\sigma}_{i,j^*}^\pl} \right)^2 \ge -\Opar{\frac{1}{\ell^2}}.$$
Otherwise if $j^* = j^c$, we know $ \left(\frac{\alpha_{b,j^c}^\pl}{\hat{\sigma}^\pl_{b,j^c}}\right)^2 - \sum_{i\neq b} \left(\frac{\alpha_{i,j^c}^\pl}{\hat{\sigma}_{i,j^c}^\pl} \right)^2 \ge 0$. Hence we prove (i). (ii) can be proved in a similar manner.



\begin{lemma} \label{lem:ratio2}
 \item $\lim\inf_{\ell\rightarrow \infty} \frac{\alpha_{i,j^c}^\pl}{\alpha_{i',j^c}^\pl} > 0$, $\forall i\neq i' \neq b,$ \text{ almost surely. }
\end{lemma}
\textbf{Proof}
Prove by contradiction. Suppose there exists $i,i'$ such that $\lim\inf_{\ell\rightarrow \infty} \frac{\alpha_{i,\jc}^\pl}{\alpha_{i',\jc}^\pl} =0 $.  For any positive constant $\varepsilon>0$, we can find a sufficiently large $\ell$ such that $(i',\jc)$ is sampled at $\ell$ and $\frac{\alpha_{i,\jc}^\pl}{\alpha_{i', \jc}^\pl} \le \varepsilon$. Since 
$\hat{\mu}_b^\pl,\hat{\mu}_i^\pl$ and $\hat{\mu}_{i'}^\pl$ all will converge to the true value $a.s.$ and $\hat{b}^\pl = b$ for $\ell$ sufficiently large. Then there exists constants $a,e$ and $U>L>0$, such that for $\ell$ sufficiently large, $0<\frac{1}{2}<\pi_\jc^\tl<1$, $0 < r \le \hat{\sigma}^\pl_{i,\jc} \le h$ and  $U>(\hat{\mu}_b^\pl-\hat{\mu}_i^\pl)^2$ and $L<(\hat{\mu}_b^\pl - \hat{\mu}_{i'}^\pl)^2$. 
Then
\begin{align}
    &\dfrac{(\hat{\mu}^\pl_{b} - \hat{\mu}_{i}^\pl)^2}{\sum\limits_{j=1}^D \frac{(\hat{\sigma}^\pl_{b,j})^2(\pi_j^\tl)^2}{\alpha_{b,j}^\pl} + \sum\limits_{j=1}^D\frac{(\hat{\sigma}^\pl_{i,j})^2(\pi_j^\tl)^2}{\alpha_{i,j}^\pl}} - \dfrac{(\hat{\mu}^\pl_{b} - \hat{\mu}_{i'}^\pl)^2}{\sum\limits_{j=1}^D \frac{(\hat{\sigma}^\pl_{b,j})^2(\pi_j^\tl)^2}{\alpha_{b,j}^\pl} + \sum\limits_{j=1}^D\frac{(\hat{\sigma}^\pl_{i',j})^2(\pi_j^\tl)^2}{\alpha_{i',j}^\pl}} \notag\\ 
    <&     \dfrac{U}{\sum\limits_{j=1}^D \frac{(\hat{\sigma}^\pl_{b,j})^2(\pi_j^\tl)^2}{\alpha_{b,j}^\pl} + \sum\limits_{j=1}^D\frac{(\hat{\sigma}^\pl_{i,j})^2(\pi_j^\tl)^2}{\alpha_{i,j}^\pl}} -    \dfrac{L}{\sum\limits_{j=1}^D \frac{(\hat{\sigma}^\pl_{b,j})^2(\pi_j^\tl)^2}{\alpha_{b,j}^\pl} + \sum\limits_{j=1}^D\frac{(\hat{\sigma}^\pl_{i',j})^2(\pi_j^\tl)^2}{\alpha_{i',j}^\pl}}\notag\\ 
    = & \frac{\sum\limits_{j=1}^D\frac{(\hat{\sigma}^\pl_{b,j})^2(\pi_j^\tl)^2}{\alpha_{b,j}^\pl}(U-L) + U \sum\limits_{j=1}^D\frac{(\hat{\sigma}^\pl_{i',j})^2(\pi_j^\tl)^2}{\alpha_{i',j}^\pl} -L\sum\limits_{j=1}^D \frac{(\hat{\sigma}^\pl_{i,j})^2(\pi_j^\tl)^2}{\alpha_{i,j}^\pl}}{(\sum\limits_{j=1}^D \frac{(\hat{\sigma}^\pl_{b,j})^2(\pi_j^\tl)^2}{\alpha_{b,j}^\pl} + \sum\limits_{j=1}^D\frac{(\hat{\sigma}^\pl_{i,j})^2(\pi_j^\tl)^2}{\alpha_{i,j}^\pl})(\sum\limits_{j=1}^D \frac{(\sigma_{b,j})^2(\pi_j^\tl)^2}{\alpha_{b,j}^\pl} + \sum\limits_{j=1}^D\frac{(\hat{\sigma}^\pl_{i',j})^2(\pi_j^\tl)^2}{\alpha_{i',j}^\pl})} \label{eq:ratio2.1}
\end{align}
By proof of Lemma \ref{lem: finite j}, we know $\frac{\pi_j^\tl}{\alpha_{i,j}^\pl} \rightarrow 0$ for $j\neq j^c$. Hence, for $\ell$ sufficiently large, we have $\frac{(\hat{\sigma}^\pl_{i,j})^2(\pi_j^\tl)^2}{\alpha_{i,j}^\pl} = \Opar{\frac{1}{\ell}}$ for $j\neq j^c$ by Lemma \ref{lem: posterior convergence}. 
 Then,
\begin{align}   &\sum\limits_{j=1}^D\frac{(\hat{\sigma}^\pl_{b,j})^2(\pi_j^\tl)^2}{\alpha_{b,j}^\pl}(U-L) + U \sum\limits_{j=1}^D\frac{(\hat{\sigma}^\pl_{i',j})^2(\pi_j^\tl)^2}{\alpha_{i',j}^\pl} -L\sum\limits_{j=1}^D \frac{(\hat{\sigma}^\pl_{i,j})^2(\pi_j^\tl)^2}{\alpha_{i,j}^\pl} \notag \\
     \le 
      &\frac{(\hat{\sigma}^\pl_{b,\jc})^2}{\alpha_{b,\jc}^\pl}(U-L) +\frac {U}{{\alpha_{i',\jc}^\pl}} {h^2} -\frac{L}{4\alpha_{i,\jc}^\pl} r^2 + \Opar{\frac{1}{\ell}} \notag \\
      \le &\frac{\hat{\sigma}^\pl_{b,\jc}\hat{\sigma}^\pl_{i',\jc}}{\alpha_{i',\jc}^\pl}(U-L) +\frac {U}{{\alpha_{i',\jc}^\pl}} {h^2} -\frac{L}{4\alpha_{i,\jc}^\pl} r^2 + \Opar{\frac{1}{\ell}} \label{eq:ratio2.2}\\
     \le & \frac{h^2}{\alpha_{i',\jc}^\pl}(U-L) +\frac {U}{{\alpha_{i',\jc}^\pl}} {h^2} -\frac{L}{4\varepsilon\alpha_{i',\jc}^\pl} r^2 + \Opar{\frac{1}{\ell}} \notag\\
     =  & \frac{1}{\alpha_{i',\jc}^\pl} \left\{ 
     {h^2}(U-L) + {U}h^2 -\frac{L}{4\varepsilon} r^2 + \Opar{\frac{1}{\ell}} 
     \right\}\label{eq:ratio2.3}
   \end{align}
\eqref{eq:ratio2.2} holds because: $i',j^c$ is sampled at iteration $\ell$, by Lemma \ref{lem: global O(1/l)} we have 
$\left(\frac{\alpha_{b,j}^\pl}{\hat{\sigma}_{b,j}}\right)^2 \ge \sum_i \left(\frac{\alpha_{i,\jc}^\pl}{\hat{\sigma}_{i,\jc}}\right)^2 - \Opar{\frac{1}{\ell^2}} \ge \left(\frac{\alpha_{i',\jc}^\pl}{\hat{\sigma}_{i',\jc}}\right)^2 -\Opar{\frac{1}{\ell^2}} $. Hence, $\frac{\alpha_{b,\jc}^\pl}{\hat{\sigma}_{b,\jc}} \ge \frac{\alpha_{i,\jc}^\pl}{\hat{\sigma}_{i,\jc}} - \Opar{\frac{1}{\ell}}, \forall i\neq b$. Furthermore, since we also have $\sum_{i\in\mathcal{I}} \frac{\alpha_{i,j^c}^\pl}{\hat{\sigma}_{i,j^c}^\pl} \ge \sum_{i\in\mathcal{I}} \frac{\alpha_{i,j^c}^\pl}{h} = 1/h -\Opar{\frac{1}{\ell}}$ since $(i,j),j\neq j^c$ will only be sampled finitely many times.  By pigeon hole principle we obtain $\frac{\alpha_{b,j^c}^\pl}{\hat{\sigma}_{b,j^c}\pl} \ge \frac{1}{Kh} - \Opar{\frac{1}{\ell}}$. Then $ \frac{\hat{\sigma}_{i',j^c}}{\alpha_{i',j^c}^\pl} \ge \frac{1}{\frac{\hat{\sigma}_{b,j^c}^\pl}{\alpha_{b,\jc}^\pl} + \Opar{\frac{1}{\ell}}} = \frac{\hat{\sigma}_{b,j^c}^\pl}{\alpha_{b,\jc}^\pl} - \Opar{\frac{1}{\ell}}$.

From \eqref{eq:ratio2.3}, if we choose $\varepsilon$ that makes \eqref{eq:ratio2.3}$< 0$. We obtain   $\dfrac{(\hat{\mu}^\pl_{b} - \hat{\mu}_{i}^\pl)^2}{\sum\limits_{j=1}^D \frac{(\hat{\sigma}^\pl_{b,j})^2(\pi_j^\tl)^2}{\alpha_{b,j}^\pl} + \sum\limits_{j=1}^D\frac{(\hat{\sigma}^\pl_{i,j})^2(\pi_j^\tl)^2}{\alpha_{i,j}^\pl}} -     \dfrac{(\hat{\mu}^\pl_{b} - \hat{\mu}_{i'}^\pl)^2}{\sum\limits_{j=1}^D \frac{(\hat{\sigma}^\pl_{b,j})^2(\pi_j^\tl)^2}{\alpha_{b,j}^\pl} + \sum\limits_{j=1}^D\frac{(\hat{\sigma}^\pl_{i',j})^2(\pi_j^\tl)^2}{\alpha_{i',j}^\pl}} <0$, which implies $(i',\jc)$ cannot be sampled at $\ell$, a contradiction.

\begin{lemma}\label{lem:ratio3}
(i) $\lim\inf_{\ell\rightarrow \infty} \frac{\alpha_{b,\jc}^\pl}{\alpha_{i,\jc}^\pl} > 0\  \forall i\neq b$ \text{ almost surely }; \ (ii) $\lim\inf_{\ell\rightarrow \infty} \frac{\alpha_{i,\jc}^\pl}{\alpha_{b,\jc}^\pl} > 0,\  \forall i\neq b $ \text{ almost surely. }
\end{lemma}

\textbf{proof of (i)} Prove by contradiction. By Lemma \ref{lem:ratio2}, there exists a positive constant $c>0$ such that $\frac{\alpha_{k,\jc}^\pl}{\alpha_{i,\jc}^\pl} \ge c $ for all $k\neq b$ and $\ell$ sufficiently large. Hence, if  $\lim\inf_{\ell\rightarrow \infty} \frac{\alpha_{b,\jc}^\pl}{\alpha_{i,\jc}^\pl} = 0$ for some $i$, then it holds for all $i$. For $ \varepsilon < \frac{\underline{\sigma}}{2\Bar{\sigma}K}$, we have there exists $\ell$ sufficiently large, $(i,j^c)$ is simulated at $\ell$ for some $i\neq b$ and $\alpha_{b,j^c}^\pl < \varepsilon$. Furthermore since $(i,j)$ can only be sampled finitely many times and $1 = \sum_{i,j} \alpha_{i,j}^\pl$. We can find $i_0\neq b$, $\alpha_{i_0,j^c}^\pl \ge \frac{1}{K}$. Furthermore, by consistency result, for $\ell$ sufficiently large we have 
$0<\underline{\sigma} \le \hat{\sigma}_{i,j^c}^\pl \le \Bar{\sigma} \forall i\in\mathcal{I}$.
Then we have 
$$  (\frac{\alpha_{b,\jc}^\pl}{\hat{\sigma}^\pl_{b,\jc}})^2 - {\sum\limits_{k\neq b} (\frac{\alpha_{k,\jc}^\pl}{\hat{\sigma}^\pl_{k,\jc}})^2} \le (\frac{\alpha_{b,\jc}^\pl}{\hat{\sigma}^\pl_{b,\jc}})^2 - (\frac{\alpha_{i_0,\jc}^\pl}{\hat{\sigma}^\pl_{i_0,\jc}})^2 \le \frac{\varepsilon^2}{\underline{\sigma}^2} - \frac{1}{K^2\Bar{\sigma}^2} < -\frac{1}{2K^2\Bar{\sigma}^2}.$$
This contradicts $(i,j^c)$ is simulated at $\ell$. Hence, we prove (i).

\textbf{Proof of (ii)} (ii) can be proved in a similar way. \hfill $\blacksquare$
\begin{lemma}\label{lem:ratio4}
  $ \lim\inf_{\ell\rightarrow \infty} \alpha_{i,j^c}^\pl > 0,\ \forall i,j$ \text{ almost surely. }
\end{lemma}
\textbf{ Proof } This is a direct result of Lemma \ref{lem: finite j}, \ref{lem:ratio2} and \ref{lem:ratio3}. By Lemma \ref{lem: finite j}, for $\ell$ sufficiently large, $\sum_{i\in\mathcal{I}} \alpha_{i,j^c}^\pl \ge 1/2$. By Lemma \ref{lem:ratio2} and \ref{lem:ratio3}, $\forall i_0 \in \mathcal{I}$, there exists $c>0$, such that for $\ell$ sufficiently large, $\alpha_{i_0,j^c} \ge c{\alpha_{i,j^c}^\pl}, i\neq i_0$. Then $2\alpha_{i_0,j^c}^\pl \ge \frac{\alpha_{i_0,j^c}^\pl}{\sum\limits_{i}\alpha_{i,j^c}^\pl}\ge\frac{c}{K} >0.$ The proof is complete. \hfill $\blacksquare$

Lemma \ref{lem:mu} and \ref{lem:variance} guarantee the convergence rate of the estimated expected performance and estimated variance, respectively. 

\begin{lemma} (Inter-changeable $O(\cdot)$ notation).  
    Suppose Assumptions in Theorem \ref{thm: convergence speed} hold true. Given a function $f(\ell)$, the following statements are equivalent almost surely:
    \begin{enumerate}
        \item $f(\ell) = O(\sqrt{\frac{\log \log\ell}{\ell}})$;
        \item $f(\ell) = O(\sqrt{\frac{\log \log t_\ell}{t_\ell}})$;
        \item $f(\ell) = O(\sqrt{\frac{\log \log N_{i,j}^\pl}{N_{i,j}^\pl}})$ for some $(i,j)$.
    \end{enumerate}
\end{lemma}
\textbf{Proof.}
1 and 2 are equivalent since by Assumption \ref{assump: optimality balance}, we know for $\ell$ large enough, $\frac{1}{2} \bar{n}t_\ell \le \ell \le 2\bar{n} t_\ell$. 1 implies 3 since $N_{i,j}^\pl \le N^\pl = \ell + KD n_0$ and for large $\ell$, $\sqrt{\frac{\log \log\ell}{\ell}}$ decreases in $\ell$. To show 3 implies 1, by Lemma \ref{lem:ratio4}, we know almost surely, there exists $C>0$, $\alpha_{i,j}^\pl \ge C$ for $\ell$ large, which implies $N_{i,j}^\pl \ge C N^\pl \ge C\ell$. This completes the proof. 
We will arbitrarily use the three $O(\cdot)$ notations in the following proofs. $\hfill\blacksquare$

\begin{lemma} \label{lem:mu}
$ |\hat{\mu}_i^\pl - {\mu}_i | = O(\sqrt{\frac{\log \log\ell}{\ell}})$ almost surely.
\end{lemma}
\textbf{Proof}
Notice that 
$\hat{\mu}_i^\pl - {\mu}_i =(\pi_\jc^\tl\hat{\mu}_{i,\jc}^\pl -  \mu^c\jc) + \sum_{j\neq j^c} \pi_j^\tl\hat{\mu}_{i,j}^\pl $. 
For the first term, 
By LIL, $|\hat{\mu}_{i,\jc}^\pl -\mu_{i,\jc}| = O(\sqrt{\frac{\log\log N_{i,j}^\pl}{N_{i,j}^\pl}}) = O(\sqrt{\frac{\log \log N_{i,j}^\pl}{N_{i,j}^\pl}}) = O(\sqrt{\frac{\log \log\ell}{\ell}}) $. Furthermore, by Lemma \ref{lem: posterior convergence}, $|\pi_j^\tl-1| = \Opar{\exp (-\kappa \tl)} = O(\sqrt{\frac{\log \log t_\ell}{t_\ell}}) = O(\sqrt{\frac{\log \log\ell}{\ell}})$. Hence,
$$|\pi_\jc^\tl\hat{\mu}_{i,\jc}^\pl -\mu_{i,\jc}| = |\pi_\jc^\tl(\hat{\mu}_{i,\jc}^\pl -\mu_{i,\jc}) +(\pi_\jc^\tl- 1) \mu_{i,\jc}| \le  \pi_\jc^\tl|\hat{\mu}_{i,\jc}^\pl -\mu_{i,\jc}| + |\pi_\jc^\tl- 1| |\mu_{i,\jc}| =  O(\sqrt{\frac{\log \log\ell}{\ell}}).$$ 
For the second term, $ \pi_j^\tl \hat{\mu}_{i,j}^\pl =  \hat{\mu}_{i,j}^\pl \cdot \Opar{\exp(-\kappa\tl)} = \Opar{\logll}$.  \hfill $\blacksquare$
\begin{lemma} \label{lem:variance}
$ |(\hat{\sigma}_{i,\jc}^\pl)^2 - (\sigma_{i,\jc})^2 | = O(\sqrt{\frac{\log \log\ell}{\ell}})$ almost surely. As a result, $ |\hat{\sigma}_{i,\jc}^\pl - \sigma_{i,\jc} | = O(\sqrt{\frac{\log \log\ell}{\ell}})$ almost surely.
\end{lemma}

\textbf{Proof}
Since 
\begin{align*}
    \sum_{s=1}^{N_{i,\jc}^\pl}(X_{i,\jc}^{(s)} - \hat{\mu}_{i,\jc}^\pl)^2 =  \sum_{s=1}^{N_{i,\jc}^\pl}(X_{i,\jc}^{(s)} - {\mu}_{i,\jc})^2 -N_{i,\jc}^\pl (\hat{\mu}_{i,\jc}^\pl-\mu_{i,\jc})^2, \\
\end{align*}
We have 
$$ (\hat{\sigma}_{i,\jc}^\pl)^2 - (\sigma_{i,\jc})^2 = \frac{1}{N_{i,\jc}^\pl-1}\sum_{s=1}^{N_{i,\jc}^\pl}[(X_{i,\jc}^{(s)} - {\mu}_{i,\jc})^2-\sigma_{i,\jc}^2] - \frac{N_{i,\jc}^\pl}{N_{i,\jc}^\pl-1}(\hat{\mu}_{i,\jc}^\pl-\mu_{i,\jc})^2 +\frac{\sigma_{i,\jc}^2}{N_{i,\jc}^\pl-1}.$$
Since $(X_{i,\jc}^{(s)} - {\mu}_{i,\jc})^2-\sigma_{i,\jc}^2$ are i.i.d. with mean $0$, by LIL, we have with probability $1$,
$$\left|\frac{1}{N_{i,\jc}^\pl-1}\sum_{s=1}^{N_{i,\jc}^\pl}[(X_{i,\jc}^{(s)} - {\mu}_{i,\jc})^2-\sigma_{i,\jc}^2\right|= O\left(\sqrt{\frac{\log\log N_{i,\jc}^\pl}{N_{i,\jc}^\pl}}\right)= O(\sqrt{\frac{\log \log N_{i,\jc}^\pl}{N_{i,\jc}^\pl}}) = O\left(\sqrt{\frac{\log \log\ell}{\ell}}\right).$$
Further since $ \left|\hat{\mu}_{i,\jc}^\pl-\mu_{i,\jc}\right| = O\left(\sqrt{\frac{\log \log\ell}{\ell}}\right)$ and $\left|\frac{\sigma_{i,\jc}^2}{N_{i,\jc}^\pl-1}\right| = O\left(\sqrt{\frac{\log \log\ell}{\ell}}\right) $, we get the desired result. \hfill $\blacksquare$

Lemma \ref{lem:length} is a simple but useful result that we will use frequently in the following proof.
\begin{lemma} \label{lem:length}
Let $(i,\jc)$ be a fixed design-input pair. Suppose $(i,\jc)$ is sampled at iteration $r$. Let $\ell_r = \inf\{\ell>0: \mathbf{1}_{(i,\jc)}^{(r+\ell)} = 1\}$. Hence $r+\ell_r$ is the next iteration $(i,\jc)$ will be sampled after $r$. Then we have $r<r+\ell_r = O(r)$ \text{ almost surely }.
\end{lemma}
\textbf{Proof}
Prove by contradiction. Suppose $\forall C_0 >0$, there exists an iteration $r$ such that $t_r > C_0 r$. We have 
$$\alpha_{i,\jc}^{(r+\ell_r)} = \frac{N_{i,\jc}^{(r+\ell_r)}}{N^{(r+\ell_r)}} = \frac{N_{i,\jc}^{(r)}+1}{N^{(r)}+\ell_r} <\frac{2(n_0+r)}{KBn_0+(C_0+1) r} < \frac{3}{C_0}$$ for large $r$.  The first inequality holds since $N(r) = KBn_0 + r$ and $N_{i,\jc}^{(0)} = n_0$. By the arbitrariness of $C_0$ and the fact that if $C_0 \rightarrow \infty$, the iteration $r$ that satisfy $\ell_r > C_0 r$ must also go to $\infty$. We have $\lim\inf_{\ell\rightarrow \infty} \alpha_{i,\jc}^\pl = 0$, contradicting Lemma \ref{lem:ratio4}.(ii). \hfill $\blacksquare$

\subsection{Proof of Theorem \ref{thm: convergence speed}.1.}
\begin{lemma} \label{lem: approx optimal 1}
Suppose the assumptions in Theorem \ref{thm: convergence speed} hold true. Then we have $$|\hat{\alpha}^\pl_{i,j} - \alpha^*_{i,j}| = O\left(\sqrt{\frac{\log \log\ell}{\ell}}\right) \text{ almost surely }.$$  Here $\alpha^*_{i,j} = 0$ for $j\neq \jc$.  
\end{lemma}
\textbf{Proof.}
Notice both  $\alpha^*_{i,j}$ and $\hat{\alpha}^\pl_{i,j}$, $i\neq b$ can be computed explicitly as 
$$\alpha_{i,j}^* = \frac{w^*_{i,j}}{\sum_{i',j'} w^*_{i',j'} },
$$
where $$
    w^*_{i,j} = {\pi_j^c} \sigma_{i,j}\sum_{k=1}^D\sigma_{i,k}{\pi_k^c}/({\mu}_b - {\mu}_i)^2
$$ if $i\neq b$ and $$w^*_{b,j} = \sqrt{ \sum_{i'\neq b} \frac{\sigma_{b,j'}^2 w_{i',j}^2}{\sigma_{i',j}^2}},$$
 $\pi_j^c = 0$ for $j\neq j^c$.  
 $$\hat{\alpha}^\pl_{i,j} = \frac{\hat{w}_{i,j}^\pl}{\sum_{i',j'} w_{i',j'}^\pl}, $$
 where 
$$ \hat{w}^\pl_{i,j} = {\pi_j^\tl}{\hat{\sigma}^\pl_{i,j}\sum_{k=1}^D\hat{\sigma}^\pl_{i,k}{\pi_k^\tl}/(\hat{\mu}_b^\pl - \hat{\mu}^\pl_i)^2},$$
if $i\neq b$ and 
$$ \hat{w}^\pl_{b,j} = \sqrt{ \sum_{i'\neq b} \frac{(\hat{\sigma}^\pl_{b,j'})^2 (\hat{w}^\pl_{i',j})^2}{(\hat{\sigma}^\pl_{i',j})^2}}.$$
 Since $\pi_j^\tl = \Opar{\logll}$, $\hat{w}_{i,j}^\pl = w^*_{i,j} + \Opar{\logll}$.
Hence, we also have  $\hat{\alpha}_{i,j}^\pl = \alpha_{i,j}^* + \Opar{\logll}$.
\hfill $\blacksquare$
 

\textbf{Proof of Theorem \ref{thm: convergence speed}.1.}
We first prove
$$ \max_{i,j} \{\alpha^\pl_{i,j} - \hat{\alpha}_{i,j}^\pl\} = O\left( \sqrt{\frac{\log \log\ell}{\ell}}\right).$$ 

Let $(i_\ell,j_\ell) = \arg \max_{i,j} \{\alpha^\pl_{i,j} - \hat{\alpha}_{i,j}^\pl\}$. Then since $\sum_{i,j} \alpha^\pl = \sum_{i,j}\hat{\alpha}^\pl = 1$, we know  $\alpha^\pl_{i_\ell,j_\ell} - \hat{\alpha}_{i_\ell,j_\ell}^\pl \ge 0$ and $(i_\ell,j_\ell)$ cannot be simulated at $\ell$.  Then, let $r$ be the last time $(i_\ell,j_\ell)$ is simulated, where we must have 
$\alpha^{(r)}_{i_\ell,j_\ell} - \hat{\alpha}_{i_\ell,j_\ell}^{(r)} \le 0$. By Lemma \ref{lem: approx optimal 1}, we have there exists $C_1>0$, such that for both $r$ and $\ell$ sufficiently large, $|\hat{\alpha}^{(r)}_{i_\ell,j_\ell} - \alpha^*_{i_\ell,j_\ell}| = O\left(\sqrt{\frac{\log\log r}{r}}\right)$ and $|\hat{\alpha}^\pl_{i_\ell,j_\ell} - \alpha^*_{i_\ell,j_\ell}| = O\left( \sqrt{\frac{\log \log\ell}{\ell}}\right)$. 
Furthermore by Lemma \ref{lem:length}, we know $\ell =O(r)$ and hence $\sqrt{\frac{\log \log\ell}{\ell}} = O\left(\sqrt{\frac{\log\log r}{r}}\right) $. This implies $|\alpha^{(r)}_{i_\ell,j_\ell} - \alpha^*_{i_\ell,j_\ell}| = O\left( \sqrt{\frac{\log \log\ell}{\ell}}\right)$ and $|\alpha^{(r)}_{i_\ell,j_\ell} -\alpha^{(\ell)}_{i_\ell,j_\ell} | = O\left( \sqrt{\frac{\log \log\ell}{\ell}}\right)$. Hence, there exists $C_2 > 0$,

$$ 0\le \alpha^\pl_{i_\ell,j_\ell} - \hat{\alpha}^\pl_{i_\ell,j_\ell} \le  \alpha^{(r)}_{i_\ell,j_\ell} + \frac{1}{N^\pl}- \hat{\alpha}^{(r)}_{i_\ell,j_\ell} + C_1 \logll \le \frac{1}{N^\pl}+  C_1 \logll = O\left(\sqrt{\frac{\log \log\ell}{\ell}}\right),$$
where $\alpha^\pl_{i_\ell,j_\ell}- \alpha^{(r)}_{i_\ell,j_\ell} = \frac{N_{i_\ell,j_\ell}^{(r)}+1}{N^\pl} - \frac{N_{i_\ell,j_\ell}^{(r)}}{N^{(r)}}\le \frac{N_{i_\ell,j_\ell}^{(r)}+1}{N^\pl} - \frac{N_{i_\ell,j_\ell}^{(r)}}{N^\pl} = \frac{1}{N^\pl}$.

Second, we prove $\forall (i,j)$, $|\alpha_{i,j}^\pl - \hat{\alpha}_{i,j}^\pl| = O\left(\sqrt{\frac{\log \log\ell}{\ell}}\right)$. Let $I_1 = \{(i,j), \alpha_{i,j}^\pl - \hat{\alpha}_{i,j}^\pl \ge 0\}$ and $I_2 = \{(i,j), \alpha_{i,j}^\pl - \hat{\alpha}_{i,j}^\pl \le 0\}$. We have if $(i,j) \in I_1$, $ 0 \le \alpha_{i,j}^\pl - \hat{\alpha}_{i,j}^\pl \le \alpha_{i_\ell,j_\ell}^\pl - \hat{\alpha}_{i_\ell,j_\ell}^\pl  = O\left(\sqrt{\frac{\log \log\ell}{\ell}}\right)$. Else if $(i,j)\in I_2$, 
$0\ge \alpha_{i,j}^\pl - \hat{\alpha}_{i,j}^\pl \ge \sum_{(i',j')\in I_2} \alpha_{i',j'}^\pl - \hat{\alpha}_{i',j'}^\pl = - \sum_{(i',j')\in I_1} ( \alpha_{i',j'}^\pl - \hat{\alpha}_{i',j'}^\pl) \ge -|I_2|(\alpha^\pl_{i_\ell,j_\ell} - \hat{\alpha}^\pl_{i_\ell,j_\ell}) \ge -KD(\alpha^\pl_{i_\ell,j_\ell} - \hat{\alpha}^\pl_{i_\ell,j_\ell}) = O\left(\sqrt{\frac{\log \log\ell}{\ell}}\right).  $

Third, noticing $|\hat{\alpha}_{i,j}^\pl - \alpha^*_{i,j}| = O\left(\sqrt{\frac{\log \log\ell}{\ell}}\right) $, we obtain 
$$ |{\alpha}_{i,j}^\pl - \alpha^*_{i,j}| \le |\alpha_{i,j}^\pl - \hat{\alpha}_{i,j}^\pl |+  |\hat{\alpha}_{i,j}^\pl - \alpha^*_{i,j}| = O\left(\sqrt{\frac{\log \log\ell}{\ell}}\right) + O\left(\sqrt{\frac{\log \log\ell}{\ell}}\right) = O\left(\sqrt{\frac{\log \log\ell}{\ell}}\right).$$
This completes the proof. \hfill $\blacksquare$
\subsection{Proof of Theorem \ref{thm: convergence speed}.2}
We need some more lemmas to prove the result.

The following Lemma \ref{lem:NofS_bi} bounds the amount of budget allocated to a non-optimal design-input pair $(i,\jc)$ between two successive samples of the best design-input pair $(b,\jc)$ under the same input realization.

\begin{lemma} \label{lem:NofS_bi}
Suppose $(b,\jc)$ is sampled at iteration $r$. Let $\ell_r=\inf\{\ell>0:\mathbf{1}_{(b,\jc)}^{(r+\ell)}=1\}$. $r+\ell_r$ is the next iteration at which $(b,\jc)$ is sampled. Then between the two samples of $(b,\jc)$, the number of samples that can be allocated to $(i,\jc), i\neq b$ is at most $O(\sqrt{r\log \log r})$ \text{ almost surely }. 
\end{lemma}
\textbf{Proof}
Fix a non-optimal design $i$. Let $s_r = \sup\{\ell<\ell_r:\mathbf{1}_{(i,\jc)}^{(r+\ell)}=1\}$. $r+s_r$ is the last time before $r+\ell_r$ at which $(i,\jc)$ is sampled. If $s_r<0$, then the lemma holds true, otherwise, assume $s_r > 0$. Since $(b,\jc)$ is sampled at $r$, then $$ \left( \frac{N_{b,\jc}^{(r)}}{\hat{\sigma}^{(r)}_{b,\jc }   }\right)^2 \le \sum\limits_{k\neq b} \left( \frac{N_{k,\jc}^{(r)}}{\hat{\sigma}^{(r)}_{k,\jc}}\right)^2 + \Opar{1} \le \sum\limits_{k\neq b} \left( \frac{N_{k,\jc}^{(r)}}{\hat{\sigma}^{(r)}_{k,\jc}}\right)^2 + \Opar{\sqrt{r{\log\log r}}} $$ by Lemma \ref{lem: global O(1/l)}.  
Similarly, since $(i,\jc)$ is sampled at $r+s_r$,
$$ \left( \frac{N_{b,\jc}^{(r+s_r)}}{\hat{\sigma}^{(r+s_r)}_{b,\jc }}\right)^2 \ge \sum\limits_{k\neq b} \left( \frac{N_{k,\jc}^{(r+s_r)}}{\hat{\sigma}^{(r+s_r)}_{k,\jc}}\right)^2 - \Opar{\sqrt{r{\log\log r}}}, $$
where $r+s_r = O(r)$ by Lemma \ref{lem:length}. Furthermore notice $\hat{\sigma}_{i,j^c}^{(r)} = \hat{\sigma}_{i,j^c}^{(r+s_r)} + \Opar{\sqrt{\frac{\log\log r}{r}}}$. Then
\begin{align}
    0 &\le \left( \frac{N_{b, \jc}^{(r+s_r)}}{\hat{\sigma}^{(r+s_r)}_{b, \jc} }\right)^2 - \sum\limits_{k\neq b} \left( \frac{N_{k, \jc}^{(r+s_r)}}{\hat{\sigma}^{(r+s_r)}_{k, \jc} }\right)^2 + \Opar{\sqrt{r{\log\log r}}} \notag\\
    &=\left( \frac{N_{b, \jc}^{(r)}+1}{\hat{\sigma}^{(r+s_r)}_{b, \jc} }\right)^2 - \sum\limits_{k\neq b} \left( \frac{N_{k, \jc}^{(r+s_r)}}{\hat{\sigma}^{(r+s_r)}_{k, \jc} }\right)^2 + \Opar{\sqrt{r{\log\log r}}} \notag\\
    &\le \left( \frac{N_{b, \jc}^{(r)}}{\hat{\sigma}^{(r+s_r)}_{b, \jc} }\right)^2 - \sum\limits_{k\neq b} \left( \frac{N_{k, \jc}^{(r+s_r)}}{\hat{\sigma}^{(r+s_r)}_{k, \jc} }\right)^2+ \Opar{r} \notag\\
    &\le \sum\limits_{k\neq b} \left( \frac{N_{k,   \jc}^{(r)}}{\hat{\sigma}^{(r)}_{k,   \jc}}\right)^2 -\sum\limits_{k\neq b} \left( \frac{N_{k, \jc}^{(r+s_r)}}{\hat{\sigma}^{(r+s_r)}_{k, \jc} }\right)^2+\Opar{r\sqrt{r{\log\log r}}} \notag\\
    &\le \sum\limits_{k\neq b} \left( \frac{N_{k, \jc}^{(r)}}{\hat{\sigma}^{(r)}_{k, \jc}}\right)^2 -\sum\limits_{k\neq b} \left( \frac{N_{k, \jc}^{(r+s_r)}}{\hat{\sigma}^{(r)}_{k, \jc}}\right)^2+\Opar{r\sqrt{r{\log\log r}}}\label{eq:NofS_bi.5}
\end{align}
 Then, since for each $k\neq b$,  $\left( \frac{N_{k,\jc}^{(r+s_r)}}{\hat{\sigma}^{(r)}_{k,\jc}}\right)^2 - \left( \frac{N_{k,\jc}^{(r)}}{\hat{\sigma}^{(r)}_{k,\jc}}\right)^2 \ge 0 $. We obtain
\begin{equation*}
    \left( \frac{N_{i,\jc}^{(r+s_r)}}{\hat{\sigma}^{(r)}_{i,\jc}}\right)^2 - \left( \frac{N_{i,\jc}^{(r)}}{\hat{\sigma}^{(r)}_{i,\jc}}\right)^2 \le \Opar{r\sqrt{r\log \log r}}.
\end{equation*}
 Since
\begin{equation} \label{eq:NofS_bi.6}
    \Opar{r\sqrt{r\log \log r}} \ge (N_{i,\jc}^{(r+s_r)})^2 - (N_{i,\jc}^{(r)})^2 = (N_{i,\jc}^{(r+s_r)} -N_{i,\jc}^{(r)})(N_{i,\jc}^{(r+s_r)} + N_{i,\jc}^{(r)})\ge 2N_{i,\jc}^{(r)}(N_{i,\jc}^{(r+s_r)} -N_{i,\jc}^{(r)})
\end{equation}
By Lemma \ref{lem:ratio4}.(ii), we have $\frac{1}{N_{i,j^c}^{(r)} } \Opar{r\sqrt{r\log \log r}} =  \Opar{\sqrt{r\log \log r}}$. This implies  $N_{i,\jc}^{(r+\ell_r)} -N_{i,\jc}^{(r)} = N_{i,\jc}^{(r+s_r)} -N_{i,\jc}^{(r)} = O(\sqrt{r\log \log r})$. \hfill $\blacksquare$

Conversely, the following Lemma \ref{lem:NofS_banyi} bounds the amount of budget allocated to the optimal design-input pair $(b,\jc)$ between two successive samples of any two non-optimal design-input pair  under the same input realization.  

\begin{lemma}\label{lem:NofS_banyi}
 Suppose at iteration $r$ a non-optimal design $i_1$ is sampled. Let $\ell_r = \inf\{\ell>0:\exists i \neq b, \mathbf{1}_{(i,\jc)}^{(r+\ell)} =1\}$. $r+\ell_r$ is the next iteration at which a non-optimal design is sampled. Then between iteration $r$ and $r+\ell_r$, the number of samples that can be allocated to $(b,\jc)$ is $O(\sqrt{r\log \log r})$ \text{ almost surely }.
\end{lemma}
\textbf{Proof} The proof is similar to Lemma \ref{lem:NofS_bi}.
Define $s_r = \sup\{ \ell<\ell_r: \mathbf{1}_{(b,\jc)}=1\}$. $r+s_r$ is the last time before $r+\ell_r$ the optimal design is sampled. If $s_r < 0$, then the lemma holds. Otherwise assume $s_r>0$. Since $(b,\jc)$ is sampled at $r+s_r$,$ \left( \frac{N_{b,\jc}^{(r+s_r)}}{\hat{\sigma}^{(r+s_r)}_{b,\jc }}\right)^2 \le \sum\limits_{k\neq b} \left( \frac{N_{k,\jc}^{(r+s_r)}}{\hat{\sigma}^{(r+s_r)}_{k,\jc}}\right)^2 + \Opar{\sqrt{r{\log\log r}}}$. Since $(i_1,j^c)$ is sampled at $r$, we  have 
$ \left( \frac{N_{b,\jc}^{(r)}}{\hat{\sigma}^{(r)}_{b,\jc }   }\right)^2 \ge \sum\limits_{k\neq b} \left( \frac{N_{k,\jc}^{(r)}}{\hat{\sigma}^{(r)}_{k,\jc}}\right)^2 - \Opar{\sqrt{r{\log\log r}}}$.
Further by Lemma \ref{lem:length}, we have $s_r < t_r = O(r)$. 
We have
\begin{align}
    0 &\ge \left( \frac{N_{b, \jc}^{(r+s_r)}}{\hat{\sigma}^{(r+s_r)}_{b, \jc} }\right)^2 -\sum\limits_{k\neq b} \left( \frac{N_{k, \jc}^{(r+s_r)}}{\hat{\sigma}^{(r+s_r)}_{k, \jc} }\right)^2 -\Opar{\rlogr}\notag\\
    &=\left( \frac{N_{b,\jc}^{(r+s_r)}}{\hat{\sigma}^{(r+s_r)}_{b,\jc} }\right)^2 -\sum\limits_{k\neq i_1 \neq b} \left( \frac{N_{k,\jc}^{(r)}}{\hat{\sigma}^{(r+s_r)}_{k,\jc} }\right)^2 - \left( \frac{N_{i_1,\jc}^{(r)}+1}{\hat{\sigma}^{(r+s_r)}_{i_1,\jc} }\right)^2 -\Opar{\rlogr} \notag\\
    &\ge \left( \frac{N_{b,\jc}^{(r+s_r)}}{\hat{\sigma}^{(r+s_r)}_{b,\jc} }\right)^2 -\sum\limits_{k\neq i_1 \neq b} \left( \frac{N_{k,\jc}^{(r)}}{\hat{\sigma}^{(r+s_r)}_{k,\jc} }\right)^2 - \left( \frac{N_{i_1,\jc}^{(r)}}{\hat{\sigma}^{(r+s_r)}_{i_1,\jc} }\right)^2 - \Opar{r} \notag\\
    & \ge \left( \frac{N_{b,\jc}^{(r+s_r)}}{\hat{\sigma}^{(r)}_{b,\jc} }\right)^2 -\sum\limits_{k \neq b} \left( \frac{N_{k,\jc}^{(r)}}{\hat{\sigma}^{(r)}_{k,\jc} }\right)^2 -\Opar{r\rlogr}  \notag\\
    &\ge \left( \frac{N_{b,\jc}^{(r+s_r)}}{\hat{\sigma}^{(r)}_{b,\jc} }\right)^2 -\left( \frac{N_{b,\jc}^{(r)}}{\hat{\sigma}^{(r)}_{b,\jc} }\right)^2 - \Opar{r\rlogr} \notag
\end{align}
With similar reason as in proof of Lemma \ref{lem:NofS_bi}, this implies
$$ N_{b,\jc}^{(r+s_r)}-N_{b,\jc}^{(r)} \le \Opar{\rlogr}$$ The proof is complete. \hfill $\blacksquare$

\begin{lemma} \label{lem:total_balance}
$\left|\left(\frac{\alpha_{b,\jc}^\pl}{\hat{\sigma}^\pl_{b,\jc}}\right)^2 - \sum\limits_{i\neq b} \left(\frac{\alpha_{i,\jc}^\pl}{\hat{\sigma}^\pl_{i,\jc}}\right)^2\right| = \Opar{\logll}  \quad a.s.$. 
\end{lemma}
\textbf{Proof}
First notice by Lemma \ref{lem: finite j}, when $\ell$ is large enough, only design-input pair under the true parameter will be simulated.  Let $ |\Delta^\pl| :=\left|\left(\frac{\alpha_{b,\jc}^\pl}{\hat{\sigma}^\pl_{b,\jc}}\right)^2 - \sum\limits_{i\neq b} \left(\frac{\alpha_{i,\jc}^\pl}{\hat{\sigma}^\pl_{i,\jc}}\right)^2\right|$. Fix an iteration.

\textbf{Case 1.} $\Delta\pl \le 0$, let $\ell_s = \sup\{r < \ell: \exists i\neq b,\  \mathbf{1}_{(i,\jc)}^{(r)} =1\}$. Then $\ell_s$ is the last time before $\ell$ that a non-optimal design, denote by $i_0$, is sampled under $\theta_\jc$.  We have 
\begin{align}
    0 &\ge \left(\frac{N_{b,\jc}^{(\ell)}}{\hat{\sigma}^{(\ell)}_{b,\jc}}\right)^2 - \sum\limits_{i\neq b} \left(\frac{N_{i,\jc}^{(\ell)}}{\hat{\sigma}^{(\ell)}_{i,\jc}}\right)^2 \notag\\
    & \ge \left(\frac{N_{b,\jc}^{(\ell_{s})}}{\hat{\sigma}^{(\ell)}_{b,\jc}}\right)^2 - \sum\limits_{i\neq i_0\neq b} \left(\frac{N_{i,\jc}^{(\ell_{s})}}{\hat{\sigma}^{(\ell)}_{i,\jc}}\right)^2 - \left(\frac{N_{i_0,\jc}^{(\ell_{s})}+1}{\hat{\sigma}^{(\ell)}_{i,\jc}}\right)^2 \notag\\
    &\ge \left(\frac{N_{b,\jc}^{(\ell_{s})}}{\hat{\sigma}^{(\ell)}_{b,\jc}}\right)^2 - \sum\limits_{i\neq i_0} \left(\frac{N_{i,\jc}^{(\ell_{s})}}{\hat{\sigma}^{(\ell)}_{i,\jc}}\right)^2 -  \Opar{\ell} \notag
\end{align}
Divide both sides by $(N^{(\ell_s)})^2$ and notice $\ell = O(\ell_s)$, we have

\begin{align}
    0 &> \left(\frac{\alpha_{b,\jc}^{(\ell_{s})}}{\hat{\sigma}^{(\ell)}_{b,\jc}}\right)^2 - \sum\limits_{i\neq i_0} \left(\frac{\alpha_{i,\jc}^{(\ell_{s})}}{\hat{\sigma}^{(\ell)}_{i,\jc}}\right)^2 - \Opar{\frac{1}{\ell}} \notag\\
    &\ge \left(\frac{\alpha_{b,\jc}^{(\ell_{s})}}{\hat{\sigma}^{(\ell_{s})}_{b,\jc}}\right)^2 - \sum\limits_{i\neq i_0} \left(\frac{\alpha_{i,\jc}^{(\ell_{s})}}{\hat{\sigma}^{(\ell_{s})}_{i,\jc}}\right)^2 - \Opar{\logll} \notag.
\end{align} 
Which completes the proof for case $1$.

\textbf{Case 2.} $\Delta\pl \ge 0$. Let  $\ell_{s} = \sup\{r < \ell :\mathbf{1}_{(b,\jc)}^{(r)} =1\}$. $\ell_{s}$ is the last time before $\ell$ at which  $(b,\jc)$ is sampled. By Lemma \ref{lem:length}, $\ell = O(\ell_s)$ since $(b,\jc)$ is not sampled between $\ell_{s}$ and $\ell$; by Lemma \ref{lem:NofS_bi}, $N_{i,\jc}^\pl - N_{i,\jc}^{(\ell_s)} =\Opar{\llogl}, \ \forall i\neq b$. Then,
\begin{align}
    0 &\le \left(\frac{N_{b,\jc}^{(\ell)}}{\hat{\sigma}^{(\ell)}_{b,\jc}}\right)^2 - \sum\limits_{i\neq i_0} \left(\frac{N_{i,\jc}^{(\ell)}}{\hat{\sigma}^{(\ell)}_{i,\jc}}\right)^2  \notag\\
    & \le \left(\frac{N_{b,\jc}^{(\ell_{s})}+1}{\hat{\sigma}^{(\ell)}_{b,\jc}}\right)^2 - \sum\limits_{i\neq i_0} \left(\frac{N_{i,\jc}^{(\ell_{s})}}{\hat{\sigma}^{(\ell)}_{i,\jc}}\right)^2 +  \Opar{\ell \llogl} \notag\\
    &\le  \left(\frac{N_{b,\jc}^{(\ell_{s})}}{\hat{\sigma}^{(\ell_{s})}_{b,\jc}}\right)^2 - \sum\limits_{i\neq i_0} \left(\frac{N_{i,\jc}^{(\ell_{s})}}{\hat{\sigma}^{(\ell_{s})}_{i,\jc}}\right)^2 +  \Opar{\ell \llogl}  \notag \\
    &\le \Opar{\ell \llogl}  \notag
\end{align}
 Divide both sides by $(N^\pl)^2$, we obtain $ 0 \le \Delta_{\jc}\pl \le \Opar{\logll}$. The proof is complete. \hfill $\blacksquare$

\textbf{Proof of \eqref{thmeq:error total balance} in Theorem \ref{thm: convergence speed}.2.}

\textbf{Proof.}
By Lemma \ref{lem:variance} and Lemma \ref{lem:total_balance}, 

\begin{align*}
  \left|\left(\frac{\alpha_{b,\jc}^{(r)}}{\sigma_{b,\jc}}\right)^2 - \sum\limits_{i\neq b} \left(\frac{\alpha_{i,\jc}^{(r)}}{\sigma_i(\theta^c)}\right)^2\right| = & \left|\left(\frac{\alpha_{b,\jc}^{(r)}}{\hat{\sigma}^{(r)}_{b,\jc}}\right)^2 - \sum\limits_{i\neq b} \left(\frac{\alpha_{i,\jc}^{(r)}}{\hat{\sigma}^{(r)}_{i,\jc}}\right)^2\right|+O\left(\sqrt{\frac{\log\log r}{r}}\right) \\
  =&O\left(\sqrt{\frac{\log\log r}{r}}\right) 
\end{align*} \hfill $\blacksquare$

The next Lemma \ref{lem:intermediate_lemma} is a little technical, which is used to bound the amount of budget allocated to a non-optimal design-input pair $(b,\jc)$ between two successive samples of a non-optimal design-input pair $(i,\jc)$, as shown in Lemma \ref{lem:NofS_ib}.  
\begin{lemma}\label{lem:intermediate_lemma}
Under $\jc$, suppose a non-optimal design $(k,\jc)$ is sampled at iteration $r$. Define 
$$
\left\{\begin{aligned}
&\ell_r := \inf_\ell\{\ell>0: \mathbf{1}_{(k,\jc)}^{(r+\ell)} =1\}\\
&s_r':=\sup_\ell \{  \ell<\ell_r: \mathbf{1}_{(b,\jc)}^{(r+\ell)} = 1\}\\
&s_r :=\sup_\ell \{\ell<s_r':\mathbf{1}_{(i,\jc)}^\pl =1 \text{ for some } i\neq b\}\\
&d_{i,\jc}^{(r,q)} = N_{i,\jc}^{(r+q)} - N_{i,\jc}^{(r)}
\end{aligned} \right.
$$
For all $C_1>0$, if there exists $C_2$ sufficiently large (depend on $C_1$ but not on $r$), such that $C_2 \sqrt{r\log \log r} \le d_{b,\jc}^{(r,s_r)}$ holds for infinitely many $r$'s, then for such sufficiently large $r$, there exists another sub-optimal design $i\neq k\neq b$ and a $u \le s_r$, $i$ is sampled at $r+u$ and 
\begin{equation}
    \left(1+C_1\logrr\right)\frac{N_{i,\jc}^{(r)}}{N_{b,\jc}^{(r)}} \le \frac{N_{i,\jc}^{(r+u)}}{N_{b,\jc}^{(r+s_r)}} \le \frac{N_{i,\jc}^{(r+u)}}{N_{b,\jc}^{(r+u)}}
\end{equation}
holds \text{ almost surely }.
\end{lemma}
\textbf{Proof}
By Lemma \ref{lem:length}, $t_r = O(r)$, which implies there exists $C_0 >0$, $d_{b,\jc}^{(r+s_r)} \le C_0 r$. Hence, for any fixed $C_1$, there exists $C_2$ such that for infinitely many $r$'s $C_2 \sqrt{r\log \log r}  \le d_{b,\jc}^{(r,s_r)} \le C_0 r$. Let $ \Delta_{\jc}^{(r)} =\left(\frac{\alpha_{b,\jc}^{(r)}}{\hat{\sigma}^{(r)}_{b,\jc}}\right)^2 - \sum\limits_{i\neq b} \left(\frac{\alpha_{i,\jc}^\pl}{\hat{\sigma}^{(r)}_{i,\jc}}\right)^2$. By the definition of $s_r$, $\Delta_{\jc}^{(r+s_r)} \ge -\Opar{\logrr}$ and $\Delta_{\jc}^{(r+s_r+1)} =  \Delta_{\jc}^{(r+s'_r)} + O\left(\logrr\right)$. Since $(b,\jc)$ is sampled at $r+s'_r$, $ \Delta_{\jc}^{(r+s'_r)} \le \Opar{\logrr}$, then there exists $C_3> C_3' > 0$, 
\begin{equation*}
    \Delta_{\jc}^{(r+s_r+1)} \le  \Delta_{\jc}^{(r+s'_r)} +C'_3\logrr  \le C_3\sqrt{\frac{\log\log r}{r}}.
\end{equation*}
Then one can choose $C_4>C_3$, 
\begin{equation}\label{eq:intermediate_lemma_1}
    |\Delta_{\jc}^{(r+s_r)}| \le  C_4\sqrt{\frac{\log\log r}{r}}
\end{equation}
Since $N_{k,\jc}^{(r+s_r)} - N_{k,\jc}^{(r)} =1$,  \eqref{eq:intermediate_lemma_1} implies
\begin{equation*}
    \left( \frac{N_{b,\jc}^{(r+s_r)}/\hat{\sigma}^{(r+s_r)}_{b,\jc}}{N^{(r+s_r)}}\right)^2 - \sum_{i\neq b\neq k} \left( \frac{N_{i,\jc}^{(r+s_r)}/\hat{\sigma}^{(r+s_r)}_{i,\jc}}{N^{(r+s_r)}}\right)^2 -\left( \frac{(N_{k,\jc}^{(r)}+1)/\hat{\sigma}^{(r+s_r)}_{k,\jc}}{N^{(r+s_r)}}\right)^2 \le C_4\sqrt{\frac{\log\log r}{r}}.
\end{equation*}
By some simple algebraic calculation we get
\begin{equation*}
    \sum_{i\neq b\neq k}\left( \frac{N_{i,\jc}^{(r+s_r)}/\hat{\sigma}^{(r+s_r)}_{i,\jc}}{N_{b,\jc}^{(r+s_r)}/\hat{\sigma}^{(r+s_r)}_{b,\jc}}  \right)^2 + \left( \frac{(N_{k,\jc}^{(r)}+1)/\hat{\sigma}^{(r+s_r)}_{k,\jc}}{N_{b,\jc}^{(r+s_r)}/\hat{\sigma}^{(r+s_r)}_{b,\jc}}  \right)^2 + C_4\sqrt{\frac{\log\log r}{r}} \left( \frac{N^{(r+s_r)}}{N_{b,\jc}^{(r+s_r)}/\hat{\sigma}^{(r+s_r)}_{b,\jc}}  \right)^2 \ge 1. 
\end{equation*}
By Lemma \ref{lem:ratio3}, there exists $C_5>0$, $C_5 N_{b,\jc}^{(r)}/\hat{\sigma}^{(r)}_{b,\jc} >  N^{(r)}$ for all large $r$. Hence,
\begin{equation*}
    \sum_{i\neq b\neq k}\left( \frac{N_{i,\jc}^{(r+s_r)}/\hat{\sigma}^{(r+s_r)}_{i,\jc}}{N_{b,\jc}^{(r+s_r)}/\hat{\sigma}^{(r+s_r)}_{b,\jc}}  \right)^2 \ge 1 - \left( \frac{(N_{k,\jc}^{(r)+1})/\hat{\sigma}^{(r+s_r)}_{k,\jc}}{N_{b,\jc}^{(r+s_r)}/\hat{\sigma}^{(r+s_r)}_{b,\jc}}  \right)^2 - C_4 C^2_5\sqrt{\frac{\log\log r}{r}}
\end{equation*}
Furthermore, there exists $C_6',C_6,,C_7>0$, such that
\begin{align}
     & \sum_{i\neq b\neq k}\left( \frac{N_{i,\jc}^{(r+s_r)}/\hat{\sigma}^{(r+s_r)}_{i,\jc}}{N_{b,\jc}^{(r+s_r)}/\hat{\sigma}^{(r+s_r)}_{b,\jc}}  \right)^2 -\sum_{i\neq b\neq k}\left( \frac{N_{i,\jc}^{(r)}/\hat{\sigma}^{(r+s_r)}_{i,\jc}}{N_{b,\jc}^{(r)}/\hat{\sigma}^{(r+s_r)}_{b,\jc}}  \right)^2 \notag \\
    \ge& 1 - \sum_{i\neq b\neq k}\left( \frac{N_{i,\jc}^{(r)}/\hat{\sigma}^{(r+s_r)}_{i,\jc}}{N_{b,\jc}^{(r)}/\hat{\sigma}^{(r+s_r)}_{b,\jc}}  \right)^2 - \left( \frac{(N_{k,\jc}^{(r)}+1)/\hat{\sigma}^{(r+s_r)}_{k,\jc}}{N_{b,\jc}^{(r+s_r)}/\hat{\sigma}^{(r+s_r)}_{b,\jc}}  \right)^2 -C_4 C^2_5\sqrt{\frac{\log\log r}{r}} \notag\\
    \ge & 1 - \sum_{i\neq b\neq k}\left( \frac{N_{i,\jc}^{(r)}/\hat{\sigma}^{(r)}_{i,\jc}}{N_{b,\jc}^{(r)}/\hat{\sigma}^{(r)}_{b,\jc}}  \right)^2 - \left( \frac{(N_{k,\jc}^{(r)}+1)/\hat{\sigma}^{(r+s_r)}_{k,\jc}}{N_{b,\jc}^{(r+s_r)}/\hat{\sigma}^{(r+s_r)}_{b,\jc}}  \right)^2 -(C_4 C^2_5+C_6')\sqrt{\frac{\log\log r}{r}} \notag\\
    \ge &  \left( \frac{(N_{k,\jc}^{(r)})/\hat{\sigma}^{(r)}_{k,\jc}}{N_{b,\jc}^{(r)}/\hat{\sigma}^{(r)}_{b,\jc}}  \right)^2 -  \left( \frac{(N_{k,\jc}^{(r)}+1)/\hat{\sigma}^{(r+s_r)}_{k,\jc}}{N_{b,\jc}^{(r+s_r)}/\hat{\sigma}^{(r+s_r)}_{b,\jc}}  \right)^2 - (C_4 C^2_5+C_6)\sqrt{\frac{\log\log r}{r}} \label{eq:intermediate_lemma_2} \\
    \ge &  \left( \frac{(N_{k,\jc}^{(r)})/\hat{\sigma}^{(r+s_r)}_{k,\jc}}{N_{b,\jc}^{(r)}/\hat{\sigma}^{(r+s_r)}_{b,\jc}}  \right)^2 -  \left( \frac{(N_{k,\jc}^{(r)}+1)/\hat{\sigma}^{(r+s_r)}_{k,\jc}}{N_{b,\jc}^{(r+s_r)}/\hat{\sigma}^{(r+s_r)}_{b,\jc}}  \right)^2 - (C_4 C^2_5+C_7)\sqrt{\frac{\log\log r}{r}} \notag \\
    = & \left(  \frac{\hat{\sigma}^{(r+s_r)}_{b,\jc}}{\hat{\sigma}^{(r+s_r)}_{k,\jc}} \right) ^2 \left[\frac{(N_{k,\jc}^{(r)})^2 (N_{b,\jc}^{(r+s_r)})^2 - (N_{k,\jc}^{(r)} + 1)^2(N_{b,\jc}^{(r)})^2}{(N_{b,\jc}^{(r)})^2 (N_{b,\jc}^{(r+s_r)})^2} \right] - (C_4 C^2_5+C_7)\sqrt{\frac{\log\log r}{r}} \notag \\
    =& \left(  \frac{\hat{\sigma}^{(r+s_r)}_{b,\jc}}{\hat{\sigma}^{(r+s_r)}_{k,\jc}} \right) ^2 \left[\frac{(N_{k,\jc}^{(r)})^2 (N_{b,\jc}^{(r)}+d_{b,\jc}^{(r,s_r)})^2 - (N_{k,\jc}^{(r)} + 1)^2(N_{b,\jc}^{(r)})^2}{(N_{b,\jc}^{(r)})^2 (N_{b,\jc}^{(r+s_r)})^2} \right] - (C_4 C^2_5+C_7)\sqrt{\frac{\log\log r}{r}} \notag\\
    =& \left(  \frac{\hat{\sigma}^{(r+s_r)}_{b,\jc}}{\hat{\sigma}^{(r+s_r)}_{k,\jc}} \right) ^2 \left[\frac{ (2N_{b,\jc}^{(r)}N_{k,\jc}^{(r)}+N_{k,\jc}^{(r)}d_{b,\jc}^{(r,s_r)}+N_{b,\jc}^{(r)}) (N_{k,\jc}^{(r)}d_{b,\jc}^{(r,s_r)} -N_{b,\jc}^{(r)})}{(N_{b,\jc}^{(r)})^2 (N_{b,\jc}^{(r+s_r)})^2} \right] - (C_4 C^2_5+C_7)\sqrt{\frac{\log\log r}{r}} \label{eq:intermediate_lemma_3}.
\end{align}
\eqref{eq:intermediate_lemma_2} holds since $(k,\jc)$ is sampled at $r$. By Lemma \ref{lem:ratio3} and $s_r = O(r)$, there exists $C_8, C_9, C_{10}>0$, $\frac{N_{k,\jc}^{(r)}}{N_{b,\jc}^{(r)}} \ge C_8$, $N_{b,\jc}^{(r)} \ge C_9  r$ and $N_{b,\jc}^{(r+s_r)} \le C_{10} r$. We have \eqref{eq:intermediate_lemma_3} is lower bounded by
\begin{align}
    &\left(  \frac{\hat{\sigma}^{(r+s_r)}_{b,\jc}}{\hat{\sigma}^{(r+s_r)}_{k,\jc}} \right) ^2 \left[\frac{C_8^2(2C_9 r + d_{b,\jc}^{(r,s_r)}+1/C_8)(d_{b,\jc}^{(r,s_r)} -1/C_8)}{C_{10}^2 r^2} \right] - (C_4 C^2_5+C_7)\sqrt{\frac{\log\log r}{r}} \notag \\
    = &  \left(  \frac{\hat{\sigma}^{(r+s_r)}_{b,\jc}}{\hat{\sigma}^{(r+s_r)}_{k,\jc}} \right) ^2 \left(\frac{C_8}{C_{10}}\right)^2 \left[  \frac{2C_9d_{b,\jc}^{(r,s_r)}}{r} +\left(\frac{d_{b,\jc}^{(r,s_r)}}{r}\right)^2 -\frac{2C_9}{C_8r} -\frac{1}{C_8^2r^2}\right] - (C_4 C^2_5+ C_7)\sqrt{\frac{\log\log r}{r}}\notag\\
    \ge & \left(  \frac{\hat{\sigma}^{(r+s_r)}_{b,\jc}}{\hat{\sigma}^{(r+s_r)}_{k,\jc}} \right) ^2 \left(\frac{C_8}{C_{10}}\right)^2 2C_9 \left[  \frac{d_{b,\jc}^{(r,s_r)}}{r}  -\frac{1}{C_8r} -\frac{1}{2C_9C_8^2r^2}\right] - (C_4 C^2_5+C_7)\sqrt{\frac{\log\log r}{r}} \label{eq:intermediate_lemma_4}
\end{align}
There exists $C_{11} > 0 $, $\frac{1}{C_8r} + \frac{1}{2C_9C_8^2r^2} \le C_{11} \sqrt{\frac{\log\log r}{r}}$. Choose $ 0 < C_{12} \le \left(  \frac{\hat{\sigma}^{(r+s_r)}_{b,\jc}}{\hat{\sigma}^{(r+s_r)}_{k,\jc}} \right) ^2 \left(\frac{C_8}{C_{10}}\right)^2 2C_9 $ for all large $r$. Furthermore since $d_{b,\jc}^{(r,s_r)} \ge C_2 \sqrt{r\log \log r}$, there exists $C_{13} >0$
\begin{align}
    \eqref{eq:intermediate_lemma_4} &\ge C_{12} \left[ (C_2 - C_{11}) \sqrt{\frac{\log\log r}{r}}\right] - (C_4 C^2_5+C_7)\sqrt{\frac{\log\log r}{r}} \notag \\
    & = [C_{12}(C_2-C_{11}) - (C_4C_5^2+C_7)] \sqrt{\frac{\log\log r}{r}} \notag \\
    & \ge C_{13} C_2  \sqrt{\frac{\log\log r}{r}}\label{eq:intermediate_lemma_5}
\end{align}
\eqref{eq:intermediate_lemma_5} holds for $C_2$ large enough (but not depends on $r$). For example, we can choose  $C_{13}= C_{12}/2$ and then \eqref{eq:intermediate_lemma_5} holds for all $C_2 \ge 2C_{11}+ \frac{2C_4C_5^2+C_7}{C_{12}}$. Since $C_3,C_4,\ldots,C_{13}$ are all independent of $r$, the $C_2$ here is also independent of $r$. Then,
\begin{align}
    &\sum_{i\neq b\neq k}\left( \frac{N_{i,\jc}^{(r+s_r)}/\hat{\sigma}^{(r+s_r)}_{i,\jc}}{N_{b,\jc}^{(r+s_r)}/\hat{\sigma}^{(r+s_r)}_{b,\jc}}  \right)^2 -\sum_{i\neq b\neq k}\left( \frac{N_{i,\jc}^{(r)}/\hat{\sigma}^{(r+s_r)}_{i,\jc}}{N_{b,\jc}^{(r)}/\hat{\sigma}^{(r+s_r)}_{b,\jc}}  \right)^2\ge C_{13}C_2 \sqrt{\frac{\log\log r}{r}}\notag
\end{align}
There exists a non-optimal design $h \neq k\neq b$, such that 
 \begin{equation*}
     \left( \frac{N_{h,\jc}^{(r+s_r)}/\hat{\sigma}^{(r+s_r)}_{h,\jc}}{N_{b,\jc}^{(r+s_r)}/\hat{\sigma}^{(r+s_r)}_{b,\jc}}  \right)^2 - \left( \frac{N_{h,\jc}^{(r)}/\hat{\sigma}^{(r+s_r)}_{h,\jc}}{N_{b,\jc}^{(r)}/\hat{\sigma}^{(r+s_r)}_{b,\jc}}  \right)^2 \ge \frac{1}{K-2}C_{13}C_{2}\sqrt{\frac{\log\log r}{r}}.
 \end{equation*}
 Or equivalently,
 \begin{equation*}
     \left(\frac{N_{h,\jc}^{(r+s_r)}/N_{b,\jc}^{(r+s_r)}}{N_{h,\jc}^{(r)}/N_{b,\jc}^{(r)}}\right)^2 \ge 1+ \frac{C_{13}C_2}{K-2}\left(\frac{\hat{\sigma}^{(r+s_r)}_{h,\jc}}{\hat{\sigma}^{(r+s_r)}_{b,\jc}}\right)^2 \left(\frac{N_{b,\jc}^{(r)}}{N_{h,\jc}^{(r)}}\right)^2  \sqrt{\frac{\log\log r}{r}}.
 \end{equation*}
 By Lemma \ref{lem:ratio3} and the convergence of the sample variance, there exists $0 < C_{14} <\frac{C_{13}}{K-2}\left(\frac{\hat{\sigma}^{(r+s_r)}_{h,\jc}}{\hat{\sigma}^{(r+s_r)}_{b,\jc}}\right)^2 \left(\frac{N_{b,\jc}^{(r)}}{N_{h,\jc}^{(r)}}\right)^2 $ for all large $r$. We obtain
 \begin{equation*}
      \left(\frac{N_{h,\jc}^{(r+s_r)}/N_{b,\jc}^{(r+s_r)}}{N_{h,\jc}^{(r)}/N_{b,\jc}^{(r)}}\right)^2 \ge 1+C_{14}C_2 \sqrt{\frac{\log\log r}{r}}
 \end{equation*}
 Hence,
 \begin{align}
   \frac{N_{h,\jc}^{(r+s_r)}}{N_{b,\jc}^{(r+s_r)}} \ge \sqrt{1+C_{14}C_2\sqrt{\frac{\log\log r}{r}} } \frac{N_{h,\jc}^{(r)}}{N_{b,\jc}^{(r)}}  \ge (1+\frac{C_{14}C_2}{4}\sqrt{\frac{\log\log r}{r}})\frac{N_{h,\jc}^{(r)}}{N_{b,\jc}^{(r)}}
 \end{align}
 for all large $r$ by Taylor Expansion. Let $v = \sup\{\ell\le s_r : \mathbf{1}_{h,\jc}^{(r+\ell)} =1\}$. $r+v$ is the last time before $r+s_r$ at which $h$ is sampled. We then have
 \begin{align}
     \frac{N_{h,\jc}^{(r+v)}}{N_{b,\jc}^{(r+v)}} &\ge \frac{N_{h,\jc}^{(r+v)}}{N_{b,\jc}^{(r+s_r)}} \notag \\
     & = \frac{N_{h,\jc}^{(r+s_r)}-1}{N_{b,\jc}^{(r+s_r)}} \notag \\
     & = (1 - \frac{1}{N_{b,\jc}^{(r+s_r)}}) \frac{N_{h,\jc}^{(r+s_r)}}{N_{b,\jc}^{(r+s_r)}} \notag \\
     &\ge (1 - \frac{1}{N_{b,\jc}^{(r)}})(1+\frac{C_{14}C_2}{4}\sqrt{\frac{\log\log r}{r}})\frac{N_{h,\jc}^{(r)}}{N_{b,\jc}^{(r)}} \notag \\
     &\ge (1 - \frac{1}{C_9 r})(1+\frac{C_{14}C_2}{4}\sqrt{\frac{\log\log r}{r}})\frac{N_{h,\jc}^{(r)}}{N_{b,\jc}^{(r)}} \notag \\
     & \ge (1+ \frac{C_{14}C_2}{8}\sqrt{\frac{\log\log r}{r}} )\frac{N_{h,\jc}^{(r)}}{N_{b,\jc}^{(r)}} \notag 
 \end{align}
 for all large $r$. Then, given any $C_1$, there exists $C_2 \ge \max\{ \frac{8 C_1}{C_{14}}, 2C_{11}+ \frac{2C_4C_5^2+C_7}{C_{12}}\}$. The Lemma holds for true for the $C_1,C_2$. \hfill $\blacksquare$

\begin{lemma} \label{lem:NofS_ib}
For a fixed sub-optimal design $k$, between two samples of $(k,\jc)$. Suppose $(k,\jc)$ is sampled at $r$ and let $\ell_r = \inf\{ \ell>0: \mathbf{1}_{(k,\jc)}=1\}$. $r+\ell_r$ is the next time $(k,\jc)$ being sampled. Then the number of samples that can be allocated to $(b,\jc)$ between $r$ and $r+\ell_r$ is  $O(\sqrt{r\log \log r})$ \text{ almost surely }.
\end{lemma}
\textbf{Proof}
We use the same notation of $s'_{r}, s_r, d_{i,j}^{(r,q)}$ as in Lemma \ref{lem:intermediate_lemma}. Since $d_{b,\jc}^{(r,\ell_r)} - d_{b,\jc}^{(r,s_r)} = d_{b,\jc}^{(r,s'_r)} - d_{b,\jc}^{(r,s_r)} +1 = O(\sqrt{r\log \log r})$ by Lemma \ref{lem:NofS_banyi}, it is sufficient to prove $d_{b,\jc}^{(r,s_r)} = O(\sqrt{r\log \log r})$. Prove by contradiction. Suppose the statement does not hold. Then $\forall C_2 >0$, there exists $r$ such that $d_{b,\jc}^{(r,s_r)} \ge C_2 \sqrt{r\log \log r}$.  By Lemma \ref{lem:intermediate_lemma}, $\forall C_1 >0 $ (remain to be specified), there exists an iteration $r$ at which $(k,\jc)$ is sampled, another non-optimal design $h \neq k$ and an iteration $v < s_r$, such that $(h,\jc)$ is sampled at $v$ and 
\begin{equation*}
    \frac{\alpha_{h,\jc}^{(r+v)}}{\alpha_{b,\jc}^{(r+s_r)}} \ge (1+C_1\sqrt{\frac{\log\log r}{r}})\frac{\alpha_{h,\jc}^{(r)}}{\alpha_{b,\jc}^{(r)}}
\end{equation*}
holds. We aim to show $(h,\jc)$ cannot be sampled at $r+v$ for a contradiction. It is sufficient to show
\begin{equation}
    \frac{(\hat{\mu}_b^{(r+v)} - \hat{\mu}_h^{(r+v)} )^2}{\sum\limits_{j=1}^D \frac{(\hat{\sigma}^{(r+v)}_{b,j})^2 (\pi_j^{(t_{r+v})})^2}{\alpha_{b,j}^{(r+v)}} + \sum\limits_{j=1}^D \frac{(\hat{\sigma}^{(r+v)}_{h,j})^2(\pi_j^{(t_{r+v})})^2}{\alpha_{h,j}^{(r+v)}}} > \frac{(\hat{\mu}_b^{(r+v)} - \hat{\mu}_k^{(r+v)} )^2}{\sum\limits_{j=1}^D \frac{(\hat{\sigma}^{(r+v)}_{b,j})^2 (\pi_j^{(t_{r+v})})^2}{\alpha_{b,j}^{(r+v)}} + \sum\limits_{j=1}^D \frac{(\hat{\sigma}^{(r+v)}_{k,j})^2(\pi_j^{(t_{r+v})})^2}{\alpha_{k,j}^{(r+v)}}}.
\end{equation}
Denote by $\delta_i^\pl = (\hat{\mu}_b^\pl - \hat{\mu}_i^\pl)^2 $. It is equivalent to show 
\begin{equation}\label{eq:NofS_ib_1}
\begin{aligned}
    \delta_h^{(r+v)} \left( \sum\limits_{j=1}^D \frac{(\hat{\sigma}^{(r+v)}_{b,j})^2 (\pi_j^{(t_{r+v})})^2}{\alpha_{b,j}^{(r+v)}} + \sum\limits_{j=1}^D \frac{(\hat{\sigma}^{(r+v)}_{k,j})^2(\pi_j^{(t_{r+v})})^2}{\alpha_{k,j}^{(r+v)}}\right) \\
    > \delta_k^{(r+v)} \left(\sum\limits_{j=1}^D \frac{(\hat{\sigma}^{(r+v)}_{b,j})^2 (\pi_j^{(t_{r+v})})^2}{\alpha_{b,j}^{(r+v)}} + \sum\limits_{j=1}^D \frac{(\hat{\sigma}^{(r+v)}_{h,j})^2(\pi_j^{(t_{r+v})})^2}{\alpha_{h,j}^{(r+v)}}\right)
\end{aligned}
\end{equation}

By Lemma \ref{lem: finite j}, we have $ \frac{\alpha_{i,j}^{(r)}}{\hat{\sigma}^{(r)}_{i,j}\pi_{j}^{(t_r)}}   =O(\sqrt{\frac{\log\log r}{r}})$ for $j\neq \jc$. Then, there exist $C_3 > 0$
\begin{equation} 
\text{ LHS of } \eqref{eq:NofS_ib_1} \ge \delta_h^{(r+v)} \left( \frac{(\hat{\sigma}^{(r+v)}_{b,\jc})^2 (\pi_\jc^{(t_{r+v})})^2}{\alpha_{b,\jc}^{(r+v)}}+ \frac{(\hat{\sigma}^{(r+v)}_{k,\jc})^2 (\pi_\jc^{(t_{r+v})})^2}{\alpha_{k,\jc}^{(r+v)}}\right)  - C_3 \logrr    \notag
\end{equation}
\begin{equation}
    \text{ RHS of } \eqref{eq:NofS_ib_1} \le \delta_k^{(r+v)} \left( \frac{(\hat{\sigma}^{(r+v)}_{b,\jc})^2 (\pi_\jc^{(t_{r+v})})^2}{\alpha_{b,\jc}^{(r+v)}}+ \frac{(\hat{\sigma}^{(r+v)}_{h,\jc})^2 (\pi_\jc^{(t_{r+v})})^2}{\alpha_{h,\jc}^{(r+v)}}\right)  + C_3 \logrr \notag
\end{equation}
Since $(k,\jc)$ is sampled at $r$, 

\begin{equation*}
    \delta_h^{(r)} \left( \sum\limits_{j=1}^D \frac{(\hat{\sigma}^{(r)}_{b,j})^2 (\pi_j^{(t_r)})^2}{\alpha_{b,j}^{(r)}} + \sum\limits_{j=1}^D \frac{(\hat{\sigma}^{(r)}_{k,j})^2(\pi_j^{(t_r)})^2}{\alpha_{k,j}^{(r)}}\right) \ge \delta_k^{(r)} \left(\sum\limits_{j=1}^D \frac{(\hat{\sigma}^{(r)}_{b,j})^2 (\pi_j^{(t_r)})^2}{\alpha_{b,j}^{(r)}} + \sum\limits_{j=1}^D \frac{(\hat{\sigma}^{(r)}_{h,j})^2(\pi_j^{(t_r)})^2}{\alpha_{h,j}^{(r)}}\right).
\end{equation*}
Then we obtain, 
\begin{equation}\label{eq:NofS_ib_2.5}
    \begin{aligned}
   &\delta_h^{(r+v)} \left( \frac{(\hat{\sigma}^{(r+v)}_{b,\jc})^2 (\pi_\jc^{(t_{r+v})})^2}{\alpha_{b,\jc}^{(r+v)}}+ \frac{(\hat{\sigma}^{(r+v)}_{k,\jc})^2 (\pi_\jc^{(t_{r+v})})^2}{\alpha_{k,\jc}^{(r+v)}}\right)   \\
   \ge& 
   \delta_k^{(r+v)} \left( \frac{(\hat{\sigma}^{(r+v)}_{b,\jc})^2 (\pi_\jc^{(t_{r+v})})^2}{\alpha_{b,\jc}^{(r+v)}}+ \frac{(\hat{\sigma}^{(r+v)}_{h,\jc})^2 (\pi_\jc^{(t_{r+v})})^2}{\alpha_{h,\jc}^{(r+v)}}\right)  - 2 C_3 \logrr
\end{aligned}
\end{equation}

Then, 
\begin{equation}\label{eq:NofS_ib_3}
    \begin{aligned}
    \text{ LHS of } \eqref{eq:NofS_ib_2.5}  = &\delta_h^{(r)}  \frac{(\hat{\sigma}^{(r)}_{b,\jc})^2(\pi_{\jc}^{(t_r)})^2}{\alpha_{b,\jc}^{(r)}} \cdot \frac{\delta_h^{(r+v)}}{\delta_h^{(r)}} \frac{(\hat{\sigma}^{(r+v)}_{b,\jc})^2}{(\hat{\sigma}^{(r)}_{b,\jc})^2}\frac{\alpha_{b,\jc}^{(r)}}{\alpha_{b,\jc}^{(r+v)}} \frac{(\pi_{\jc}^{(t_{r+v})})^2}{(\pi_{\jc}^{(t_r)})^2} \\
    &+
    \delta_h^{(r)}  \frac{(\hat{\sigma}^{(r)}_{k,\jc})^2(\pi_{\jc}^{(t_r)})^2}{\alpha_{k,\jc}^{(r)}} \cdot \frac{\delta_h^{(r+v)}}{\delta_h^{(r)}} \frac{(\hat{\sigma}^{(r+v)}_{k,\jc})^2}{(\hat{\sigma}^{(r)}_{k,\jc})^2}\frac{\alpha_{k,\jc}^{(r)}}{\alpha_{k,\jc}^{(r+v)}} \frac{(\pi_{\jc}^{(t_{r+v})})^2}{(\pi_{\jc}^{(t_r)})^2} - C_3\sqrt{\frac{\log\log r}{r}}.  
\end{aligned}
\end{equation}
and 
\begin{equation}\label{eq:NofS_ib_4}
    \begin{aligned}
    \text{ RHS of } \eqref{eq:NofS_ib_2.5}  = &\delta_k^{(r)}  \frac{(\hat{\sigma}^{(r)}_{b,\jc})^2(\pi_{\jc}^{(t_r)})^2}{\alpha_{b,\jc}^{(r)}} \cdot \frac{\delta_h^{(r+v)}}{\delta_h^{(r)}} \frac{(\hat{\sigma}^{(r+v)}_{b,\jc})^2}{(\hat{\sigma}^{(r)}_{b,\jc})^2}\frac{\alpha_{b,\jc}^{(r)}}{\alpha_{b,\jc}^{(r+v)}} \frac{(\pi_{\jc}^{(t_{r+v})})^2}{(\pi_{\jc}^{(t_r)})^2} \\
    &+
    \delta_h^{(r)}  \frac{(\hat{\sigma}^{(r)}_{k,\jc})^2(\pi_{\jc}^{(t_r)})^2}{\alpha_{k,\jc}^{(r)}} \cdot \frac{\delta_h^{(r+v)}}{\delta_h^{(r)}} \frac{(\hat{\sigma}^{(r+v)}_{k,\jc})^2}{(\hat{\sigma}^{(r)}_{k,\jc})^2}\frac{\alpha_{k,\jc}^{(r)}}{\alpha_{k,\jc}^{(r+v)}} \frac{(\pi_{\jc}^{(t_{r+v})})^2}{(\pi_{\jc}^{(t_r)})^2} - C_3\sqrt{\frac{\log\log r}{r}}.  
\end{aligned}
\end{equation}
$\forall 1\le i\le K$,  by Lemma \ref{lem:mu},  $\frac{\delta_i^{(r+v)}}{\delta_i^{(r)}} = 1+O(\sqrt{\frac{\log\log r}{r}})$; by Lemma \ref{lem: posterior convergence} and \ref{lem:variance}, $\frac{\pi_{\jc}^{(t_{r+v})}}{\pi_{\jc}^{(t_r)}} =1+O(\sqrt{\frac{\log\log r}{r}}) $, and $\frac{(\hat{\sigma}^{(r+v)}_{i,\jc})^2}{(\hat{\sigma}^{(r)}_{i,\jc})^2} = 1+O(\sqrt{\frac{\log\log r}{r}}) $; Furthermore $N_{k,\jc}^{(r+v)} = N_{k,\jc}^{(r)}+1$. Hence, there exist $C_5>0$, such that 
\begin{align}
    \eqref{eq:NofS_ib_3} \ge &\delta_h^{(r)}  \frac{(\hat{\sigma}^{(r)}_{b,\jc})^2(\pi_{\jc}^{(t_r)})^2}{\alpha_{b,\jc}^{(r)}} \frac{\alpha_{b,\jc}^{(r)}}{\alpha_{b,\jc}^{(r+v)}} (1-C_5\sqrt{\frac{\log\log r}{r}}) + \delta_h^{(r)}  \frac{(\hat{\sigma}^{(r)}_{k,\jc})^2(\pi_{\jc}^{(t_r)})^2}{\alpha_{k,\jc}^{(r)}}  (1-C_5\sqrt{\frac{\log\log r}{r}}) - C_3\sqrt{\frac{\log\log r}{r}} \label{eq:NofS_ib_5}
\end{align}
and 
\begin{align}
    \eqref{eq:NofS_ib_4} \le \delta_k^{(r)}  \frac{(\hat{\sigma}^{(r)}_{b,\jc})^2(\pi_{\jc}^{(t_r)})^2}{\alpha_{b,\jc}^{(r)}} \frac{\alpha_{b,\jc}^{(r)}}{\alpha_{b,\jc}^{(r+v)}} (1+C_5\sqrt{\frac{\log\log r}{r}})+
    \delta_k^{(r)}  \frac{(\hat{\sigma}^{(r)}_{h,\jc})^2(\pi_{\jc}^{(t_r)})^2}{\alpha_{h,\jc}^{(r)}}  \frac{\alpha_{h,\jc}^{(r)}}{\alpha_{h,\jc}^{(r+v)}}(1+C_5\sqrt{\frac{\log\log r}{r}}) + C_3\sqrt{\frac{\log\log r}{r}}\label{eq:NofS_ib_6}
\end{align}

Divide by $\frac{\alpha_{b,\jc}^{(r)}}{\alpha_{b,\jc}^{(r+v)}}$ by both sides and notice $1\le  \frac{\alpha_{b,\jc}^{(r+v)}}{\alpha_{b,\jc}^{(r)}} \le  \frac{1}{\alpha_{b,\jc}^{(r)}} \le C_7$ for some $C_7>0$ by Lemma \ref{lem:ratio4}. Then $\eqref{eq:NofS_ib_5} > \eqref{eq:NofS_ib_6}$ can be  implied by
\begin{equation} \label{eq:NofS_ib_7}
\begin{aligned}
    &\delta_h^{(r)}  \frac{(\hat{\sigma}^{(r)}_{b,\jc})^2(\pi_{\jc}^{(t_r)})^2}{\alpha_{b,\jc}^{(r)}}  +\delta_h^{(r)}  \frac{(\hat{\sigma}^{(r)}_{k,\jc})^2(\pi_{\jc}^{(t_r)})^2}{\alpha_{k,\jc}^{(r)}}  - C_8 \sqrt{\frac{\log\log r}{r}}\notag\\
    > &\delta_k^{(r)}  \frac{(\hat{\sigma}^{(r)}_{b,\jc})^2(\pi_{\jc}^{(t_r)})^2}{\alpha_{b,\jc}^{(r)}}+
     \delta_k^{(r)}  \frac{(\hat{\sigma}^{(r)}_{h,\jc})^2(\pi_{\jc}^{(t_r)})^2}{\alpha_{h,\jc}^{(r)}}  \left(\frac{\alpha_{h,\jc}^{(r)}}{\alpha_{b,\jc}^{(r)}}\middle /\frac{\alpha_{h,\jc}^{(r+v)}}{\alpha_{b,\jc}^{(r+v)}} \right)(1+C_5\sqrt{\frac{\log\log r}{r}}) + C_8\sqrt{\frac{\log\log r}{r}}\notag
\end{aligned}
\end{equation}
for some $C_8$ independent of $r$ (depends on $C_5,C_7,C_3$). Since $\left(\frac{\alpha_{h,\jc}^{(r)}}{\alpha_{b,\jc}^{(r)}}\middle /\frac{\alpha_{h,\jc}^{(r+v)}}{\alpha_{b,\jc}^{(r+v)}} \right) < \frac{1}{1+C_1\sqrt{\frac{\log\log r}{r}} } < 1 - \frac{C_1}{2}\sqrt{\frac{\log\log r}{r}} $ for all large $r$'s, \eqref{eq:NofS_ib_7} can be further implied by 
\begin{equation}\label{eq:NofS_ib_8}
\begin{aligned}
    &\delta_h^{(r)}  \frac{(\hat{\sigma}^{(r)}_{b,\jc})^2(\pi_{\jc}^{(t_r)})^2}{\alpha_{b,\jc}^{(r)}}  +\delta_h^{(r)}  \frac{(\hat{\sigma}^{(r)}_{k,\jc})^2(\pi_{\jc}^{(t_r)})^2}{\alpha_{k,\jc}^{(r)}}  - C_8 \sqrt{\frac{\log\log r}{r}}\\
    > & \delta_k^{(r)}  \frac{(\hat{\sigma}^{(r)}_{b,\jc})^2(\pi_{\jc}^{(t_r)})^2}{\alpha_{b,\jc}^{(r)}}+
     \delta_k^{(r)}  \frac{(\hat{\sigma}^{(r)}_{h,\jc})^2(\pi_{\jc}^{(t_r)})^2}{\alpha_{h,\jc}^{(r)}} (1- \frac{C_1}{2}\sqrt{\frac{\log\log r}{r}})(1+C_5\sqrt{\frac{\log\log r}{r}}) + C_8\sqrt{\frac{\log\log r}{r}}\\
    =  &\delta_k^{(r)}  \frac{(\hat{\sigma}^{(r)}_{b,\jc})^2(\pi_{\jc}^{(t_r)})^2}{\alpha_{b,\jc}^{(r)}}+
     \delta_k^{(r)}  \frac{(\hat{\sigma}^{(r)}_{h,\jc})^2(\pi_{\jc}^{(t_r)})^2}{\alpha_{h,\jc}^{(r)}}  + \left[C_8+ C_5   \delta_k^{(r)}  \frac{(\hat{\sigma}^{(r)}_{h,\jc})^2(\pi_{\jc}^{(t_r)})^2}{\alpha_{h,\jc}^{(r)}} -\frac{C_1}{2}   \delta_k^{(r)}  \frac{(\hat{\sigma}^{(r)}_{h,\jc})^2(\pi_{\jc}^{(t_r)})^2}{\alpha_{h,\jc}^{(r)}}\right]\sqrt{\frac{\log\log r}{r}}
\end{aligned}
\end{equation}
Let $C_9,C_{10}$ be two constants such that $ C_9 \le    \delta_k^{(r)}  \frac{(\hat{\sigma}^{(r)}_{h,\jc})^2(\pi_{\jc}^{(t_r)})^2}{\alpha_{h,\jc}^{(r)}} \le C_{10}$. Then by \eqref{eq:NofS_ib_2.5}, \eqref{eq:NofS_ib_8} is implied by
\begin{equation}\label{eq:NofS_ib_9}
    \left[2C_8+ 2C_3+ C_5 C_{10} -\frac{C_1}{2} C_{9}\right]\sqrt{\frac{\log\log r}{r}} <0.
\end{equation}
Since $C_3,C_5,C_8,C_9,C_{10}$ are all independent of $r$, we can choose $C_1$ to be large enough such that \eqref{eq:NofS_ib_9} holds, which gives us the contradiction that $(h,\jc)$ cannot be sampled at $r+v$. \hfill $\blacksquare$

We can now finally prove the ``Local Balance" optimality condition.

\textbf{Proof of \eqref{thmeq:error local balance} in Theorem \ref{thm: convergence speed}.2.}
\textbf{Proof.}
Denote by $\delta_i^\pl = (\hat{\mu}_b^\pl - \hat{\mu}_i^\pl)^2$ and $\Theta_i^\pl = \frac{\delta_i^\pl}{\sum\limits_{j=1}^D \frac{(\hat{\sigma}^{(r)}_{i,j})^2 (\pi_j^\tl)^2}{\alpha_{i,j}^\pl} + \sum\limits_{j=1}^D \frac{(\hat{\sigma}^{(r)}_{i,j})^2(\pi_j^\tl)^2}{\alpha_{i,j}^\pl}} $. For any two non-optimal designs $i\neq k$ and an iteration $r$. Without loss of generality, assume $\Theta_i^{(r)} \le \Theta_k^{(r)}$. Let $u= \sup\{\ell<r: \mathbf{1}_{k,\jc}^\pl = 1\}$ be the last time $k$ is sampled before $r$. $r = O(u)$ by Lemma \ref{lem:length}. Also notice $\frac{(\hat{\sigma}_{i,j}^\pl)^2(\pi_j^\tl)^2}{\alpha_{i,j}^\pl} = \Opar{\logll}$ for $j\neq \jc$, we have 
\begin{align}
    0 &\le \Theta_k^{(r)} -\Theta_i^{(r)} \notag\\
    &\le \frac{\delta_k^{(r)}}{\frac{(\hat{\sigma}^{(r)}_{b,\jc})^2(\pi_{\jc}^{(t_r)})^2}{\alpha_{b,\jc}^{(r)}}  + \frac{(\hat{\sigma}^{(r)}_{k,\jc})^2(\pi_{\jc}^{(t_r)})^2}{\alpha_{k,\jc}^{(r)}} - \Opar{\sqrt{\frac{\log\log r}{r}} }} - \frac{\delta_i^{(r)}}{\frac{(\hat{\sigma}^{(r)}_{b,\jc})^2(\pi_{\jc}^{(t_r)})^2}{\alpha_{b,\jc}^{(r)}}  + \frac{(\hat{\sigma}^{(r)}_{i,\jc})^2(\pi_{\jc}^{(t_r)})^2}{\alpha_{i,\jc}^{(r)}}  + \Opar{\sqrt{\frac{\log\log r}{r}}} } \\
    & \le  \frac{\delta_k^{(r)}}{\frac{(\hat{\sigma}^{(r)}_{b,\jc})^2(\pi_{\jc}^{(t_r)})^2}{\alpha_{b,\jc}^{(r)}}  + \frac{(\hat{\sigma}^{(r)}_{k,\jc})^2(\pi_{\jc}^{(t_r)})^2}{\alpha_{k,\jc}^{(r)}}} - \frac{\delta_i^{(r)}}{\frac{(\hat{\sigma}^{(r)}_{b,\jc})^2(\pi_{\jc}^{(t_r)})^2}{\alpha_{b,\jc}^{(r)}}  + \frac{(\hat{\sigma}^{(r)}_{i,\jc})^2(\pi_{\jc}^{(t_r)})^2}{\alpha_{i,\jc}^{(r)}}  } +  \Opar{\sqrt{\frac{\log\log r}{r}}} \label{eq:localprove_2}
 \end{align}
By Lemma \eqref{lem:mu}, $\delta_{i'}^{(r)} = \delta_{i'}^{(u)} + O(\sqrt{\frac{\log\log u}{u}}) = \delta_{i'}^{(u)}+ O(\sqrt{\frac{\log\log r}{r}}), \ \forall i'$; by Lemma \ref{lem: posterior convergence} and Lemma \ref{lem:variance}, $\pi_{\jc}^{(t_r)} =1 - O(\sqrt{\frac{\log\log r}{r}})$, $\hat{\sigma}_{i,j}^{(r)} = \hat{\sigma}_{i,j}^{(u)} + O(\sqrt{\frac{\log\log r}{r}})$. Furthermore $\alpha_{i',j'}^\pl$ is upper and lower bounded by positive constants $\forall i',j^c$. Then,
\begin{align}
  \eqref{eq:localprove_2} \le   \frac{\delta_k^{(u)}}{\frac{(\hat{\sigma}^{(u)}_{b,\jc})^2(\pi_{\jc}^{(t_u)})^2}{\alpha_{b,\jc}^{(r)}}  + \frac{(\hat{\sigma}^{(u)}_{k,\jc})^2(\pi_{\jc}^{(t_u)})^2}{\alpha_{k,\jc}^{(r)}}} - \frac{\delta_i^{(u)}}{\frac{(\hat{\sigma}^{(u)}_{b,\jc})^2 (\pi_{\jc}^{(t_u)})^2}{\alpha_{b,\jc}^{(r)}} + \frac{(\hat{\sigma}^{(u)}_{i,\jc})^2 (\pi_{\jc}^{(t_u)})^2}{\alpha_{i,\jc}^{(r)}} } + \Opar{ \sqrt{\frac{\log\log r}{r}}}\label{eq:localprove_3}.
\end{align}
By Lemma \ref{lem:NofS_ib}, $N_{b,\jc}^{(r)} = N_{b,j}^{(u)} + O(\sqrt{u \log u}) = N_{b,j}^{(u)} + O(\sqrt{r\log \log r})$. Hence,
\begin{align}
    \eqref{eq:localprove_3} &\le  \frac{1}{N^{(r)}} \frac{\delta_k^{(u)}}{\frac{(\hat{\sigma}^{(u)}_{b,\jc})^2(\pi_{\jc}^{(t_u)})^2}{N_{b,\jc}^{(u)} + \Opar{\sqrt{r\log \log r}}}  + \frac{(\hat{\sigma}^{(u)}_{k,\jc})^2 (\pi_{\jc}^{(t_u)})^2}{N_{k,\jc}^{(u)} +1 } } -\frac{1}{N^{(r)}}\frac{\delta_i^{(u)}}{\frac{(\hat{\sigma}^{(u)}_{b,\jc})^2 (\pi_{\jc}^{(t_u)})^2}{N_{b,\jc}^{(r)}} + \frac{(\hat{\sigma}^{(u)}_{i,\jc})^2 (\pi_{\jc}^{(t_u)})^2}{N_{i,\jc}^{(r)}} } + \Opar{ \sqrt{\frac{\log\log r}{r}}} \label{eq:localprove_4}
\end{align}
Since $r = O(u)$ and by Lemma \ref{lem:ratio4}, $ \lim\inf_{r\rightarrow}\frac{N_{i,\jc}^{(u)}}{N^{(r)}} >0 \ \forall i$. Then,
\begin{align}
\eqref{eq:localprove_4} &=   \frac{\delta_k^{(u)}}{\frac{(\hat{\sigma}^{(u)}_{b,\jc})^2(\pi_{\jc}^{(t_u)})^2}{\frac{N_{b,\jc}^{(u)}}{N^{(r)}} + \Opar{\logrr}}  + \frac{(\hat{\sigma}^{(u)}_{k,\jc})^2 (\pi_{\jc}^{(t_u)})^2}{\frac{N_{k,\jc}^{(u)} +1}{N^{(r)}} } } -\frac{1}{N^{(r)}}\frac{\delta_i^{(u)}}{\frac{(\hat{\sigma}^{(u)}_{b,\jc})^2 (\pi_{\jc}^{(t_u)})^2}{N_{b,\jc}^{(r)}} + \frac{(\hat{\sigma}^{(u)}_{i,\jc})^2 (\pi_{\jc}^{(t_u)})^2}{N_{i,\jc}^{(r)}} }+ \Opar{ \sqrt{\frac{\log\log r}{r}}} \notag \\
    &\le \frac{1}{N^{(r)}} \frac{\delta_k^{(u)}}{\frac{(\hat{\sigma}^{(u)}_{b,\jc})^2(\pi_{\jc}^{(t_u)})^2}{N_{b,\jc}^{(u)} }  + \frac{(\hat{\sigma}^{(u)}_{k,\jc})^2 (\pi_{\jc}^{(t_u)})^2}{N_{k,\jc}^{(u)}  } } -\frac{1}{N^{(r)}}\frac{\delta_i^{(u)}}{\frac{(\hat{\sigma}^{(u)}_{b,\jc})^2 (\pi_{\jc}^{(t_u)})^2}{N_{b,\jc}^{(r)}} + \frac{(\hat{\sigma}^{(u)}_{i,\jc})^2 (\pi_{\jc}^{(t_u)})^2}{N_{i,\jc}^{(r)}} } + \Opar{\logrr} \notag\\
    & \le \frac{1}{N^{(r)}} \frac{\delta_k^{(u)}}{\sum_{j=1}^D \frac{(\hat{\sigma}^{(u)}_{b,j})^2(\pi_{j}^{(t_u)})^2}{N_{b,j}^{(u)} }  + \sum_{j=1}^D \frac{(\hat{\sigma}^{(u)}_{k,j})^2 (\pi_{j}^{(t_u)})^2}{N_{k,j}^{(u)}  } } - \frac{1}{N^{(r)}}\frac{\delta_i^{(u)}}{\sum_{j=1}^D\frac{(\hat{\sigma}^{(u)}_{b,j})^2 (\pi_{j}^{(t_u)})^2}{N_{b,j}^{(r)}} +\sum_{j=1}^D \frac{(\hat{\sigma}^{(u)}_{i,\jc})^2 (\pi_{\jc}^{(t_u)})^2}{N_{i,\jc}^{(r)}} } \notag\\
    &+\Opar{\logrr}  \notag\\
    & \le\Opar{\sqrt{\frac{\log\log r}{r}}}. \notag
\end{align}
The second inequality holds because $\frac{(\hat{\sigma}^{(u)}_{i,j})^2(\pi_{j}^{(t_u)})^2}{ N_{i,j}^{(u)}/N^{(r)} } = \Opar{\logrr}$ for all $i\in\mathcal{I},j\neq \jc$.
This finalizes the proof. \hfill $\blacksquare$

\newcommand{\pr}{{(r)}}
\newcommand{\ptr}{{(t_{\ell_r})}}
\subsection{Proof of Lemma \ref{lem: consistency continuous}.}
\textbf{Proof of 1.} 
By Assumption \ref{assump: consistency continuous space}.2 and \ref{assump: consistency continuous space}.3, we know the Bayesian consistency holds by Doob's consistency theorem. Furthermore by Assumption \ref{assump: consistency continuous space}.1, we know there exists $\epsilon>0$, such that $B_\epsilon(\theta^c) \in \Theta_{j^c}$, where $B_\epsilon(\theta^c)$ is the ball centered at $\theta^c$ with radius $\epsilon$.
Hence, by Bayesian consistency, $\pi^t_\jc \rightarrow 1$ almost surely for every sequence of input data.

\textbf{Proof of 2.}
\textbf{Mean:}
Recall $X_{i,\jc}^\pr$ is the $r$th simulation output for $(i,\jc)$. Let  $\ell_r$ be the iteration at which $X_{i,\jc}^\pr$ is simulated. Here $\ell_r$ is random given $r$. 

 Let $\Bar{y}_{i,\jc}^t = \mathbb{E}_{\Tilde{\pi}^t_{\jc}}[y_{i,\jc}^{(1)}(\theta)]$, $Z_{i,\jc}^\pr = X_{i,\jc}^\pr - \Bar{y}_{i,\jc}^t$, 
 Notice
 $$ \hat{\mu}^\pl_{i,\jc} - \mu_i(\theta^c) = \underbrace{ \frac{1}{N_{i,\jc}^\pl}\sum_{r=1}^{N_{i,\jc}^\pl} Z_{i,\jc}^\pr}_{(I)} +
 \underbrace{\frac{1}{N_{i,\jc}^\pl}\sum_{r=1}^{N_{i,\jc}^\pl} (\Bar{y}_{i,\jc}^\ptr - \mu_i(\theta^c))}_{(II)}. $$
 
 For (I), we have 
 \begin{align*}
     \mathbb{E} [ Z_{i,\jc}^\pr | Z_{i,\jc}^{(1)},\ldots,Z_{i,\jc}^{(r-1)}] =& \mathbb{E} [\mathbb{E} [ Z_{i,\jc}^\pr | \Tilde{\pi}_\jc^{(t_{\ell_r})},\ell_r,Z_{i,\jc}^{(1)},\ldots,Z_{i,\jc}^{(r-1)}]|Z_{i,\jc}^{(1)},\ldots,Z_{i,\jc}^{(r-1)}].
 \end{align*}
 Notice $\ell_r = \ell,Z_{i,jc}^{(r')},r'=1,\ldots,r-1$ are determined by input data and simulation outputs up to iteration $\ell$.
By Assumption \ref{assump: consistency continuous space}.5, 
$$ \mathbb{E} [ Z_{i,\jc}^\pr | \Tilde{\pi}_\jc^{(t_{\ell_r})},\ell_r,Z_{i,\jc}^{(1)},\ldots,Z_{i,\jc}^{(r-1)}] = \mathbb{E} [ Z_{i,\jc}^\pr | \Tilde{\pi}_\jc^{(t_{\ell_r})}] = 0.$$
This implies $ \mathbb{E} [ Z_{i,\jc}^\pr | Z_{i,\jc}^{(1)},\ldots,Z_{i,\jc}^{(r-1)}] =0, $ which means $\{Z_{i,\jc}^\pr\}_{r=1}^\infty$ is a Martingale difference sequence (MDS). 
Moreover, let $\theta^{(\ell_r)}|\ell_r \sim \Tilde{\pi}^{t_{\ell_r}}_\jc$ be the input parameter under which $X_{i,\jc}^\pr$ is simulated. Here $\Tilde{\pi}^t_j$ is the conditional posterior distribution at stage $t$.

$$\frac{1}{2}\mathbb{E} [(Z_{i,\jc}^\pr)^2] \le \mathbb{E}[(X_{i,\jc}^{(r)})^2+(\Bar{y}_{i,\jc}^{t_{\ell_r}})^2] \le \mathbb{E} [\mathbb{E} [(X_{i,\jc}^{(r)})^2 |\theta^{(\ell_r)}]] + \mathbb{E}[(\Bar{y}_{i,\jc}^{t_{\ell_r}})^2] \le \max_\theta \{ (y_{i,\jc}^{(1)}(\theta))^2+y_{i,\jc}^{(2)}(\theta)\} < \infty, $$
where the boundedness comes from the fact that $y_{i,\jc}^{(k)},k=1,2$ is continuous and the parameter set is compact by Assumption \ref{assump: consistency continuous space}. Hence, 
by strong law of large number (SLLN) for MDS \citep{csorgHo1968strong}, we obtain $\lim_{\ell \rightarrow 
 \infty} \frac{1}{\ell}\sum_{r=1}^\ell Z_{i,\jc}^\pr = 0$ almost surely. Notice
Since $N_{i,\jc}^\pl \rightarrow \infty$ almost surely as $\ell \rightarrow \infty$, we have $\lim_{\ell \rightarrow 
 \infty} \frac{1}{N_{i,\jc}^\pl}\sum_{r=1}^{N_{i,\jc}^\pl} Z_{i,\jc}^\pr = 0$ almost surely. 
 


For (II), by Bayesian consistency, $\pi^t$ converges to a Dirac measure at $\theta^c$ almost surely. Further notice $y_{i,j}^{(1)}$ is continuous and bounded. This implies 
$ \Bar{y}_{i,\jc}^{(t)} \rightarrow \mu_i(\theta^c)$ almost surely. Since $N_{i,\jc}^{(\ell)} \rightarrow \infty$ almost surely, we have the average sum $(III)\rightarrow 0$ almost surely. (That is, the convergence of sequence $\{a_n\}_{n=1}^\infty$ implies the convergence of sequence $\{S_n\}_{n=1}^\infty$ to the same limit, where $S_n = \frac{1}{n}\sum_{k=1}^n a_k$ ).

 \textbf{Variance:}
 For $\hat{\sigma}_{i,\jc}^\pl$, let $\Tilde{y}_{i,\jc}^t = \mathbb{E}_{\Tilde{\pi}_\jc^t}[y_{i,\jc}^{(2)}(\theta)]$.

\begin{align*}
    \frac{N_{i,\jc}^\pl}{N_{i,\jc}^\pl-1}(\hat{\sigma}_{i,\jc}^\pl)^2 - \sigma_i(\theta^c)^2 = & \frac{1}{N_{i,\jc}^\pl} \sum_{r=1}^{N_{i,\jc}^\pl} (X_{i,\jc}^\pr-\hat{\mu}_{i,\jc}^\pl)^2 -\sigma_i(\theta^c)^2 \\
    =& \underbrace{ \frac{ 1}{N_{i,\jc}^\pl} \sum_{r=1}^{N_{i,\jc}^\pl} 
    \left((X_{i,\jc}^{(r)})^2 - \Tilde{y}_{i,\jc}^{(t_{\ell_r})} \right)}_{(I)} + \underbrace{\frac{1}{N_{i,\jc}^{(\ell)}} \sum_{r=1}^{N_{i,\jc}^{(\ell)}}\left(\Tilde{y}_{i,\jc}^{(t_{\ell_r})} -(\hat{\mu}^\pl_{i,\jc})^2 -\sigma_i(\theta^c)^2\right)}_{(II)} 
\end{align*}
(I) can be proved to converge to $0$ in a similar manner as for the mean estimator, by regarding $(X_{i,\jc}^{(r)})^2, r=1,\ldots,N_{i,\jc}^\pl$ as the simulation samples and 
$ \mathbb{E}\left[(X_{i,\jc}^{(r)})^2 | \Tilde{\pi}_\jc^{t_{\ell_r}}\right] = \mathbb{E}\left[ \mathbb{E}\left[(X_{i,\jc}^{(r)})^2 | \theta^{(\ell_r)} \right] | \Tilde{\pi}_\jc^{t_{\ell_r}}\right] = \mathbb{E}\left[ y_{i,\jc}^{(2)}|\Tilde{\pi}_\jc^{t_{\ell_r}}\right] = \Tilde{y}_{i,\jc}^{t_{\ell_r}}$. Hence, $Z_{i,\jc}^{(r)} := (X_{i,\jc}^{(r)})^2 - \Tilde{y}_{i,\jc}^{(t_{\ell_r})} $ is a MDS sequence for $r=1,2,\ldots$. Furthermore we can show $\mathbb{E}[Z_{i,\jc}^2] \le 2 \max_{\theta \in \Theta}\{ y_{i,\jc}^{(4)} + (y_{i,\jc}^{(2)})^2\} < \infty $. Hence, we can invoke the SLLN for MDS to prove $(I)\rightarrow0$ almost surely. 

For (II), again by the Bayesian consistency and the continuity of $y_{i,\jc}^{(2)}(\theta)$, we can show almost surely, $\Tilde{y}_{i,\jc}^{t_{\ell_r}} \rightarrow \mu_i(\theta^c)^2 + \sigma_i(\theta^c)^2$. Since we also have $\hat{\mu}_{i,\jc}^\pl \rightarrow \mu_i(\theta^c)$ almost surely, we obtain $\Tilde{y}_{i,\jc}^{t_{\ell_r}}  -  \hat{\mu}_{i,\jc}^\pl - \sigma_i(\theta^c)^2 \rightarrow 0$ almost surely. This implies $(II) \rightarrow 0$ almost surely. The proof is complete.
 $\hfill \blacksquare$

\end{document}